\documentclass[letterpaper, 12pt, oneside]{book}

    \usepackage[left=72pt, right=72pt, top=72pt, bottom=72pt, dvips, pdftex]{geometry}
    \usepackage{xcolor}
    \usepackage[procnames]{listings}
    \usepackage{color}
    \usepackage{fancyhdr} 
    \usepackage{amsmath,amsthm, amsfonts,amssymb}
    \usepackage{setspace} 
\usepackage{multirow}
    \usepackage[linktocpage,bookmarksopen,bookmarksnumbered,
                pdftitle={Knots and Links from Random Projection}]{hyperref}
\usepackage{dsfont}
\usepackage{float}
     \usepackage{graphicx}
\usepackage{url}
\usepackage[toc,page]{appendix}
   \newtheorem{Theorem}{Theorem}[section]
    
    \newtheorem{Lemma}[Theorem]{Lemma}
    
    \theoremstyle{definition}
        \newtheorem{Definition}[Theorem]{Definition}

    \theoremstyle{remark}
        \newtheorem{Remark}[Theorem]{Remark}
    \theoremstyle{Acknowledgments}

    \newenvironment{Proof}{\par\noindent{\sc Proof}\quad}{\hfill\qed\par\smallskip}
    
    \newenvironment{RestateTheorem}[2]{\par\vspace{12pt}\noindent{\bf Theorem \ref{#2}.\it#1}\it}{\par\vspace{2pt}}

    \numberwithin{equation}{section}
    \numberwithin{figure}{section}
    \newcommand{\n}{\vspace{12pt}} 
   
    \newcommand{\va}{\textbf{a}}
\newcommand{\vu}{\textbf{u}}
\newcommand{\vU}{\textbf{U}}
\newcommand{\vq}{\textbf{q}}
\newcommand{\ve}{\textbf{e}}
\newcommand{\vV}{\textbf{v}}
\newcommand{\vb}{\textbf{b}}
\newcommand{\vr}{\textbf{r}}
\newcommand{\vx}{\textbf{x}}
\newcommand{\vy}{\textbf{y}}
\newcommand{\vz}{\textbf{z}}
\newcommand{\proj}{\textrm{proj}}  
    \newcommand{\newchapter}[3] 
	{                           
        \chapter[#2]{#3}
        \chaptermark{#1}
        \thispagestyle{fancy}
	}

\usepackage[numbers]{natbib}

\begin{document}
    
\definecolor{b1}{RGB}{0,0,0}
\lstset{language=Python, 
        basicstyle=\ttfamily\small, 
        keywordstyle=\color{b1},
        commentstyle=\color{b1},
        stringstyle=\color{b1},
        showstringspaces=false,
        identifierstyle=\color{b1},
        procnamekeys={def,class}}
    \pagenumbering{roman}
    \pagestyle{plain}

    \singlespacing

    ~\vspace{-0.75in} 
    \begin{center}

        \begin{huge}
           Knots and Links from Random Projections
        \end{huge}\\\n
        By\\\n
        {\sc Christopher David Westenberger\\ \href{mailto:cdwestenberger@gmail.com}{cdwestenberger@gmail.com}}\\
    \end{center}
 
    \centerline{\textbf{\underline{Abstract}}}
In this paper we study a model of random knots obtained by fixing a space curve in $n$-dimensional Euclidean space with $n>3$, and orthogonally projecting the space curve on to random $3$ dimensional subspaces. By varying the space curve we obtain different models of random parametrized knots, and we will study how the expectation value of the curvature changes as a function of the initial parametrized space curve. In the case when the initial data is a pair of space curves, or more generally a pair of manifolds satisfying certain conditions on their dimension, then we obtain models of random links for which we will give methods to compute the second moment of the linking number. As an application of our computations, we will study numerous models of random knots and links, and how to recover these models by appropriately choosing the initial space curves to be projected.
    \tableofcontents
    \pagestyle{fancy}
    \pagenumbering{arabic}
       
    %
    %

    \newchapter{Introduction}{Introduction}{Introduction}
    \label{sec:Intro}
    
        \section[Informal Overview and Historical Motivation]{Informal Overview and Historical Motivation}
        \label{sec:Intro:Motivation}

Over the past few decades numerous models of random knots and links have been studied and were developed either with a specific application in mind, or with the hope that the model will be sufficiently universal in the sense that it avoids biasing certain knots and links. Models include closed random walks on a lattice, random knots with confinement, random equilateral polygons, random kinematical links, random knots from billiard diagrams, random Fourier knots (for the aforementioned see \cite{Alvarado},\cite{Arsuaga}, \cite{Cohen},\cite{Duplantier}, for example), diagrams sampled from random 4-valent graphs (by Dunfield et al. See SnapPy documentation: \cite{SnapPy}), and the recent Petaluma model studied in \cite{Hass}, \cite{Hass2}. The distributions of the invariants depend on the details of the models, and the difficulty in computing certain statistics resides in the sampling method. Given this vast list, the question of whether or not there is a universal model of random knotting and linking that may contain some of the aforementioned models as special cases remains both ill-posed and unexplored. In this paper we seek to explore this question by introducing a further model of random knots which is sufficiently general that it may be able to encompass some previously studied models of random knots and links.

            \subsection[The Model]{The Model}
            \label{sec:Intro:Motivation:Part1}
        
In this paper we will consider a model of random knots and links given by first fixing a closed space curve (or pair of closed space curves in the case of links) $\vr(t):\mathbb{S}^1\rightarrow\mathbb{R}^n$ with $n>3$, and then choosing a $3$ dimensional subspace in $\mathbb{R}^n$ at random to project the curve $\vr(t)$ on to. The image under the projection of $\vr(t)$ is a parametrized knot in $\mathbb{R}^3$ for which one may compute a number of quantities associated to the space curve, like the curvature and the writhe. Our hope is to compute the mean and variance of quantities like the curvature, writhe, and linking number with respect to the unique normalized $O(n)$ invariant measure on the Grassmannian of $3$-dimensional subspaces in $\mathbb{R}^n$. The rotational invariance of the measure along with the numerous scaling and multilinearity properties of the aforementioned quantities will allow us to greatly simplify the calculations involved. By choosing different initial curves $\vr(t)$ in $\mathbb{R}^n$, we will obtain different models of knots and links, and moreover we will expect that the distribution of invariants associated with such models to have a crucial dependence on $\vr(t)$.

            \subsection[Configuration Space Integrals and Invariants]{Configuration Space Integrals and Invariants}
            \label{sec:Intro:Motivation:Part2}

Since the random knots and links that will be generated with the model described in the previous section will be parametrized, then the sort of quantities that we will consider will be ones that have descriptions as configuration space integrals. In the case of links, one such invariant is the linking number, and when the components of the link are given by non-intersecting differentiable curves, $\gamma_1$ and $\gamma_2$, given by the parametrizations $\gamma_1=\vr_1(t)$ and $\gamma_2=\vr_2(s)$, then the linking number may be computed using the Gauss linking integral:
\begin{equation}
\mathrm{Lk}(\gamma_1,\gamma_2)=\frac{1}{4\pi}\int\displaylimits_{(s,t)\in\mathbb{T}^2}\frac{(\dot{\vr}_1(t)\times\dot{\vr}_2(s))\cdot(\vr_2(s)-\vr_1(t))}{\|\vr_2(s)-\vr_1(t)\|^3}dA
\end{equation}
where $\mathbb{T}^2=\mathbb{S}^1\times\mathbb{S}^1$ is the torus where each factor is the interval $\big[0,2\pi\big]$ with its endpoints identified, and $dA=dsdt$.
As motivation, we will use the model described in the previous section to find the expected value of a closely related quantity, the average inter-crossing number, which we will denote $\mathrm{ICN}$ (see \cite{diao} for further details):
\begin{equation}\nonumber
\mathrm{ICN}(\gamma_1,\gamma_2)=\frac{1}{4\pi}\int\displaylimits_{(s,t)\in\mathbb{T}^2}\frac{|(\dot{\vr}_1(t)\times\dot{\vr}_2(s))\cdot(\vr_2(s)-\vr_1(t))|}{\|\vr_2(s)-\vr_1(t)\|^3}dA
\end{equation}
To do so, let $\gamma_1=\vr_1(t)$ and $\gamma_2=\vr_2(s)$ be differentiable closed space curves in $\mathbb{R}^4$ and define $\vV=(v_1,v_2,v_3,v_4)$ and $d\vV=exp(-\frac{1}{2}(v_1^2+v_2^2+v_3^2+v_4^2))dv_1dv_2dv_3dv_4$. This way, we have that:
\begin{displaymath}
\langle\mathrm{ICN}\rangle=\frac{1}{4\pi(2\pi)^{2}}\int\displaylimits_{\vV\in\mathbb{R}^4}\int
\displaylimits_{\mathbb{T}^2}\frac{|\mathrm{Det}\big[\dot{\vr}_1(t),\dot{\vr}_2(s),\vr_2(s)-\vr_1(t),\frac{\vV}{\|\vV\|}\big]|}{\|\proj_{\vV^\perp}(\vr_2(s)-\vr_1(t))\|^3}dsdtd\vV
\end{displaymath}
where $\proj_{\vV^\perp}(\vr_2(s)-\vr_1(t))$ denotes the orthogonal projection on to the orthogonal complement of $\frac{\vV}{\|\vV\|}$. In order to emphasize the dependence on the initial data $\vr_1(t)$ and $\vr_2(s)$, we will reverse the order of integration:
\begin{displaymath}\label{expectedLinking}
\langle\mathrm{ICN}\rangle=\frac{1}{4\pi(2\pi)^{2}}\int
\displaylimits_{\mathbb{T}^2}\int\displaylimits_{\vV\in\mathbb{R}^4}\frac{|\mathrm{Det}\big[\dot{\vr}_1(t),\dot{\vr}_2(s),\vr_2(s)-\vr_1(t),\frac{\vV}{\|\vV\|}\big]|}{\|\proj_{\vV^\perp}(\vr_2(s)-\vr_1(t))\|^3}d\vV dA
\end{displaymath}
Next, given vectors $\va_1,\va_2,\va_3\in\mathbb{R}^n$ define a matrix $A$ by making the vectors $\va_i$ the columns, which from now on we will denote by $A=\big[\va_1,\va_2,\va_3\big]$, and then define the function:
\begin{equation}\label{Iint}
I_{\langle\mathrm{ICN}\rangle}(A)=\frac{1}{(2\pi)^2}\int\displaylimits_{\vV\in\mathbb{R}^4}\frac{|\mathrm{Det}\big[\va_1,\va_2,\va_3,\frac{\vV}{\|\vV\|}\big]|}{\|\proj_{\vV^\perp}(\va_3)\|^3}d\vV
\end{equation} where $\big[\va_1,\va_2,\va_3,\frac{\vV}{\|\vV\|}\big]\in\mathrm{Mat_{4,4}}$.
Given the definition in \eqref{Iint}, we may simplify the formula for $\langle\mathrm{ICN}\rangle$ as:
\begin{equation}\nonumber
\langle\mathrm{ICN}\rangle=\frac{1}{4\pi}\int\displaylimits_{\mathbb{T}^2}I_{\langle\mathrm{ICN}\rangle}(\big[\dot{\vr}_1(t),\dot{\vr}_2(s),\vr_2(s)-\vr_1(t)\big])dsdt
\end{equation}
It will be shown in the following sections that $I_{\langle\mathrm{ICN}\rangle}(\big[\va_1,\va_2,\va_3\big])$ has an especially nice form, and is given by:
\begin{eqnarray}\nonumber
I_{\langle\mathrm{ICN}\rangle}(\big[\va_1,\va_2,\va_3\big])=\nonumber\\
\frac{\sqrt{\mathrm{Det}\left(\big[\va_1,\va_2,\va_3\big]^T\big[\va_1,\va_2,\va_3\big]\right)}}{\|\va_3\|^3},
\end{eqnarray}
so that:
\begin{align}\label{finalC}
\langle\mathrm{ICN}\rangle&=\frac{1}{4\pi}\int\displaylimits_{\mathbb{T}^2}I_{\langle\mathrm{ICN}\rangle}(\big[\dot{\vr}_1(t),\dot{\vr}_2(s),\vr_2(s)-\vr_1(t)\big])\nonumber\\
&=\frac{1}{4\pi}\int\displaylimits_{\mathbb{T}^2}\frac{\sqrt{\mathrm{Det}(\big[\dot{\vr}_1(t),\dot{\vr}_2(s),\vr_2(s)-\vr_1(t)\big]^T\big[\dot{\vr}_1(t),\dot{\vr}_2(s),\vr_2(s)-\vr_1(t)\big])}}{\|\vr_2(s)-\vr_1(t)\|^3}dsdt
\end{align}
For the most part, the goal of this work will be to both find closed forms and understand the analytic properties of functions like $I_{\langle\mathrm{ICN}\rangle}(\big[\dot{\vr}_1(t),\dot{\vr}_2(s),\vr_2(s)-\vr_1(t)\big])$. By performing this procedure of reversing the order of integration, and then integrating over all projections, we may then find bounds on $\langle\mathrm{ICN}\rangle$ in terms of the parametrization of the initial space curves $\vr_1(t)$ and $\vr_2(s)$. The computation leading to \eqref{finalC} will be completed in the second chapter, along with a similar analysis for $\langle\kappa(C)\rangle$, where $\kappa(C)$ is the total curvature of a knot (the total curvature case is not a new result, and was already discovered in a classic paper by \cite{Fary}, \cite{Milnor}, and can be found in recent papers such as \cite{Sullivan}). 

\begin{figure}[H]\label{links}
\centering
\includegraphics[scale=0.25]{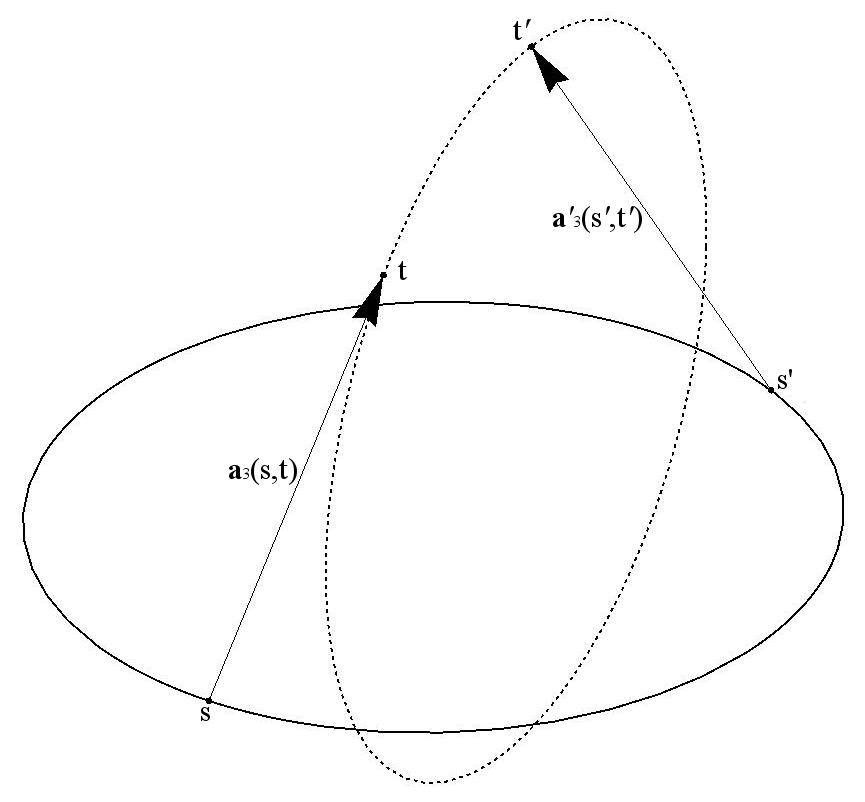}
\caption{The geometric configuration for computing $\mathrm{Lk}^2$ for the closed space curves $\vr_1(s)$ and $\vr_2(t)$ where $\va_3(s,t)=\vr_2(s)-\vr_1(t)$ and $\va_3(s',t')=\vr_2(s')-\vr_1(t')$}
\end{figure}

Since $\kappa(C)$ and $\langle\mathrm{ICN}\rangle$ are not invariants, but are only bounds on invariants of knots and links, then it will be a bit more interesting to consider instead $\langle\mathrm{Lk}^2\rangle$ and the corresponding function $I_{\langle\mathrm{Lk}^2\rangle}(A,A')$ obtained from reversing the order of integration as was demonstrated above. Given two space curves $\gamma_1$ and $\gamma_2$ with parametrizations $\vr_1(t)$ and $\vr_2(s)$, then we may compute the value of $\mathrm{Lk}^2(\gamma_1,\gamma_2)$ from the Gauss linking integral by computing:
\begin{align*}\nonumber
&\mathrm{Lk}^2(\gamma_1,\gamma_2)=\\
&\frac{1}{16\pi^2}\int\displaylimits_{(s,t)\in\mathbb{T}^2}\int\displaylimits_{(s',t')\in\mathbb{T}^2}\frac{\mathrm{Det}\big[\dot{\vr}_1(t),\dot{\vr}_2(s),\vr_2(s)-\vr_1(t)\big]\mathrm{Det}\big[\dot{\vr}_1(t'),\dot{\vr}_2(s'),\vr_2(s')-\vr_1(t')\big]}{\|\vr_2(s)-\vr_1(t)\|^3\|\vr_2(s')-\vr_1(t')\|^3}dsdtds'dt'
\end{align*}
Figure \ref{links} depicts the configurations of pairs of points that we will typically want to integrate over in order to find the value of $\mathrm{Lk}^2(\gamma_1,\gamma_2)$.
Obtaining $I_{\langle\mathrm{Lk}^2\rangle}(A,A')$ in this case is considerably more difficult, and will occupy the third chapter, which is the main part of this work. The computation involved after switching the order of integration is summarized with the following lemma:
\begin{Lemma}\label{4case}
Given matrices $A=\big[\va_1,\va_2,\va_3\big],A'=\big[\va_1',\va_2',\va_3'\big]\in\mathrm{Mat}_{4,3}$ and the function
\begin{equation}\nonumber
I_{\langle\mathrm{Lk}^2\rangle}(A,A')=\frac{1}{(2\pi)^2}\int_{\mathbb{R}^4}\frac{\mathrm{Det}([A\:\frac{\vV}{\|\vV\|}])\mathrm{Det}([A'\:\frac{\vV}{\|\vV\|}])}{\|\proj_{\vV^\perp}(\va_3)\|^3\|\proj_{\vV^\perp}(\va'_3)\|^3}e^{-\|\vV\|^2/2}d\vV
\end{equation}
where $\big[A\:\frac{\vV}{\|\vV\|}\big]$ denotes a new matrix in $\mathrm{Mat}_{4,4}$ with the column $\frac{\vV}{\|\vV\|}$ appended to the matrix $A$, then we have
\begin{center}
\begin{align*}\nonumber
&\frac{\pi}{2}I_{\langle\mathrm{Lk}^2\rangle}(A,A')=\\
&\frac{(\|\va_3\|^2\|\va_3'\|^2-(\va_3\cdot\va_3')^2)\mathrm{Det}(A^TA')+(\va_3\cdot\va_3')\mathrm{Det}(\big[\va_3,\va'_1,\va_2',\va_3'\big])\mathrm{Det}(\big[\va_3',\va_1,\va_2,\va_3\big])}{\|\va_3\|^2\|\va_3'\|^2(\|\va_3\|^2\|\va_3'\|^2-(\va_3\cdot\va_3')^2)^{3/2}}\nonumber
\end{align*}  
\end{center}
\end{Lemma}
It will be seen throughout the course of the proof of the above lemma that the requirement that the initial space curves $\vr_1(t)$ and $\vr_2(s)$ be embedded in $\mathbb{R}^4$ serves to simplify the integration over $Gr(4,3)$ since the subspaces we are integrating over are codimension $1$ so that $Gr(4,3)\cong\mathbb{S}^3$. This restrictive assumption will be removed in the second half of chapter $3$ by introducing the Stiefel manifold of orthonormal $3$-frames in $\mathbb{R}^n$, denoted $V_{n,3}$, and using a decomposition of its unique normalized $O(n)$ measure coupled with the previous result in order to find $\langle\mathrm{Lk}^2\rangle$ for any pair of curves in $\mathbb{R}^n$, for $n$ sufficiently large.

\section{Summary of the Main Results}
\begin{Theorem}\label{expInterlinking}
Given two closed, differentiable, non-intersecting space curves $\gamma_1,\gamma_2:\mathbb{S}^1\rightarrow\mathbb{R}^n$ with parametrizations $\vr_1(t)$ and $\vr_2(s)$, and let $\mathrm{proj}_U(\gamma_1),\mathrm{proj}_U(\gamma_2)$ denote the orthogonal projections of $\gamma_1$ and $\gamma_2$ on to the $3$ dimensional subspace spanned by the columns of $U$, then the expected value of $\mathrm{ICN}(\mathrm{proj}_U(\gamma_1),\mathrm{proj}_U(\gamma_2))$, denoted $\langle\mathrm{ICN}\rangle$, averaged over all orthogonal projections to $3$ dimensional subspaces with respect to the unique normalized $O(n)$-invariant measure on $Gr(n,3)$, is given by the following integral:
\begin{equation}\nonumber
\langle\mathrm{ICN}\rangle=\frac{1}{4\pi}\int\displaylimits_{(s,t)\in\mathbb{T}^2}I_{\langle\mathrm{ICN}\rangle}(s,t)dsdt
\end{equation}
where
\begin{equation}\nonumber
I_{\langle\mathrm{ICN}\rangle}(s,t)=C_{\langle\mathrm{ICN}\rangle}\frac{\sqrt{\textrm{Det}(\big[\frac{d\vr_1(t)}{dt},\frac{d\vr_2(s)}{ds},\vr_2(s)-\vr_1(t)\big]^T\big[\frac{d\vr_1(t)}{dt},\frac{d\vr_2(s)}{ds},\vr_2(s)-\vr_1(t)\big])}}{\|\vr_2(s)-\vr_1(t)\|^3},
\end{equation}
where $C_{\langle\mathrm{ICN}\rangle}$ is a constant.
\end{Theorem}
After proving the above theorem we then focus on the computation of $\langle\mathrm{ICN}\rangle$:
\begin{Theorem}\label{GaussianConfigLk}
Let $\vr_1(t):\mathbb{S}^1\rightarrow\mathbb{R}^{2n+1}$ and $\vr_2(s):\mathbb{S}^1\rightarrow\mathbb{R}^{2n+1}$ be given by:
\begin{align*}\nonumber
\vr_1(t)&=c_0\textbf{e}_0+\sum_{k=1}^{n}(c_k\cos(kt)\ve_{2k-1}+c_k\sin(kt)\ve_{2k})\\
&=(c_0,c_1\cos(t),c_1\sin(t),c_2\cos(2t),c_2\sin(2t),...,c_n\cos(nt),c_n\sin(nt))\textrm{, and}\\
\vr_2(s)&=d_0\textbf{e}_0+\sum_{k=1}^{n}(d_k\cos(ks)\ve_{2k-1}+d_k\sin(ks)\ve_{2k})\\
&=(d_0,d_1\cos(s),d_1\sin(s),d_2\cos(2s),d_2\sin(2s),...,d_n\cos(ns),d_n\sin(ns))\nonumber,
\end{align*}
where $\{\ve_i\}_{i=1}^n$ are the standard basis vectors, then $\langle\mathrm{ICN}\rangle$ is finite and satisfies the following bound:
\begin{equation}\nonumber
\langle\mathrm{ICN}\rangle\leq C_{\langle\mathrm{ICN}\rangle}\frac{\sqrt{(\sum_{j=0}^n j^2c_j^2)(\sum_{j=0}^n j^2d_j^2)}}{\mathrm{min}_{(s,t)\in\mathbb{T}^2}\|\vr_2(s)-\vr_1(t)\|^2}
\end{equation}
\end{Theorem}
In chapter $3$ we move on to proving the main theorem of the paper concerning the second moment of the linking number.
\begin{Theorem}\label{mainResultTheorem}
Given two closed, differentiable, non-intersecting space curves $\gamma_1,\gamma_2:\mathbb{S}^1\rightarrow\mathbb{R}^n$ with parametrizations $\vr_1(t)$ and $\vr_2(s)$, and let $\mathrm{proj}_U(\gamma_1),\mathrm{proj}_U(\gamma_2)$ denote the orthogonal projections of $\gamma_1$ and $\gamma_2$ on to the $3$ dimensional subspace spanned by the columns of $U$, then the expected value of $\mathrm{Lk}^2(\mathrm{proj}_U(\gamma_1),\mathrm{proj}_U(\gamma_2))$, denoted $\langle\mathrm{Lk}^2\rangle$, averaged over all orthogonal projections to $3$ dimensional subspaces with respect to the unique normalized $O(n)$-invariant measure on $Gr(n,3)$, is given by the following integral:
\begin{equation}\nonumber
\langle\mathrm{Lk}^2\rangle=\frac{1}{16\pi^2}\int\displaylimits_{(s,t)\in\mathbb{T}^2}\int\displaylimits_
{(s',t')\in\mathbb{T}^2}I_{\langle\mathrm{Lk}^2\rangle}(A(s,t),A'(s',t'))dsdtds'dt'
\end{equation}
where
\begin{align*}\label{eq:result3THM}
&\frac{\pi}{4}I_{\langle\mathrm{Lk}^2\rangle}(A,A')=\nonumber\\
&\frac{(\va_3\cdot\va_3')-\sqrt{\|\va_3\|^2\|\va_3'\|^2-(\va_3\cdot\va_3')^2}\sin^{-1}\left(\frac{\va_3\cdot\va_3'}{\|\va_3\|\|\va_3'\|}\right)}{(\va_3\cdot\va_3')^3\sqrt{\|\va_3\|^2\|\va_3'\|^2-(\va_3\cdot\va_3')^2}}\mathrm{Det}\left(\begin{bmatrix}
    \va_1\cdot\va_1'       & \va_1\cdot\va_2' & \va_1\cdot\va_3'   \\
    \va_2\cdot\va_1'       & \va_2\cdot\va_2' & \va_2\cdot\va_3'   \\
    \va_3\cdot\va_1'       & \va_3\cdot\va_2' & \va_3\cdot\va_3'   \\
\end{bmatrix}\right)
\end{align*} 
and
\begin{align*}\nonumber
\va_1&=\dot{\vr}_1(t)\mathrm{, }\va_2=\dot{\vr}_2(s)\mathrm{, }\va_3=\vr_2(s)-\vr_1(t)\\
\va_1'&=\dot{\vr}_1(t')\mathrm{, }\va_2'=\dot{\vr}_2(s')\mathrm{, }\va_3'=\vr_2(s')-\vr_1(t')\\
A&=\big[\dot{\vr}_1(t),\dot{\vr}_2(s),\vr_2(s)-\vr_1(t)\big]\\
A'&=\big[\dot{\vr}_1(t'),\dot{\vr}_2(s'),\vr_2(s')-\vr_1(t')\big]\\
\end{align*}
\end{Theorem}
After integrating over all the projections, we then turn our attention to integrating over the configuration space. 
\begin{Theorem}\label{Lk2Bound}
With the definitions above, let $v_1=\textrm{max}_{t\in\mathbb{S}^1}\|\dot{\vr}_1(t)\|$, $v_2=\textrm{max}_{s\in\mathbb{S}^1}\|\dot{\vr}_2(s)\|$, $k=\mathrm{min}_{(s,t)\in\mathbb{T}^2}\|\vr_2(t)-\vr_1(s)\|$ and 
\begin{equation}\nonumber
C=\int_{0}^{2\pi}\int_{0}^{2\pi}\int_{0}^{2\pi}\int_{0}^{2\pi}\frac{1}{\sqrt{\|\va_3(s,t)\|^2\|\va_3'(s',t')\|^2-(\va_3(s,t)\cdot\va_3'(s',t'))^2}}dsdtds'dt',
\end{equation}
then
\begin{equation}\nonumber
\langle\mathrm{Lk}^2\rangle\leq\frac{1}{(4\pi)^2}\frac{4Cv_1^2v_2^2}{\pi k^2}
\end{equation}
\end{Theorem}
That is, we obtain a bound on the second moment of the linking number provided the initial space curves are chosen so that $v_1,v_2$ and $C$ are finite. In the fourth chapter we will prove theorems very analogous to theorems in chapter 3 concerning the second moment of linking numbers associated to manifolds as defined in \cite{higherlink}. Lastly, in order to show the applicability of some of our results, in the penultimate chapter we consider a model that is very similar to the model in \cite{Hass}.
        \label{sec:Intro:Summary}

    %

    \newchapter{Expectation Values}{Expectation Values}{Expectation Values}
    \label{sec:LabelForChapter2}
    
        \section[Background]{Background}
        \label{sec:LabelForChapter2:Section1} In the first part of this chapter we will find the expectation value of $\langle\mathrm{ICN}\rangle$ when we fix a pair of space curves in $\mathbb{R}^4$ and project the curves randomly on to $3$ dimensional subspaces. As discussed in the introduction, this will require two steps, the first will be integrating over all $3$ dimensional subspaces, and then integrating that result over the configuration space of two points on distinct components of a link. The first step can be done exactly, whereas the second will be approached by making specific assumptions about the initial data. In the second part of this chapter we will find the expected curvature when we fix one space curve in $\mathbb{R}^n$, and project the curve randomly on to $3$ dimensional subspaces, where we pick random subspaces by sampling the spans of the columns of Gaussian random matrices in $\mathrm{Mat}_{n,3}$.

        \section[Integration Over $3$-Dimensional Subspaces]{Integration Over $3$-Dimensional Subspaces}
        \label{sec:LabelForChapter2:Section2} 
            
	\subsection{Average Intercrossing Number: $\langle\mathrm{ICN}\rangle$}
We start with the following lemma:
\begin{Lemma}\label{firstlemma}
Given a matrix $A=\big[\va_1,\va_2,\va_3\big]\in\mathrm{Mat}_{4,3}$ and the function
\begin{equation}\nonumber
I_{\langle\mathrm{ICN}\rangle}(A)=\frac{1}{(2\pi)^2}\int_{\mathbb{R}^4}\frac{|\mathrm{Det}([A\:\frac{\vV}{\|\vV\|}])|}{\|\proj_{\vV^\perp}(\va_3)\|^3}e^{-\|\vV\|^2/2}d\vV
\end{equation}
where $\big[A\:\frac{\vV}{\|\vV\|}\big]$ denotes a new matrix in $\mathrm{Mat}_{4,4}$ with the column $\frac{\vV}{\|\vV\|}$ appended to the matrix $A$, then we have
\begin{equation}\nonumber
I_{\langle\mathrm{ICN}\rangle}(\big[\va_1,\va_2,\va_3\big])=\frac{\sqrt{\mathrm{Det}\left(\big[\va_1,\va_2,\va_3\big]^T\big[\va_1,\va_2,\va_3\big]\right)}}{\|\va_3\|^3},
\end{equation}
\end{Lemma}
\begin{Proof}\label{prevlemma}
We first observe the following properties of $I_{\langle\mathrm{ICN}\rangle}(\big[\va_1,\va_2,\va_3\big])$:
\begin{align}\label{props}
I_{\langle\mathrm{ICN}\rangle}\left(A\left( \begin{array}{ccc}
k_1 & 0 & 0 \\
0 & k_2 & 0 \\
0 & 0 & k_3 \end{array} \right)\right)=\frac{\mid k_1 \mid\mid k_2 \mid}{\mid k_3 \mid^2}I_{\langle\mathrm{ICN}\rangle}(A)\\
I_{\langle\mathrm{ICN}\rangle}\left(A\left( \begin{array}{ccc}
1 & 0 & 0 \\
\frac{-\va_1\cdot(\va_2-\frac{\va_2\cdot\va_3}{\|\va_3\|^2}\va_3)}{\| (\va_2-\frac{\va_2\cdot\va_3}{\|\va_3\|^2}\va_3 \|} & 1 & 0 \\
 \frac{-\va_1\cdot\va_3}{\|\va_3\|^2} & \frac{-\va_2\cdot\va_3}{\|\va_3\|^2} & 1 \end{array} \right)\right)=I_{\langle\mathrm{ICN}\rangle}(A)
\end{align}
Using (2.2.2) we may make the columns orthogonal, and then using (2.2.1) we may further make them orthonormal. Next, using the fact that the measure $e^{-\|\vV\|^2/2}d\vV$ is $O(4)$ invariant, we may take the orthonormal set of vectors to coincide with the standard basis vectors, $\ve_1,\ve_2,\ve_3$ so that $I_{\langle\mathrm{ICN}\rangle}(A)$ may be rewritten as:
\begin{align*}
I_{\langle\mathrm{ICN}\rangle}(A)&=\frac{f(\va_1,\va_2,\va_3)}{(2\pi)^2}\int_{\vV\in\mathbb{R}^4}\frac{\mid\mathrm{Det}([\ve_1\:\ve_2\:\ve_3\:\frac{\vV}{\|\vV\|}])\mid}{\|\proj_{\vV^\perp}(\va_3)\|^3}\mathrm{exp}(-\frac{\|\vV\|^2}{2})d\vV\\
&=f(\va_1,\va_2,\va_3)I_{|\mathrm{Lk}|}(\big[\ve_1,\ve_2,\ve_3\big]),
\end{align*}
where
\begin{align*}
f([\va_1\:\va_2\:\va_3])&=\frac{\|\va_1-\frac{\va_1\cdot(\va_2-\frac{\va_2\cdot\va_3}{\|\va_3\|^2}\va_3)}{\| (\va_2-\frac{\va_2\cdot\va_3}{\|\va_3\|^2}\va_3 \|}(\va_2-\frac{\va_2\cdot\va_3}{\|\va_3\|^2}\va_3)-\frac{\va_1\cdot\va_3}{\|\va_3\|^2}\va_3\|\| \va_2-\frac{\va_2\cdot\va_3}{\|\va_3\|^2}\va_3 \|\|\va_3\|}{\|\va_3\|^3}\\
&=\frac{\sqrt{\mathrm{Det}(A^TA)}}{\|\va_3\|^3}\textrm{, and}
\end{align*} 
Now we will compute the integral:
\begin{align*}
I_{\langle\mathrm{ICN}\rangle}(A)&=\frac{1}{(2\pi)^2}f([\va_1\:\va_2\:\va_3])\int_{\mathbb{R}^4}\frac{\mid\mathrm{Det}([\ve_1\:\ve_2\:\ve_3\:\frac{\vV}{\|\vV\|}])\mid}{\|\proj_{v^\perp}(\ve_3)\|^3}\mathrm{exp}(-\frac{\|\vV\|^2}{2})d\vV
\\
&=\frac{1}{(2\pi)^2}f([\va_1\:\va_2\:\va_3])\int_{\mathbb{R}^4}\frac{\mid v_4/\|\vV\| \mid}{\|\proj_{\vV^\perp}(\ve_3)\|^3}\mathrm{exp}(-\frac{\|\vV\|^2}{2})d\vV\\
&=\frac{1}{(2\pi)^2}f([\va_1\:\va_2\:\va_3])\int_{\mathbb{R}^4}\frac{\mid v_4/\|\vV\| \mid}{(\frac{\|\vV\|^2-v_3^2}{\|v\|^2})^{3/2}}\mathrm{exp}(-\frac{\|\vV\|^2}{2})d\vV=\\
&=\frac{1}{(2\pi)^2}f([\va_1\:\va_2\:\va_3])\int_{\mathbb{R}^4}\frac{\|\vV\|^2\mid v_4\mid}{(\|\vV\|^2-v_3^2)^{3/2}}\mathrm{exp}(-\frac{\|\vV\|^2}{2})d\vV=f([\va_1\:\va_2\:\va_3])
\end{align*}
\end{Proof}
Given two space curves  $\vr_1(t)$ and $\vr_2(s)$ in $\mathbb{R}^4$, if we set 
\begin{align*}\nonumber
\va_1(t)&=\frac{d\vr_1(t)}{dt}\\
\va_2(s)&=\frac{d\vr_2(s)}{ds}\\
\va_3(t,s)&=\vr_2(s)-\vr_1(t)\textrm{ and }\\
A(s,t)&=\big[\va_1(t),\va_2(s),\va_3(s,t)\big]
\end{align*}  in the previous lemma, then we will have the following lemma:
\begin{Lemma}
Given two differentiable, non-intersecting space curves $\vr_1(t)$ and $\vr_2(s)$ in $\mathbb{R}^4$, then the value of $ICN$ averaged over all orthogonal projections to $3$ dimensional subspaces with respect to the unique normalized $O(4)$-invariant measure on $Gr(4,3)$, is given by the following integral:
\begin{equation}\nonumber
\langle\mathrm{ICN}\rangle=\frac{1}{4\pi}\int\displaylimits_{(s,t)\in\mathbb{T}^2}I_{\langle\mathrm{ICN}\rangle}(A(s,t))dsdt
\end{equation}
\end{Lemma}
\begin{Proof}
Consider the integral, $I_{\langle\mathrm{ICN}\rangle}(A)$, from the previous lemma. We may integrate over $\vV=\{v_1,v_2,v_3,v_4\}$ by changing to $4$-dimensional spherical coordinates where 
\begin{align*}
v_1&=r\cos(\phi_1)\\
v_2&=r\sin(\phi_1)\cos(\phi_2)\\
v_3&=r\sin(\phi_1)\sin(\phi_2)\cos(\phi_3)\\
v_4&=r\sin(\phi_1)\sin(\phi_2)\sin(\phi_3)
\end{align*}
so that $e^{-\|\vV\|^2/2}d\vV=e^{-r^2/2}r^3\sin^2(\phi_1)\sin(\phi_2)drd\phi_1d\phi_2d\phi_3$ and $r\in\big[0,\infty)$, $\phi_1\in\big[0,\pi\big]$, $\phi_2\in\big[0,\pi\big]$, $\phi_3\in\big[0,2\pi)$. Since the function we are integrating against $e^{-r^2/2}r^3\sin^2(\phi_1)\sin(\phi_2)drd\phi_1d\phi_2d\phi_3$ is invariant under scaling $\vV$, then we may integrate out the radial coordinate $r$ in order to obtain an integral over a $3$-sphere, and the result follows from the previous lemma since $\mathbb{S}^3\cong Gr(4,3)$ and since $Gr(4,3)$ has a unique normalized $O(4)$ invariant measure on it.
\end{Proof}
Using the techniques in \cite{Zhang} and that we will develop in the proof of \ref{mainResultTheorem}, it is straightforward to show that:
\begin{RestateTheorem}{Given two closed, differentiable, non-intersecting space curves $\gamma_1,\gamma_2:\mathbb{S}^1\rightarrow\mathbb{R}^n$ with parametrizations $\vr_1(t)$ and $\vr_2(s)$, and let $\mathrm{proj}_U(\gamma_1),\mathrm{proj}_U(\gamma_2)$ denote the orthogonal projections of $\gamma_1$ and $\gamma_2$ on to the $3$ dimensional subspace spanned by the columns of $U$, then the expected value of $\mathrm{ICN}(\mathrm{proj}_U(\gamma_1),\mathrm{proj}_U(\gamma_2))$, denoted $\langle\mathrm{ICN}\rangle$, averaged over all orthogonal projections to $3$ dimensional subspaces with respect to the unique normalized $O(n)$-invariant measure on $Gr(n,3)$, is given by the following integral:
\begin{equation}\nonumber
\langle\mathrm{ICN}\rangle=\frac{1}{4\pi}\int\displaylimits_{(s,t)\in\mathbb{T}^2}I_{\langle\mathrm{ICN}\rangle}(s,t)dsdt
\end{equation}
where
\begin{equation}\nonumber
I_{\langle\mathrm{ICN}\rangle}(s,t)=C_{\langle\mathrm{ICN}\rangle}\frac{\sqrt{\mathrm{Det}(\big[\frac{d\vr_1(t)}{dt},\frac{d\vr_2(s)}{ds},\vr_2(s)-\vr_1(t)\big]^T\big[\frac{d\vr_1(t)}{dt},\frac{d\vr_2(s)}{ds},\vr_2(s)-\vr_1(t)\big])}}{\|\vr_2(s)-\vr_1(t)\|^3},
\end{equation}}{expInterlinking}
where $C_{\langle\mathrm{ICN}\rangle}$ is a constant. 
\end{RestateTheorem}
\begin{Proof}
We will focus only on integrating over $\mathrm{Mat}_{n,3}$, since the integration over $\mathrm{Gr}(n,3)$ follows easily from the following computation and the discussion in the next chapter. Define that:
\begin{equation}\nonumber
I_{\langle\mathrm{ICN}\rangle}(\big[\va_1,\va_2,\va_3\big])=\frac{1}{(2\pi)^{3n/2}}\int\displaylimits_{U\in\mathrm{Mat}_{n,3}}\frac{\sqrt{\mathrm{Det}(\big[\va_1,\va_2,\va_3\big]^TU(U^TU)^{-1}U^T\big[\va_1,\va_2,\va_3\big])}}{(\va_3^TU(U^TU)^{-1}U^T\va_3)^{3/2}}dU,
\end{equation}
where we think of $U\in\mathrm{Mat}_{n,3}$ as a matrix with columns $\vu_1,\vu_2$ and $\vu_3$, and 
\begin{equation}\nonumber
dU=\mathrm{exp}(-\frac{1}{2}\mathrm{Tr}(U^TU))d\vu_1d\vu_2d\vu_3
\end{equation}
Next, apply the Gram-Schmidt process to the ordered set of vectors $\{\va_3,\va_2,\va_1\}$ to obtain the orthonormal set $\{\vq_3,\vq_2,\vq_1\}$, so that we have:
\begin{align*}\nonumber
I_{\langle\mathrm{ICN}\rangle}(\big[\va_1,\va_2,\va_3\big])&=\frac{1}{(2\pi)^{3n/2}}\frac{\mathrm{Det}(A^TA)^{1/2}}{\|\va_3\|^3}\int\displaylimits_{U\in\mathrm{Mat}_{n,3}}\frac{\sqrt{\mathrm{Det}(\big[\vq_1,\vq_2,\vq_3\big]^TU(U^TU)^{-1}U^T\big[\vq_1,\vq_2,\vq_3\big])}}{(\vq_3^TU(U^TU)^{-1}U^T\vq_3)^{3/2}}dU\\
&=\frac{\mathrm{Det}(A^TA)^{1/2}}{\|\va_3\|^3}I_{\langle\mathrm{ICN}\rangle}(\big[\vq_1,\vq_2,\vq_3\big]),
\end{align*}
where in the above $A=\big[\va_1,\va_2,\va_3\big]$. Since the measure $dU$ is $O(n)$ invariant, then we may take the set of vectors $\{\vq_1,\vq_2,\vq_3\}$ to coincide with the first $3$ standard basis vectors: $\{\ve_1,\ve_2,\ve_3\}$. The result follows by setting $C_{\langle\mathrm{ICN}\rangle}=I_{\langle\mathrm{ICN}\rangle}(\ve_1,\ve_2,\ve_3)$
\end{Proof}
\begin{Remark}
Notice that if we considered $\langle\mathrm{Lk}\rangle$ instead, the integral would have vanished.
\end{Remark}

\subsection{Total Curvature}
For pedagogical purposes, we will apply  a similar analysis to find the expected total curvature of a knot. For the most part, the following is due to Fary \cite{Fary}, Milnor \cite{Milnor}, and appears in a more generalized context by Sullivan \cite{Sullivan}. In this section, we will simply show that the aforementioned arguments easily fit in to our framework discussed previously.  To do so, recall that if we are given a twice continuously differentiable space curve, C, with parametrization $\vr(t):\mathbb{S}^1\rightarrow\mathbb{R}^3$, we may compute the total curvature, $\kappa(\mathrm{C})$, using a Gauss-like integral:
\begin{equation}\nonumber
\kappa(\mathrm{C})=\int_{0}^{2\pi}\frac{|\vr'(t)\times\vr''(t)|}{\|\vr'(t)\|^2}dt=\int_{0}^{2\pi}\frac{\sqrt{\|\vr'(t)\|^2\|\vr''(t)\|^2-(\vr'(t)\cdot\vr''(t))^2}}{\|\vr'(t)\|^2}dt,
\end{equation}
First we will prove a lemma analogous to \ref{prevlemma}:
\begin{Lemma}
Given a matrix $A=\big[\va_1,\va_2\big]\in\mathrm{Mat}_{n,2}$ and the function
\begin{equation}\nonumber
I_{\kappa(\mathrm{C})}(A)=\frac{1}{(2\pi)^{3n/2}}\int_{\mathbb{R}^n}\int_{\mathbb{R}^n}\int_{\mathbb{R}^n}\frac{\mid\mathrm{Det}(A^TU(U^TU)^{-1}U^TA)\mid^{1/2}}{(\va_1^TU(U^TU)^{-1}U^T\va_1)^{2}}dU
\end{equation}
where $U$ denotes the matrix $U=\big[\vu_1,\vu_2,\vu_3\big]$ and $dU=exp(-\frac{1}{2}(\|\vu_1\|^2+\|\vu_2\|^2+\|\vu_3\|^2))d\vu_1d\vu_2d\vu_3$, then
\begin{equation}\nonumber
I_{\kappa(\mathrm{C})}(A)=k\frac{\sqrt{\|\va_1\|^2\|\va_2\|^2-(\va_1\cdot\va_2)^2}}{\|\va_1\|^2}
\end{equation} where $k$ is a constant.  
\end{Lemma}
\begin{Proof}
As in the proof of lemma \ref{firstlemma}, we may use the numerous multilinearity and scaling properties to write $I_{\kappa(\mathrm{C})}(A)$ as:
\begin{equation}\nonumber
I_{\kappa(\mathrm{C})}(A)=\frac{1}{(2\pi)^{3n/2}}\frac{\sqrt{|\mathrm{Det}(A^TA)|}}{\|\va_1\|^2}\int_{\mathbb{R}^n}\int_{\mathbb{R}^n}\int_{\mathbb{R}^n}\frac{\mid\mathrm{Det}(\widetilde{U}(U^TU)^{-1}\widetilde{U}^T)\mid^{1/2}}{(\widetilde{U_{(1)}}(U^TU)^{-1}\widetilde{U_{(1)}}^T)^{2}}dU
\end{equation}
where $U$ is an $n$ by $3$ matrix with columns $\vu_1,\vu_2$ and $\vu_3$, and
\begin{equation}\nonumber
\widetilde{U}=\left(\begin{array}{ccc}
u_{11} & u_{21} & u_{31} \\
u_{12} & u_{22} & u_{32}\end{array}\right)
\end{equation}
and
\begin{equation}\nonumber
\widetilde{U_{(1)}}=\left(\begin{array}{ccc}
u_{11} & u_{21} & u_{31}\end{array}\right),
\end{equation}
where $u_{ij}$ is the $j^{th}$ entry of the vector $\vu_i$.

The result follows by setting $k=\frac{1}{(2\pi)^{3n/2}}\int_{\mathbb{R}^n}\int_{\mathbb{R}^n}\int_{\mathbb{R}^n}\frac{\mid\mathrm{Det}(\widetilde{U}(U^TU)^{-1}\widetilde{U}^T)\mid^{1/2}}{(\widetilde{U_{(1)}}(U^TU)^{-1}\widetilde{U_{(1)}}^T)^{2}}dU.$
\end{Proof}
The following is essentially due to Fary and Milnor:
\begin{Theorem}
Let $C$ be a twice-differentiable space curve with parametrization $\vr(t):\mathbb{S}^1\rightarrow\mathbb{R}^n$ and $\textrm{proj}_{U}(\vr(t))$ be the orthogonal projection of $\vr(t)$ to the $3$-dimensional subspace spanned by the columns of $U$, then the value of the total curvature, $\kappa(\textrm{proj}_{U}\vr(t)))$, averaged over all projections to $3$ dimensional subspaces with respect to the Gaussian measure on $\mathrm{Mat}_{n,3}$ is given by the integral:
\begin{equation}\nonumber
\langle\kappa(\mathrm{C})\rangle=\int\displaylimits_{t\in\mathbb{T}^1}I_{\kappa(\mathrm{C})}(t)dt
\end{equation}
where
\begin{eqnarray*}\label{eq:result1THM}
I_{\kappa(\mathrm{C})}(t)=\frac{\sqrt{\|\vr'(t)\|^2\|\vr''(t)\|^2-(\vr'(t)\cdot\vr''(t))^2}}{\|\vr'(t)\|^2}\nonumber
\end{eqnarray*}
\end{Theorem}
\begin{Proof}
Given the twice differentiable curve $\vr(t):\mathbb{S}^1\rightarrow\mathbb{R}^n$, we set 
\begin{eqnarray*}\nonumber
\va_1(t)=\frac{d\vr_1(t)}{dt}\\
\va_2(t)=\frac{d^2\vr_1(t)}{dt^2}\\
A(t)=\big[\va_1(t),\va_2(t)\big]
\end{eqnarray*}
This way we will have that:
\begin{equation}\nonumber
\langle\kappa(\mathrm{C})\rangle=\int_0^{2\pi}I_{\kappa(C)}(A(t))dt
\end{equation} 
Using the previous lemma we have that:
\begin{equation}\nonumber
\langle\kappa(\mathrm{C})\rangle=\int_0^{2\pi}k\frac{\sqrt{\|\va_1(t)\|^2\|\va_2(t)\|^2-(\va_1(t)\cdot\va_2(t))^2}}{\|\va_1(t)\|^2}dt
\end{equation}
Since the above holds for any $\vr(t)$, then in particular we may take $\vr(t)=\cos(t)\ve_1+\sin(t)\ve_2$, so that $\langle\kappa(\mathrm{C})\rangle=2\pi$ since the image of almost all projections will be a planar ellipse. It is also easy to compute that $\int_0^{2\pi}k\frac{\sqrt{\|\va_1(t)\|^2\|\va_2(t)\|^2-(\va_1(t)\cdot\va_2(t))^2}}{\|\va_1(t)\|^2}dt=2\pi$, so that we see:
\begin{equation}\nonumber
\frac{1}{(2\pi)^{3n/2}}\int_{\mathbb{R}^n}\int_{\mathbb{R}^n}\int_{\mathbb{R}^n}\frac{\mid\mathrm{Det}(\tilde{U}(U^TU)^{-1}\tilde{U}^T)\mid^{1/2}}{(U_{(1)}^T(U^TU)^{-1}U^TU_{(1)}^T)^{2}}\mathrm{exp}^{-\frac{1}{2}(\|\vu_1\|^2+\|\vu_2\|^2+\|\vu_3\|^2)}d\vu_1d\vu_2d\vu_3=1
\end{equation} 
and therefore
\begin{equation}\label{dataintegral}
\langle\kappa(\mathrm{C})\rangle=\int_0^{2\pi}\frac{\sqrt{\|\vr'(t)\|^2\|\vr''(t)\|^2-(\vr'(t)\cdot\vr''(t))^2}}{\|\vr'(t)\|^2}dt
\end{equation}
\end{Proof}
That is, the expected total curvature of one of the knot projections is simply the total curvature of the space curve C that we are randomly projecting to lower dimensions. 
    
        \section[Integration Over Configuration Spaces]{Integration Over the Configuration Space}
        \label{sec:LabelForChapter2:Sectionn}     
        
Up to this point we have only considered the integration over all the $3$ dimensional subspaces, and now we will consider how to specify the initial data, so as to bound $\langle\kappa(C)\rangle$ and $\langle\mathrm{ICN}\rangle$ and therefore bound the complexity of the curvature and the average inter-crossing number. We will now prove the following theorem:
\begin{RestateTheorem}{
Let $\vr_1(t):\mathbb{S}^1\rightarrow\mathbb{R}^{2n+1}$ and $\vr_2(s):\mathbb{S}^1\rightarrow\mathbb{R}^{2n+1}$ be given by:
\begin{align*}\nonumber
\vr_1(t)&=c_0\textbf{e}_0+\sum_{k=1}^{n}(c_k\cos(kt)\ve_{2k-1}+c_k\sin(kt)\ve_{2k})\\
&=(c_0,c_1\cos(t),c_1\sin(t),c_2\cos(2t),c_2\sin(2t),...,c_n\cos(nt),c_n\sin(nt))\textrm{, and}\\
\vr_2(s)&=d_0\textbf{e}_0+\sum_{k=1}^{n}(d_k\cos(ks)\ve_{2k-1}+d_k\sin(ks)\ve_{2k})\\
&=(d_0,d_1\cos(s),d_1\sin(s),d_2\cos(2s),d_2\sin(2s),...,d_n\cos(ns),d_n\sin(ns))\nonumber
\end{align*}
then $\langle ICN\rangle$ is finite and satisfies the following bound:
\begin{equation}\nonumber
\langle ICN\rangle\leq C_{\langle ICN\rangle}\frac{\sqrt{(\sum_{j=0}^n j^2c_j^2)(\sum_{j=0}^n j^2d_j^2)}}{\mathrm{min}_{(s,t)\in\mathbb{T}^2}\|\vr_2(s)-\vr_1(t)\|^2}
\end{equation}}{GaussianConfigLk}
\end{RestateTheorem}
\begin{Proof}
By \eqref{expInterlinking} we have:
\begin{align}\nonumber
\langle\mathrm{ICN}\rangle&\leq C_{\langle\mathrm{ICN}\rangle}\int_0^{2\pi}\int_0^{2\pi}\frac{\|\frac{d\vr_1(t)}{dt}\|\|\frac{d\vr_2(s)}{ds}\|}{4\pi \mathrm{min}_{(s,t)\in\mathbb{T}^2}\|\vr_2(s)-\vr_1(t)\|^2}dtds\\
&\leq C_{\langle\mathrm{ICN}\rangle}\int_0^{2\pi}\int_0^{2\pi}\frac{\sqrt{(\sum_{j=0}^n j^2c_j^2)(\sum_{j=0}^n j^2d_j^2)}}{4\pi \mathrm{min}_{(s,t)\in\mathbb{T}^2}\|\vr_2(s)-\vr_1(t)\|^2}dtds\nonumber\\
&=C_{\langle\mathrm{ICN}\rangle}\pi\frac{\sqrt{(\sum_{j=0}^n j^2c_j^2)(\sum_{j=0}^n j^2d_j^2)}}{\mathrm{min}_{(s,t)\in\mathbb{T}^2}\|\vr_2(s)-\vr_1(t)\|^2}dtds,\nonumber
\end{align}
so that we may bound $\langle\mathrm{ICN}\rangle$:
\begin{equation}\nonumber
\langle\mathrm{ICN}\rangle\leq C_{\langle\mathrm{ICN}\rangle}\pi\frac{\sqrt{(\sum_{j=0}^n j^2c_j^2)(\sum_{j=0}^n j^2d_j^2)}}{\mathrm{min}_{(s,t)\in\mathbb{T}^2}\|\vr_2(s)-\vr_1(t)\|^2}
\end{equation}
\end{Proof}
If we make even stronger assumptions, for example that the initial space curves be orthogonal, then we can significantly improve the above bound. For example, if the initial space curves $\vr_1(t),\vr_2(s):\mathbb{S}^1\rightarrow\mathbb{R}^{4N+2}$ are orthogonal:
\begin{align}\label{orthogonalData}
\vr_1(t)&=c_0\ve_1+\sum\displaylimits_{k=1}^{N}(c_k\cos(kt)\ve_{4k-1}+c_k\sin(kt)\ve_{4k})\\
\vr_2(s)&=d_0\ve_2+\sum\displaylimits_{k=1}^{N}(d_k\cos(ks)\ve_{4k+1}+d_k\sin(ks)\ve_{4k+2}),
\end{align}
then by Hadamard's inequality we have that:
\begin{equation}\nonumber
\langle\mathrm{ICN}\rangle\leq\frac{1}{4\pi}\int\displaylimits_{s=0}^{2\pi}\int\displaylimits_{t=0}^{2\pi} C_{\mathrm{ICN}}\frac{\|\vr'_1(t)\|\vr'_2(s)\|}{\|\vr_2(s)-\vr_1(t)\|^2}dsdt=\pi C_{\mathrm{ICN}}\frac{\sqrt{(\sum\displaylimits_{k=1}^Nk^2c_k^2)(\sum\displaylimits_{k=1}^Nk^2d_k^2)}}{\sum\displaylimits_{k=0}^Nc_k^2+\sum\displaylimits_{k=0}^Nd_k^2}
\end{equation}
Similarly, for the total curvature we have the following theorem:
\begin{Theorem}
Let $\vr(t):\mathbb{S}^1\rightarrow\mathbb{R}^{2n+1}$ be given by:
\begin{align}\label{coeffs}
\vr(t)&=c_0\textbf{e}_0+\sum_{k=1}^{n}(c_k\cos(kt)\ve_{2k-1}+c_k\sin(kt)\ve_{2k})\\
&=(c_0,c_1\cos(t),c_1\sin(t),c_2\cos(2t),c_2\sin(2t),...,,c_n\cos(nt),c_n\sin(nt)),\nonumber
\end{align}
then $\langle\kappa(\mathrm{C})\rangle$ is finite and satisfies the following bound:
\begin{equation}\nonumber
\langle\kappa(\mathrm{C})\rangle\leq2\pi\sqrt{\frac{\sum\displaylimits_{k=0}^nc_k^2k^2}{\sum\displaylimits_{k=0}^nc_k^2k^4}}
\end{equation}
\end{Theorem}
\begin{Proof}
For $\langle\kappa(C)\rangle$, notice first that:
\begin{align}\label{curveIn}
\langle\kappa(C)\rangle&\leq\int_0^{2\pi}\frac{\|\vr''(t)\|}{\|\vr'(t)\|}dt\\
\|\vr'(t)\|^2&=\sum\displaylimits_{k=1}^nc_k^2k^2\nonumber\\
\|\vr''(t)\|^2&=\sum\displaylimits_{k=1}^nc_k^2k^4\nonumber
\end{align}
The quotient $\frac{\|\vr''(t)\|}{\|\vr'(t)\|}$ has no dependence on $t$, and so the result follows from \eqref{curveIn} by integrating over $t$.
\end{Proof}
\begin{Remark}
As mentioned in the introduction, since we are considering arbitrary space curves in $\mathbb{R}^n$, and since the formulas involved did not have any explicit dependence on $n$, then we may allow $n$ to approach infinity and consider a space curve of the form:
\begin{eqnarray}
\vr(t)=c_0\ve_0+\sum_{n=1}^{\infty}(c_n\cos(nt)\ve_{2n-1}+c_n\sin(nt)\ve_{2n})\\
=(c_0,c_1\cos(t),c_1\sin(t),c_2\cos(2t),c_2\sin(2t),...,,c_n\cos(nt),c_n\sin(nt),...)\nonumber
\end{eqnarray}
In this way, if we define the sequences $\textbf{c}=\{c_n\}_{n=0}^\infty$, $\textbf{c}'=\{nc_n\}_{n=0}^\infty$, $\textbf{c}''=\{n^2c_n\}_{n=0}^\infty$  and stipulate that $\|\textbf{c}''\|_{l^2}<\infty$ (that is, the sequence \textbf{c} decays faster than $n^{(-5-\epsilon)/2}$), then the quotient $\frac{\|\vr''(t)\|}{\|\vr'(t)\|}$ will be finite and will give a bound for $\langle\kappa(C)\rangle$ so that with high probability the projections will have bounded total curvature given by:
\begin{equation}\nonumber
\langle\kappa(C)\rangle\leq\int_0^{2\pi}\frac{\|\vr''(t)\|}{\|\vr'(t)\|}dt=\int_0^{2\pi}\sqrt{\frac{\sum_{n=0}^\infty n^4c_n^2}{\sum_{n=0}^\infty n^2c_n^2}}dt
\end{equation} 
and therefore:
\begin{equation}\nonumber
\langle\kappa(C)\rangle\leq2\pi\frac{\|\textbf{c}''\|_{l^2}}{\|\textbf{c}'\|_{l^2}}
\end{equation}
A similar result holds for $\langle\mathrm{ICN}\rangle$. In particular, using the data \eqref{orthogonalData}, we will have that:
\begin{equation}
\langle\mathrm{ICN}\rangle\leq \pi C_{\mathrm{ICN}}\frac{\|\textbf{c}'\|_{l^2}\|\textbf{d}'\|_{l^2}}{\|\textbf{c}\|_{l^2}^2+\|\textbf{d}\|_{l^2}^2},
\end{equation}
so that the sequence \textbf{c} must decay faster than $n^{(-3-\epsilon)/2}$. 
\end{Remark}

    \newchapter{Second Moment of the Linking Number}{Second Moment of the Linking Number}{Second Moment of the Linking Number}
    \label{sec:LabelForChapter3}

        \section[Integration Over 2-Dimensional Subspaces and the Winding Number]{Integration Over 2-Dimensional Subspaces and the Winding Number}
        \label{sec:LabelForChapter3:Section1}

Before we go on to compute the second moment of the linking number, we would like to first consider a simpler case: that of computing the second moment of the winding number for a model of random plane curves closely related to the model of random knots discussed in the previous sections.  In particular, we fix a space curve in $\mathbb{R}^n$, and project the space curve onto random 2 dimensional subspaces of $\mathbb{R}^n$, chosen from the unique $O(n)$ invariant measure on $Gr(n,2)$. 

For a differentiable plane curve, $\vr(t)$, recall that the winding number $W$ can be computed as a degree of a map, which may be written as an integral of the form:

\begin{equation}\label{windingDEF}
W=\frac{1}{2\pi}\int\displaylimits_{\mathbb{S}^1}\frac{\mathrm{Det\left[\vr(t),\vr'(t)\right]}}{\|\vr(t)\|^2}dt
\end{equation}
With this definition, we can compute the second moment $\langle W^2\rangle$ of the winding number $W$, for a given fixed differentiable map $\vr:\mathbb{S}^1\rightarrow\mathbb{R}^n$. In particular, we compute:
\begin{equation}\label{second_moment}
	\langle W^2\rangle=\frac{1}{(2\pi)^n}\int\displaylimits_{U\in\mathrm{Mat}_{n,2}}\left(\frac{1}{(2\pi)^2}\int\displaylimits_{\mathbb{S}^1\times\mathbb{S}^1}\frac{Det(U^T\left[\vr'(t),\vr(t)\right])Det(U^T\left[\vr'(s),\vr(s)\right])}{\|U^T\vr(t)\|^2\|U^T\vr(s)\|^2}dsdt\right)e^{\frac{-Tr(U^TU)}{2}}dU
\end{equation}
where $U$ is a Gaussian random matrix. Next, switch the order of integration in (\ref{second_moment}) to obtain:
\begin{equation}\label{second_moment2}
	\langle W^2\rangle=\frac{1}{(2\pi)^2}\int\displaylimits_{\mathbb{S}^1\times\mathbb{S}^1}\left(\frac{1}{(2\pi)^n}\int\displaylimits_{U\in\mathrm{Mat}_{n,2}}\frac{Det(U^T\left[\vr'(t),\vr(t)\right])Det(U^T\left[\vr'(s),\vr(s)\right])}{\|U^T\vr(t)\|^2\|U^T\vr(s)\|^2}e^{\frac{-Tr(U^TU)}{2}}dU\right)dsdt
\end{equation}
Now we focus on the Gaussian integral in parentheses, and define:
\begin{equation}\label{Gaussian}
I(\va_1,\va_1',\va_2,\va_2')=\frac{1}{(2\pi)^n}\int\displaylimits_{U\in\mathrm{Mat}_{n,2}}\frac{Det(U^T\left[\va_1,\va_2\right])Det(U^T\left[\va_1',\va_2'\right])}{\|U^T\va_2\|^2\|U^T\va_2'\|^2}e^{\frac{-Tr(U^TU)}{2}}dU
\end{equation}
where $\va_1,\va_1',\va_2,\va_2'\in\mathbb{R}^n$ so that (\ref{second_moment}) becomes:
\begin{equation}\label{second_moment_simplified}
	\langle W^2\rangle=\int\displaylimits_{\mathbb{S}^1\times\mathbb{S}^1}I\left(\vr'(t),\vr'(s),\vr(t),\vr(s)\right)dsdt
\end{equation}
To compute (\ref{Gaussian}), first normalize $\va_2$ and $\va_2'$ to obtain the unit vectors $\vq_2$ and $\vq_2'$, and choose a basis $B=\{\vb_1,\vb_2,...,\vb_n\}$ where $\vb_1$ and $\vb_2$ are in the span of $\vq_2$ and $\vq_2'$. In what follows, for convenience we decide to take :
\begin{align*}
\vb_1&=\frac{\va_2\|\va_2'\|+\va_2'\|\va_2\|}{\sqrt{2\left(\|\va_2\|^2\|\va_2'\|^2+\|\va_2\|\|\va_2'\|(\va_2\cdot\va_2')\right)}}\\
\vb_2&=\frac{\va_2\|\va_2'\|-\va_2'\|\va_2\|}{\sqrt{2\left(\|\va_2\|^2\|\va_2'\|^2-\|\va_2\|\|\va_2'\|(\va_2\cdot\va_2')\right)}}
\end{align*}
and we further define that:
\begin{align}\label{angle}
a&=\frac{\va_2}{\|\va_2\|}\cdot\vb_1=|\cos(x/2)|\nonumber\\
b&=\frac{\va_2'}{\|\va_2'\|}\cdot\vb_2=|\sin(x/2)|,
\end{align}
where $x$ is the angle between $\va_2$ and $\va_2'$. In this way we have $a$ and $b$ are given just in terms of dot products between $\va_1,\va_1',\va_2,\va_2'$, and moreover $a^2+b^2=1$. With these definitions, we then expand the integrand using the multilinearity of the determinant to obtain:
\begin{equation}\nonumber
I(\va_1,\va_1',\va_2,\va_2')=\frac{1}{\|\va_2\|\|\va_2'\|}\sum\displaylimits_{i,j=1}^na_{1i}a_{1j}'I(\ve_i,\ve_j,a\ve_1+b\ve_2,a\ve_1-b\ve_2),
\end{equation}
where hereinafter we define $a_{1i}=\va_1\cdot\vb_i$ and similarly $a'_{1i}=\va_1'\cdot\vb_i$, are the components of $\va_1$ and $\va'_1$ in the basis $B$.
There are a few cases to consider when we compute $I(\ve_i,\ve_j,a\ve_1+b\ve_2,a\ve_1-b\ve_2)$. The integral $I(\ve_i,\ve_j,a\ve_1+b\ve_2,a\ve_1-b\ve_2)$ has a simple form:

\begin{align*}\label{gaussian}
I(\ve_i,\ve_j,&a\ve_1+b\ve_2,a\ve_1-b\ve_2)\nonumber\\
&=\int\displaylimits_{U\in\mathrm{Mat}_{n,2}}\frac{Det(\begin{bmatrix}
    u_{1i}       & au_{11}+bu_{12} \\
    u_{2i}       & au_{21}+bu_{22}
\end{bmatrix})Det(\begin{bmatrix}
    u_{1j}       & au_{11}-bu_{12} \\
    u_{2j}       & au_{21}-bu_{22}
\end{bmatrix})}{\left((au_{11}+bu_{12})^2+(au_{21}+bu_{22})^2\right)\left((au_{11}-bu_{12})^2+(au_{21}-bu_{22})^2\right)}e^{\frac{-Tr(U^TU)}{2}}dU
\end{align*}
Applying a combination of symmetries of the integrand and a few coordinate changes we obtain that:
\begin{align*}
\|\va_2\|\|\va_2'\|I(\va_1,\va_1',\va_2,\va_2')&=\left(a_{11}a_{12}'-a_{12}a_{11}'\right)I(\ve_1,\ve_2,a\ve_1+b\ve_2,a\ve_1-b\ve_2)\\
&+\left(a_{11}a_{11}'-\frac{a^2}{b^2}a_{12}a_{12}'\right)I(\ve_1,\ve_1,a\ve_1+b\ve_2,a\ve_1-b\ve_2)\\
&+\left(\sum\displaylimits_{k=3}^na_{1k}a_{1k}'\right)I(\ve_3,\ve_3,a\ve_1+b\ve_2,a\ve_1-b\ve_2)\\
\end{align*}

In the above, notice that $\sum\displaylimits_{k=3}^na_{1k}a_{1k}'=\va_1\cdot\va_1'-a_{11}a_{11}'-a_{12}a_{12}'$, so we do not need to know $\{\vb_3,...,\vb_n\}$ explicitly, which is a computational convenience. We can in fact simplify the above even further by writing $I(\ve_1,\ve_2,a\ve_1+b\ve_2,a\ve_1-b\ve_2)=\frac{a}{b}I(\ve_1,\ve_1,a\ve_1+b\ve_2,a\ve_1-b\ve_2)$ so that the computation has now simplified to computing two Gaussian integrals: $I(\ve_1,\ve_1,a\ve_1+b\ve_2,a\ve_1-b\ve_2)$ and $I(\ve_3,\ve_3,a\ve_1+b\ve_2,a\ve_1-b\ve_2)$.

Both integrals can be simplified to integrals over $\mathbb{R}^4$ (two variables of $I(\ve_3,\ve_3,a\ve_1+b\ve_2,a\ve_1-b\ve_2)$ can be easily integrated out), and the radial portion may be integrated out leaving two integrals over a $3$-sphere. From here, we choose a coordinate system on the $3$-sphere:
\begin{align*}
u_{11}&=\cos(\sigma)\cos(\theta)\\
u_{12}&=\cos(\sigma)\sin(\theta)\\
u_{21}&=\sin(\sigma)\sin(\phi)\\
u_{22}&=\sin(\sigma)\sin(\phi)\\
\end{align*}
where $\sigma\in(0,\pi/2)$, $\eta,\phi\in(0,2\pi)$. This particular coordinate system is convenient, and was chosen with the form of the integrand's denominator in mind. After integrating over the sphere we obtain:
\begin{equation}\label{final}
I(\va_1,\va_1',\va_2,\va_2')=\frac{\log(4a^2b^2)\left((b^2-a^2)(\va_1\cdot\va_1')\|\va_2\|\|\va_2'\|+(\va_1\cdot\va_2')(\va_1'\cdot\va_2)\right)}{2(b^2-a^2)^2\|\va_2\|^2\|\va_2'\|^2}
\end{equation}
where $a$ and $b$ are from (\ref{angle}).  Now we will write (\ref{final}) in terms of $\vr(t)$ and $\vr(s)$:
\begin{equation}\nonumber
I\left(\vr'(t),\vr'(s),\vr(t),\vr(s)\right)=\frac{\log(\sin(x)^2)\left(-\cos(x)\left(\vr'(t)\cdot\vr'(s)\right)\|\vr(t)\|\|\vr(s)\|+(\vr'(t)\cdot\vr(s))(\vr'(s)\cdot\vr(t)\right)}{2\cos(x)^2\|\vr(t)\|^2\|\vr(s)\|^2}
\end{equation}
where $x(s,t)$ is the angle between $\vr(t)$ and $\vr(s)$.  Simplifying even further we obtain:
\begin{equation}\label{final}
I\left(\vr'(t),\vr'(s),\vr(t),\vr(s)\right)=\frac{\log(\frac{\|\vr(t)\|^2\|\vr(s)\|^2-\left(\vr(t)\cdot\vr(s)\right)^2}{\|\vr(t)\|^2\|\vr(s)\|^2})\left((\vr'(t)\cdot\vr(s))(\vr'(s)\cdot\vr(t))-(\vr(t)\cdot\vr(s))(\vr'(t)\cdot\vr'(s)))\right)}{2\left(\vr(t)\cdot\vr(s)\right)^2}
\end{equation}

\begin{align}\label{final2}
I\left(\vr'(t),\vr'(s),\vr(t),\vr(s)\right)&=\frac{\log(\frac{\|\vr(t)\|^2\|\vr(s)\|^2-\left(\vr(t)\cdot\vr(s)\right)^2}{\|\vr(t)\|^2\|\vr(s)\|^2})\left((\vr'(t)\cdot\vr(s))(\vr'(s)\cdot\vr(t))-(\vr(t)\cdot\vr(s))(\vr'(t)\cdot\vr'(s)))\right)}{2\left(\vr(t)\cdot\vr(s)\right)^2}\\
&=\frac{\log(\frac{\|\vr(t)\|^2\|\vr(s)\|^2-\left(\vr(t)\cdot\vr(s)\right)^2}{\|\vr(t)\|^2\|\vr(s)\|^2})\mathrm{Det}\left(\begin{bmatrix}
    \vr'(t)\cdot\vr(s) & \vr(t)\cdot\vr(s)\\
    \vr'(t)\cdot\vr'(s) & \vr'(s)\cdot\vr(t) 
\end{bmatrix}\right)}{2\left(\vr(t)\cdot\vr(s)\right)^2}
\end{align}
In summary, we simplified a high dimensional Gaussian integral into several smaller integrals using the symmetries of the integrand, in order to obtain a function that depended only on the pairwise dot products of a set of vectors. With this technique developed, in the next section, we generalize this argument in order to find the second moment of the linking number for the model of random links discussed in previous sections.

        \section[Integration Over 3-Dimensional Subspaces]{Integration Over 3-Dimensional Subspaces}
        \label{sec:LabelForChapter3:Section2}

As stated in the introduction, in order to find the second moment of the linking number for a link whose components $\gamma_1$ and $\gamma_2$ are parametrized by $\vr_1(t)$ and $\vr_2(s)$, we must consider the expectation value of an integral of the form:
\begin{align*}\nonumber
&\mathrm{Lk}^2(\vr_1,\vr_2)=\\
&\frac{1}{(4\pi)^2}\int\displaylimits_{(s,t)\in\mathbb{T}^2}\int\displaylimits_{(s',t')\in\mathbb{T}^2}\frac{\mathrm{Det}\big[\dot{\vr}_1(t),\dot{\vr}_2(s),\vr_2(s)-\vr_1(t)\big]\mathrm{Det}\big[\dot{\vr}_1(t'),\dot{\vr}_2(s'),\vr_2(s')-\vr_1(t')\big]}{\|\vr_2(s)-\vr_1(t)\|^3\|\vr_2(s')-\vr_1(t')\|^3}dsdtds'dt',
\end{align*}
So, to begin, given the vectors: $\va_1,\va_2,\va_3,\va_1',\va_2',\va_3'\in\mathbb{R}^n$ define the matrices:
\begin{align*}
A&=\left[\va_1,\va_2,\va_3\right]\nonumber\\
A'&=\left[\va_1',\va_2',\va_3'\right]\nonumber
\end{align*}
and define a function $I_{\langle\mathrm{Lk}^2\rangle}:\mathrm{Mat}_{n,3}\times\mathrm{Mat}_{n,3}\rightarrow\mathbb{R}$ by:
\begin{equation}\label{LK2def}
I_{\langle\mathrm{Lk}^2\rangle}(A,A')=\frac{1}{(2\pi)^{\frac{3n}{2}}}\int\displaylimits_{U\in\mathrm{Mat}_{n,3}}\frac{\mathrm{Det}(U^TA)\mathrm{Det}(U^TA')}{\|U^T\va_3\|^3\|U^T\va_3'\|^3}e^{-\mathrm{Tr}\left(U^TU\right)/2}dU.
\end{equation}
With this definition we have that:
\begin{align*}\nonumber
\langle\mathrm{Lk}^2(\vr_1,\vr_2)\rangle&=\frac{1}{(4\pi)^2}\int\displaylimits_{(s,t)\in\mathbb{T}^2}\int\displaylimits_{(s',t')\in\mathbb{T}^2}I_{\langle\mathrm{Lk}^2\rangle}(\textbf{R}(s,t),\textbf{R}(s',t'))dsdtds'dt'\\
\textrm{ where }\textbf{R}(s,t)&=\big[\dot{\vr}_1(t),\dot{\vr}_2(s),\vr_2(s)-\vr_1(t)\big]\textrm{ and }\textbf{R}(s',t')=\big[\dot{\vr}_1(t'),\dot{\vr}_2(s'),\vr_2(s')-\vr_1(t')\big].
\end{align*}
\begin{Remark}
For conciseness, hereinafter we will suppress the dependence on the initial space curves $\vr_1(t)$ and $\vr_2(s)$ and will henceforth write $\langle\mathrm{Lk}^2\rangle$ to denote $\langle\mathrm{Lk}^2(\vr_1,\vr_2)\rangle$ when the dependence on  $\vr_1(t)$ and $\vr_2(s)$ is clear.
\end{Remark}
In order to find $\langle\mathrm{Lk}^2\rangle$, it will be helpful to determine an explicit form for the integral $I_{\langle\mathrm{Lk}^2\rangle}(A,A')$ and then we will need to integrate over the configuration space of pairs of tuples of points on distinct components of the link.  We will now focus on performing the integral in the expression \eqref{LK2def} for $I_{\langle\mathrm{Lk}^2\rangle}(A,A')$:

\begin{Lemma}\label{codim1lemma}
Given matrices $A=\big[\va_1,\va_2,\va_3\big],A'=\big[\va_1',\va_2',\va_3'\big]\in\mathrm{Mat}_{n,3}$ and the function
\begin{equation}\nonumber
I_{\langle\mathrm{Lk}^2\rangle}(A,A')=\frac{1}{(2\pi)^{\frac{3n}{2}}}\int\displaylimits_{U\in\mathrm{Mat}_{n,3}}\frac{\mathrm{Det}(U^TA)\mathrm{Det}(U^TA')}{\|U^T\va_3\|^3\|U^T\va_3'\|^3}e^{-\mathrm{Tr}\left(U^TU\right)/2}dU.
\end{equation}
then we have
\begin{align}\label{GaussianResult}
&\frac{\pi}{4}I_{\langle\mathrm{Lk}^2\rangle}(A,A')=\nonumber\\
&\frac{(\va_3\cdot\va_3')-\sqrt{\|\va_3\|^2\|\va_3'\|^2-(\va_3\cdot\va_3')^2}\sin^{-1}\left(\frac{\va_3\cdot\va_3'}{\|\va_3\|\|\va_3'\|}\right)}{(\va_3\cdot\va_3')^3\sqrt{\|\va_3\|^2\|\va_3'\|^2-(\va_3\cdot\va_3')^2}}\mathrm{Det}\left(\begin{bmatrix}
    \va_1\cdot\va_1'       & \va_1\cdot\va_2' & \va_1\cdot\va_3'   \\
    \va_2\cdot\va_1'       & \va_2\cdot\va_2' & \va_2\cdot\va_3'   \\
    \va_3\cdot\va_1'       & \va_3\cdot\va_2' & \va_3\cdot\va_3'   \\
\end{bmatrix}\right)
\end{align}  
\end{Lemma}
\begin{Proof}
To begin, normalize the vectors $\va_3$ and $\va_3'$ to form the vectors $\vq_3$ and $\vq_3'$ respectively, and then express \eqref{LK2def} as:
\begin{align}\label{unnorm_to_norm}
I_{\langle\mathrm{Lk}^2\rangle}(A,A')&=\frac{1}{(2\pi)^{\frac{3n}{2}}}\frac{1}{\|\va_3\|^2\|\va_3'\|^2}\int\displaylimits_{U\in\mathrm{Mat}_{n,3}}\frac{\mathrm{Det}(U^T\left[\va_1,\va_2,\vq_3\right])\mathrm{Det}(U^T\left[\va_1',\va_2',\vq_3'\right])}{\|U^T\vq_3\|^3\|U^T\vq_3'\|^3}e^{-\mathrm{Tr}\left(U^TU\right)/2}dU\nonumber\\
&=\frac{1}{\|\va_3\|^2\|\va_3'\|^2}I_{\langle\mathrm{Lk}^2\rangle}(\left[\va_1,\va_2,\vq_3\right],\left[\va_1',\va_2',\vq_3'\right])
\end{align}
Next, we expand the product $\mathrm{Det}(U^T\left[\va_1,\va_2,\vq_3\right])\mathrm{Det}(U^T\left[\va_1',\va_2',\vq_3'\right])$ using the Cauchy-Binet formula \cite{Binet}, \cite{Cauchy}. To do so, we use the notation given in \cite{wiki:cbf} and let $[n]$ denote the set $\{1,2,...,n\}$, and let $\binom{[n]}{3}$ denote the set of 3-combinations of $[n]$. Moreover, for a given $3$ element subset $S=\{s_{i_1},s_{i_2},s_{i_3}\}$ of $[n]$ we let $U^T_{[3],S}$ denote the $3$-by-$3$ submatrix of $U^T$ whose columns are the $s_{i_1}^{th}$,$s_{i_2}^{th}$, and $s_{i_3}^{th}$ columns of $U^T$, and finally where $\left[\va_1,\va_2,\vq_3\right]_{S,[3]}$ denotes the $3$-by-$3$ submatrix of $\left[\va_1,\va_2,\vq_3\right]$ whose rows are the $s_{i_1}^{th}$,$s_{i_2}^{th}$, and $s_{i_3}^{th}$ rows of $\left[\va_1,\va_2,\vq_3\right]$. With the notation set up, we expand $\mathrm{Det}(U^T\left[\va_1,\va_2,\vq_3\right])\mathrm{Det}(U^T\left[\va_1',\va_2',\vq_3'\right])$ using the Cauchy-Binet formula:
\begin{align}\label{cbEXPAND}
&I_{\langle\mathrm{Lk}^2\rangle}(\left[\va_1,\va_2,\vq_3\right],\left[\va_1',\va_2',\vq_3'\right])=\nonumber\\
&\frac{1}{(2\pi)^{\frac{3n}{2}}}\int\displaylimits_{U\in\mathrm{Mat}_{n,3}}\frac{\sum\displaylimits_{S\in\binom{[n]}{3}}\mathrm{Det}\left(U^T_{[n],S}\right)\mathrm{Det}\left(\left[\va_1,\va_2,\vq_3\right]_{S,[n]}\right)\cdot\sum\displaylimits_{S\in\binom{[n]}{3}}\mathrm{Det}\left(U^T_{[n],S}\right)\mathrm{Det}\left(\left[\va_1',\va_2',\vq_3'\right]_{S,[n]}\right)}{\|U^T\vq_3\|^3\|U^T\vq_3'\|^3}dU
\end{align}
To simplify the denominator of the integrand, we define a basis, $B$, of $\mathbb{R}^n$ given by:
\begin{equation}\label{specialbasis}
\vb_1=\frac{\vq_3+\vq'_3}{\| \vq_3+\vq'_3 \|}\textrm{ , }
\vb_2=\frac{\vq_3-\vq'_3}{\| \vq_3-\vq'_3 \|}\textrm{ and }\vb_3\textrm{,}\vb_4,...,\vb_n\in\mathrm{span}\{\vq_3,\vq_3'\}^\perp,
\end{equation}
and also define the one parameter family:
\begin{align*}
a&=\vq_3\cdot\vb_1=\cos(\phi/2)\\
b&=\vq_3\cdot\vb_2=\sin(\phi/2)\\
\end{align*}
where $\phi$ is the angle between $\vq_3$ and $\vq_3'$. Notice also that $\vq_3'\cdot\vb_1=a$ and $\vq_3'\cdot\vb_1=-b$, and that in this basis $\vq_3=a\ve_1+b\ve_2$ and $\vq_3'=a\ve_1-b\ve_2$. With this choice of basis, the decomposition of $\mathrm{Det}(U^T\left[\va_1,\va_2,\vq_3\right])$ appearing in \eqref{cbEXPAND} can be simplified significantly:
\begin{align*}\nonumber
&\sum\displaylimits_{S\in\binom{[n]}{3}}\mathrm{Det}\left(U^T_{[n],S}\right)\mathrm{Det}\left(\left[\va_1,\va_2,\vq_3\right]_{S,[n]}\right)=\\
&\sum\displaylimits_{k=3}^n\mathrm{Det}\begin{pmatrix}
    u_{11}      & u_{12} & u_{1k} \\
    u_{21}       & u_{22} & u_{2k} \\
    u_{31}      & u_{32} & u_{3k}   \\
\end{pmatrix}\mathrm{Det}\begin{pmatrix}
    \va_1\cdot\vb_1      & \va_2\cdot\vb_1 & a \\
    \va_1\cdot\vb_2       & \va_2\cdot\vb_2 & b \\
    \va_1\cdot\vb_k  & \va_2\cdot\vb_k & 0 \\
\end{pmatrix}\\
&+\sum\displaylimits_{3\leq j<k\leq n}
\mathrm{Det}\begin{pmatrix}
u_{1j} & u_{1k} &  au_{11}+bu_{12}   \\
u_{2j} & u_{2k} &  au_{21}+bu_{22}  \\
u_{3j} & u_{3k}  &    au_{31}+bu_{32}  \\
\end{pmatrix}
\mathrm{Det}\begin{pmatrix}
    \va_1\cdot\vb_j       & \va_2\cdot\vb_j \\
    \va_1\cdot\vb_k  & \va_2\cdot\vb_k \\
\end{pmatrix},
\end{align*}
and a similar identity can be written for  $\mathrm{Det}(U^T\left[\va_1',\va_2',\vq_3'\right])$:
\begin{align*}\nonumber
&\sum\displaylimits_{S\in\binom{[n]}{3}}\mathrm{Det}\left(U^T_{[n],S}\right)\mathrm{Det}\left(\left[\va_1',\va_2',\vq_3'\right]_{S,[n]}\right)=\\
&\sum\displaylimits_{k=3}^n\mathrm{Det}\begin{pmatrix}
    u_{11}      & u_{12} & u_{1k} \\
    u_{21}       & u_{22} & u_{2k} \\
    u_{31}      & u_{32} & u_{3k}   \\
\end{pmatrix}\mathrm{Det}\begin{pmatrix}
    \va_1'\cdot\vb_1      & \va_2'\cdot\vb_1 & a \\
    \va_1'\cdot\vb_2       & \va_2'\cdot\vb_2 & -b \\
    \va_1'\cdot\vb_k  & \va_2'\cdot\vb_k & 0 \\
\end{pmatrix}\\
&+\sum\displaylimits_{3\leq j<k\leq n}
\mathrm{Det}\begin{pmatrix}
u_{1j} & u_{1k} &  au_{11}-bu_{12}   \\
u_{2j} & u_{2k} &  au_{21}-bu_{22}  \\
u_{3j} & u_{3k}  &    au_{31}-bu_{32}  \\
\end{pmatrix}
\mathrm{Det}\begin{pmatrix}
    \va_1'\cdot\vb_j       & \va_2'\cdot\vb_j \\
    \va_1'\cdot\vb_k  & \va_2'\cdot\vb_k \\
\end{pmatrix},
\end{align*}
In the previous identities, the dependence on the basis $B$ has been suppressed in writing the matrix elements of $U$ since the Gaussian measure is invariant under $O(N)$. From here, using the decompositions above along with numerous symmetries of the Gaussian measure, we write $I_{\langle\mathrm{Lk}^2\rangle}(A,A')$ in the basis \eqref{specialbasis}:
\begin{align}\label{cbEXPAND_final}
&I_{\langle\mathrm{Lk}^2\rangle}(\left[\va_1,\va_2,\vq_3\right],\left[\va_1',\va_2',\vq_3'\right])=\nonumber\\
&-\frac{I_{\langle\mathrm{Lk}^2\rangle}(\left[\ve_1,\ve_3,a\ve_1+b\ve_2\right],\left[\ve_1,\ve_3,a\ve_1-b\ve_2\right])}{b^2}\sum\displaylimits_{k=3}^n\mathrm{Det}\begin{pmatrix}
    \va_1'\cdot\vb_1      & \va_2'\cdot\vb_1 & a \\
    \va_1'\cdot\vb_2       & \va_2'\cdot\vb_2 & -b \\
    \va_1'\cdot\vb_k  & \va_2'\cdot\vb_k & 0 \\
\end{pmatrix}\mathrm{Det}\begin{pmatrix}
    \va_1\cdot\vb_1      & \va_2\cdot\vb_1 & a \\
    \va_1\cdot\vb_2       & \va_2\cdot\vb_2 & b \\
    \va_1\cdot\vb_k  & \va_2\cdot\vb_k & 0 \\
\end{pmatrix}\nonumber\\
&+I_{\langle\mathrm{Lk}^2\rangle}(\left[\ve_3,\ve_4,a\ve_1+b\ve_2\right],\left[\ve_3,\ve_4,a\ve_1-b\ve_2\right])\sum\displaylimits_{3\leq j<k\leq n}^n\mathrm{Det}\begin{pmatrix}
    \va_1'\cdot\vb_j       & \va_2'\cdot\vb_j \\
    \va_1'\cdot\vb_k  & \va_2'\cdot\vb_k \\
\end{pmatrix}\mathrm{Det}\begin{pmatrix}
    \va_1\cdot\vb_j       & \va_2\cdot\vb_j \\
    \va_1\cdot\vb_k  & \va_2\cdot\vb_k \\
\end{pmatrix}
\end{align}
where in \eqref{cbEXPAND_final} we have implicitly used the identities:
\begin{align*}
\frac{1}{(2\pi)^{\frac{3\cdot n}{2}}}\int\displaylimits_{U\in\mathrm{Mat}_{n,3}}&\frac{\mathrm{Det}\begin{pmatrix}
    u_{11}      & u_{12} & u_{13} \\
    u_{21}       & u_{22} & u_{23} \\
    u_{31}      & u_{32} & u_{33}   \\
\end{pmatrix}^2}{\|U^T(a\ve_1+b\ve_2)\|^3\|U^T(a\ve_1-b\ve_2)\|^3}dU\nonumber\\
=-&\frac{I_{\langle\mathrm{Lk}^2\rangle}(\left[\ve_1,\ve_3,a\ve_1+b\ve_2\right],\left[\ve_1,\ve_3,a\ve_1-b\ve_2\right])}{b^2}\nonumber\\
=&\frac{I_{\langle\mathrm{Lk}^2\rangle}(\left[\ve_2,\ve_3,a\ve_1+b\ve_2\right],\left[\ve_2,\ve_3,a\ve_1-b\ve_2\right])}{a^2}\nonumber\\
=-&\frac{I_{\langle\mathrm{Lk}^2\rangle}(\left[\ve_1,\ve_3,a\ve_1+b\ve_2\right],\left[\ve_2,\ve_3,a\ve_1-b\ve_2\right])}{ab}\nonumber\\
\end{align*}
The final step is to evaluate the two Gaussian integrals appearing in \eqref{cbEXPAND_final}, which for brevity's sake we will define the shorthand for:
\begin{eqnarray*}
I_{\langle\mathrm{Lk}^2\rangle}(\left[\ve_1,\ve_3,a\ve_1+b\ve_2\right],\left[\ve_1,\ve_3,a\ve_1-b\ve_2\right])=I_{\langle\mathrm{Lk}^2\rangle}^{1313}(a,b)\nonumber\\
I_{\langle\mathrm{Lk}^2\rangle}(\left[\ve_3,\ve_4,a\ve_1+b\ve_2\right],\left[\ve_3,\ve_4,a\ve_1-b\ve_2\right])=I_{\langle\mathrm{Lk}^2\rangle}^{3434}(a,b)\nonumber\\
\end{eqnarray*}
First, integrate over all variables appearing \textit{at most} in the numerator of the integrand of $I_{\langle\mathrm{Lk}^2\rangle}^{1313}(a,b)$ except for $u_{13}, u_{23}$ and $u_{33}$ to obtain:
\begin{eqnarray*}
I_{\langle\mathrm{Lk}^2\rangle}^{1313}(a,b)=\nonumber\\
&=\frac{-b^2}{(2\pi)^{9/2}}\int\displaylimits_{\mathbb{R}^9}\frac{Det\left(\begin{bmatrix}
    u_{11}       & u_{12}&  u_{13}\\
    u_{21}       & u_{22}& u_{23}\\
    u_{31}       & u_{32}& u_{33}\\
\end{bmatrix}\right)^2e^{-\frac{1}{2}\left(u_{11}^2+u_{21}^2+u_{31}^2+u_{12}^2+u_{22}^2+u_{32}^2+u_{13}^2+u_{23}^2+u_{33}^2\right)}}{\left((au_{11}+bu_{12})^2+(au_{21}+bu_{22})^2+(au_{31}+bu_{32})^2\right)^{3/2}\left((au_{11}-bu_{12})^2+(au_{21}-bu_{22})^2+(au_{31}-bu_{32})^2\right)^{3/2}}\nonumber\\
\end{eqnarray*}
Using the rotational invariance of the Gaussian measure, we can align the vector $\{u_{13},u_{23},u_{33}\}$ along the x-axis, and integrate to obtain:
\begin{eqnarray*}
I_{\langle\mathrm{Lk}^2\rangle}^{1313}(a,b)=\nonumber\\
&=\frac{-3b^2}{8\pi^3}\int\displaylimits_{\mathbb{R}^6}\frac{Det\left(\begin{bmatrix}
    u_{11}       & u_{12}\\
    u_{21}       & u_{22}\\
\end{bmatrix}\right)^2e^{-\frac{1}{2}\left(u_{11}^2+u_{21}^2+u_{31}^2+u_{12}^2+u_{22}^2+u_{32}^2\right)}}{\left((au_{11}+bu_{12})^2+(au_{21}+bu_{22})^2+(au_{31}+bu_{32})^2\right)^{3/2}\left((au_{11}-bu_{12})^2+(au_{21}-bu_{22})^2+(au_{31}-bu_{32})^2\right)^{3/2}}\nonumber\\
\end{eqnarray*}
Next, by a straightforward change of variables we obtain:
\begin{eqnarray*}
I_{\langle\mathrm{Lk}^2\rangle}^{1313}(a,b)=\nonumber\\
&=\frac{-3}{8\pi^3a^2|a^3b^3|}\int\displaylimits_{\mathbb{R}^6}\frac{Det\left(\begin{bmatrix}
    u_{11}       & u_{12}\\
    u_{21}       & u_{22}\\
\end{bmatrix}\right)^2e^{-\frac{1}{2}\left(\frac{u_{11}^2+u_{21}^2+u_{31}^2}{a^2}+\frac{u_{12}^2+u_{22}^2+u_{32}^2}{b^2}\right)}}{\left((u_{11}+u_{12})^2+(u_{21}+u_{22})^2+(u_{31}+u_{32})^2\right)^{3/2}\left((u_{11}-u_{12})^2+(u_{21}-u_{22})^2+(u_{31}-u_{32})^2\right)^{3/2}}\nonumber\\
\end{eqnarray*}
Now use the fact that:
\begin{equation}\nonumber
Det\left(\begin{bmatrix}
    u_{11}       & u_{21} & u_{31}\\
    u_{12}       & u_{22} & u_{32}\\
\end{bmatrix}\begin{bmatrix}
    u_{11}       & u_{12}\\
    u_{21}       & u_{22} \\
    u_{31}       & u_{32} \\
\end{bmatrix}\right)=Det\left(\begin{bmatrix}
    u_{11}       & u_{12}\\
    u_{21}       & u_{22}\\
\end{bmatrix}\right)^2+Det\left(\begin{bmatrix}
    u_{11}       & u_{12}\\
    u_{31}       & u_{32}\\
\end{bmatrix}\right)^2+Det\left(\begin{bmatrix}
    u_{21}       & u_{22}\\
    u_{31}       & u_{32}\\
\end{bmatrix}\right)^2,
\end{equation}
along with the symmetry\footnote{Notice that the factor of 3 multiplying the integral has vanished} of the Gaussian measure to obtain that:
\begin{eqnarray}\label{1313_int}
I_{\langle\mathrm{Lk}^2\rangle}^{1313}(a,b)=\nonumber\\
&=\frac{-1}{8\pi^3a^2|a^3b^3|}\int\displaylimits_{\mathbb{R}^6}\frac{Det\left(\begin{bmatrix}
    u_{11}       & u_{21} & u_{31}\\
    u_{12}       & u_{22} & u_{32}\\
\end{bmatrix}\begin{bmatrix}
    u_{11}       & u_{12}\\
    u_{21}       & u_{22} \\
    u_{31}       & u_{32} \\
\end{bmatrix}\right)e^{-\frac{1}{2}\left(\frac{u_{11}^2+u_{21}^2+u_{31}^2}{a^2}+\frac{u_{12}^2+u_{22}^2+u_{32}^2}{b^2}\right)}}{\left((u_{11}+u_{12})^2+(u_{21}+u_{22})^2+(u_{31}+u_{32})^2\right)^{3/2}\left((u_{11}-u_{12})^2+(u_{21}-u_{22})^2+(u_{31}-u_{32})^2\right)^{3/2}}\nonumber\\
&=\frac{-1}{8\pi^3a^2|a^3b^3|}\int\displaylimits_{\mathbb{R}^6}\frac{\left((u_{11}^2+u_{21}^2+u_{31}^2)(u_{12}^2+u_{22}^2+u_{32}^2)-(u_{11}u_{12}+u_{21}u_{22}+u_{31}u_{32})^2\right)e^{-\frac{1}{2}\left(\frac{u_{11}^2+u_{21}^2+u_{31}^2}{a^2}+\frac{u_{12}^2+u_{22}^2+u_{32}^2}{b^2}\right)}}{\left((u_{11}+u_{12})^2+(u_{21}+u_{22})^2+(u_{31}+u_{32})^2\right)^{3/2}\left((u_{11}-u_{12})^2+(u_{21}-u_{22})^2+(u_{31}-u_{32})^2\right)^{3/2}}\nonumber\\
\end{eqnarray}
Now define that:
\begin{align*}
\vU_1&=\{u_{11},u_{21},u_{31}\}\nonumber\\
\vU_2&=\{u_{12},u_{22},u_{32}\}\nonumber\\
\end{align*}
so that \eqref{1313_int} the above can be simplified to:
\begin{eqnarray*}
I_{\langle\mathrm{Lk}^2\rangle}^{1313}(a,b)=\nonumber\\
&=\frac{-1}{8\pi^3a^2|a^3b^3|}\int\displaylimits_{\mathbb{R}^3\times\mathbb{R}^3}\frac{\left(\|\vU_1\|^2\|\vU_2\|^2-(\vU_1\cdot\vU_2)^2\right)e^{-\frac{1}{2}\left(\|\vU_1\|^2+\|\vU_2\|^2\right)}}{\left((\|\vU_1\|^2+\|\vU_2\|^2)^2-4\cdot(\vU_1\cdot\vU_2)^2\right)^{3/2}}d\vU_1d\vU_2\nonumber\\
\end{eqnarray*}
Next write $\mathbb{R}^6=\mathbb{R}^3\times\mathbb{R}^3$, and change to spherical coordinates in each factor:
\begin{align*}\nonumber
u_{11}&=R_1\sin(\theta_1)\cos(\phi_1)\nonumber\\
u_{21}&=R_1\sin(\theta_1)\sin(\phi_1)\nonumber\\
u_{31}&=R_1\cos(\theta_1)\nonumber\\
u_{12}&=R_2\sin(\theta_2)\cos(\phi_2)\nonumber\\
u_{22}&=R_2\sin(\theta_2)\sin(\phi_2)\nonumber\\
u_{32}&=R_2\cos(\theta_2)\nonumber\\
\end{align*}
to obtain:
\begin{eqnarray*}
I_{\langle\mathrm{Lk}^2\rangle}^{1313}(a,b)=\nonumber\\
&=\frac{-1}{8\pi^3a^2|a^3b^3|}\int\displaylimits_{\mathbb{R}^6}\frac{R_1^2R_2^2\left(1-\cos(\psi)^2\right)e^{-\frac{1}{2}\left(\frac{R_1^2}{a^2}+\frac{R_2^2}{b^2}\right)}}{\left((R_1^2+R_2^2)^2-4R_1^2R_2^2\cos(\psi)^2\right)^{3/2}}R_1^2R_2^2\sin(\theta_1)\sin(\theta_2)\nonumber\\
\end{eqnarray*}
where 
\begin{equation}\nonumber
\cos(\psi)=\left(\sin(\theta_1) \sin (\theta_2) \cos (\phi_1-\phi_2)+\cos (\theta_1) \cos(\theta_2)\right)\\\nonumber
\end{equation}
is the cosine of the angle between $\frac{\vU_1}{\|\vU_1\|}$ and $\frac{\vU_2}{\|\vU_2\|}$. Lastly, convert $R_1$ and $R_2$ to polar coordinates and integrate over the resulting radial coordinate to obtain:
\begin{equation}\label{penultimate}
\frac{a^2 b^4}{64\pi ^3\left| ab\right|^3}\int\displaylimits_{T=0}^{\frac{\pi}{2}}\int\displaylimits_{\mathbb{S}^2\times\mathbb{S}^2}\frac{\left(1-\cos ^2(\psi )\right) \sin ^4(2 T)}{\left(a^2 \sin ^2(T)+b^2 \cos ^2(T)\right)^2\left(1-\cos ^2(\psi) \sin ^2(2 T)\right)^{3/2}}\sin(\theta_1)\sin(\theta_2)
\end{equation}
There is a rotational symmetry here that allows us to easily integrate over one of the 2-sphere factors. Lastly, we integrate over $\phi_1$, $\theta_1$ and lastly $T$ using a symbolic integration in Wolfram Mathematica \cite{Mathematica} to get:
\begin{eqnarray}\label{1313_final}
I_{\langle\mathrm{Lk}^2\rangle}^{1313}(a,b)=\nonumber\\
&=-\frac{2}{\pi}b^2\frac{(a^2-b^2)+|ab|\left(\pi-4\tan^{-1}(|a/b|)\right)}{(a^2-b^2)^3|a||b|}\nonumber\\
&=-\frac{2}{\pi}b^2\frac{(a^2-b^2)-2|ab|\sin^{-1}(a^2-b^2)}{(a^2-b^2)^3|a||b|}\nonumber\\
\end{eqnarray}
Using similar, and in fact simpler arguments, we obtain
\begin{equation}\nonumber
I_{\langle\mathrm{Lk}^2\rangle}^{3434}(a,b)=\frac{2}{\pi}\frac{(a^2-b^2)-2|ab|\sin^{-1}(a^2-b^2)}{(a^2-b^2)^2|a||b|}
\end{equation}
From which it follows that:
\begin{equation}\label{int_relation}
I_{\langle\mathrm{Lk}^2\rangle}^{1313}(a,b)=\frac{b^2}{b^2-a^2}I_{\langle\mathrm{Lk}^2\rangle}^{3434}(a,b)
\end{equation}
Now substitute \eqref{1313_final} and \eqref{int_relation} into \eqref{cbEXPAND_final}:
\begin{align}
&\frac{a^2-b^2}{I_{\langle\mathrm{Lk}^2\rangle}^{3434}(a,b)}I_{\langle\mathrm{Lk}^2\rangle}(\left[\va_1,\va_2,\vq_3\right],\left[\va_1',\va_2',\vq_3'\right])=\nonumber\\
&\sum\displaylimits_{k=3}^n\mathrm{Det}\begin{pmatrix}
    \va_1'\cdot\vb_1      & \va_2'\cdot\vb_1 & a \\
    \va_1'\cdot\vb_2       & \va_2'\cdot\vb_2 & -b \\
    \va_1'\cdot\vb_k  & \va_2'\cdot\vb_k & 0 \\
\end{pmatrix}\mathrm{Det}\begin{pmatrix}
    \va_1\cdot\vb_1      & \va_2\cdot\vb_1 & a \\
    \va_1\cdot\vb_2       & \va_2\cdot\vb_2 & b \\
    \va_1\cdot\vb_k  & \va_2\cdot\vb_k & 0 \\
\end{pmatrix}\nonumber\\
&+(a^2-b^2)\sum\displaylimits_{3\leq j<k\leq n}^n\mathrm{Det}\begin{pmatrix}
    \va_1'\cdot\vb_j       & \va_2'\cdot\vb_j \\
    \va_1'\cdot\vb_k  & \va_2'\cdot\vb_k \\
\end{pmatrix}\mathrm{Det}\begin{pmatrix}
    \va_1\cdot\vb_j       & \va_2\cdot\vb_j \\
    \va_1\cdot\vb_k  & \va_2\cdot\vb_k \\
\end{pmatrix}\nonumber\\
&=\mathrm{Det}\begin{pmatrix}
    \va_1\cdot\va_{1}'      & \va_1\cdot\va_{2}' & \va_1\cdot\vq_{3}' \\
    \va_2\cdot\va_{1}'      & \va_2\cdot\va_{2}' & \va_2\cdot\vq_{3}'\\
    \vq_3\cdot\va_{1}'      & \vq_3\cdot\va_{2}' & \vq_3\cdot\vq_{3}' \\
\end{pmatrix}\nonumber\\
\end{align}
Rearranging, all of the above may be summarized as:
\begin{equation}\nonumber
I_{\langle\mathrm{Lk}^2\rangle}(\left[\va_1,\va_2,\vq_3\right],\left[\va_1',\va_2',\vq_3'\right])=\frac{I_{\langle\mathrm{Lk}^2\rangle}^{3434}(a,b)}{a^2-b^2}\mathrm{Det}\begin{pmatrix}
    \va_1\cdot\va_{1}'      & \va_1\cdot\va_{2}' & \va_1\cdot\vq_{3}' \\
    \va_2\cdot\va_{1}'      & \va_2\cdot\va_{2}' & \va_2\cdot\vq_{3}'\\
    \vq_3\cdot\va_{1}'      & \vq_3\cdot\va_{2}' & \vq_3\cdot\vq_{3}' \\
\end{pmatrix},
\end{equation}
and the result follows by rewriting any instances of $a$, $b$ in terms of dot products between $\va_3$ and $\va_3'$, and then substituting into \eqref{unnorm_to_norm}.
\end{Proof}
      
\section[Integration Over the Stiefel Manifold]{Integration Over the Stiefel Manifold}
        \label{sec:LabelForChapter3:Section3}
In the previous section we averaged the value of the square of the linking number over all projections of a pair of space curves in $\mathbb{R}^n$ to $3$ dimensional subspaces, to obtain the result \eqref{GaussianResult} which is a function only of the inner products of the initial data. The main theorem of this paper is the following:
\begin{RestateTheorem}{
Given two closed, differentiable, non-intersecting space curves $\gamma_1,\gamma_2:\mathbb{S}^1\rightarrow\mathbb{R}^n$ with parametrizations $\vr_1(t)$ and $\vr_2(s)$, and let $\mathrm{proj}(\gamma_1),\mathrm{proj}(\gamma_2)$ denote the orthogonal projections of $\gamma_1$ and $\gamma_2$ on to the $3$ dimensional subspace spanned by the columns of $U$, then the expected value of $\mathrm{Lk}^2(\mathrm{proj}(\gamma_1),\mathrm{proj}(\gamma_2))$, denoted $\langle\mathrm{Lk}^2\rangle$, averaged over all orthogonal projections to $3$ dimensional subspaces with respect to the unique normalized $O(n)$-invariant measure on $Gr(n,3)$, is given by the following integral:
\begin{equation}\label{theRESULT}
\langle\mathrm{Lk}^2\rangle=\frac{1}{16\pi^2}\int\displaylimits_{(s,t)\in\mathbb{T}^2}\int\displaylimits_
{(s',t')\in\mathbb{T}^2}I_{\langle\mathrm{Lk}^2\rangle}(A(s,t),A'(s',t'))dsdtds'dt'
\end{equation}
where
\begin{align*}\label{eq:result3THM}
&\frac{\pi}{4}I_{\langle\mathrm{Lk}^2\rangle}(A,A')=\nonumber\\
&\frac{(\va_3\cdot\va_3')-\sqrt{\|\va_3\|^2\|\va_3'\|^2-(\va_3\cdot\va_3')^2}\sin^{-1}\left(\frac{\va_3\cdot\va_3'}{\|\va_3\|\|\va_3'\|}\right)}{(\va_3\cdot\va_3')^3\sqrt{\|\va_3\|^2\|\va_3'\|^2-(\va_3\cdot\va_3')^2}}\mathrm{Det}\left(\begin{bmatrix}
    \va_1\cdot\va_1'       & \va_1\cdot\va_2' & \va_1\cdot\va_3'   \\
    \va_2\cdot\va_1'       & \va_2\cdot\va_2' & \va_2\cdot\va_3'   \\
    \va_3\cdot\va_1'       & \va_3\cdot\va_2' & \va_3\cdot\va_3'   \\
\end{bmatrix}\right)
\end{align*} 
and
\begin{align*}\nonumber
\va_1&=\dot{\vr}_1(t)\mathrm{, }\va_2=\dot{\vr}_2(s)\mathrm{, }\va_3=\vr_2(s)-\vr_1(t)\\
\va_1'&=\dot{\vr}_1(t')\mathrm{, }\va_2'=\dot{\vr}_2(s')\mathrm{, }\va_3'=\vr_2(s')-\vr_1(t')\\
A&=\big[\dot{\vr}_1(t),\dot{\vr}_2(s),\vr_2(s)-\vr_1(t)\big]\\
A'&=\big[\dot{\vr}_1(t'),\dot{\vr}_2(s'),\vr_2(s')-\vr_1(t')\big]\\
\end{align*}}{mainResultTheorem}
\end{RestateTheorem}
Before proving this, we first introduce the Stiefel manifold $V_{n,m}$ of orthonormal $m$-frames in $\mathbb{R}^n$ (in what follows, we will use the normalization conventions in \cite{Radon} and \cite{Zhang}):
\begin{Definition}
For $n\geq m$, the Stiefel manifold, $V_{n,m}$ of orthonormal $m$-frames in $\mathbb{R}^n$ is the set of matrices $M\in\mathrm{Mat}_{n,m}$ such that $M^TM=\mathds{1}_m$. 
\end{Definition}
There is a projection map, $\pi$, from $V_{n,m}$ to $Gr(n,m)$ where the set of vectors, $M$, is mapped to $\mathrm{span}\{M\}$, and moreover any right $O(m)$ invariant function $f(x)$ on $V_{n,m}$ gives a function $F(\pi(x))=f(x)$ on $Gr(n,m)$. Now we will fix the unique normalized $O(n)$ invariant measures $dv$ on $V_{n,k}$ and $dy$ on $Gr(n,k)$, the former which is normalized by:
\begin{equation}\nonumber
\sigma_{n,k}=\int\displaylimits_{V_{n,k}}dv=\frac{\pi^{kn/2}}{\Gamma_k(n/2)}
\end{equation}
where $\Gamma_k(n/2)=(2\pi)^{(k^2-k)/4}\prod\limits_{i=1}^k\Gamma(n/2-\frac{1}{2}(i-1))$ and $\Gamma(x)=\int\displaylimits_0^\infty t^{x-1}e^{-t}dt$ is the gamma function, and the latter, $dy$, which is normalized by stipulating that:
\begin{equation}\label{VtoG}
\int\displaylimits_{V_{n,m}}f(x)dx=\int\displaylimits_{Gr(n,m)}F(y)dy.
\end{equation}
We will now quote a result concerning a decomposition of the measure $dx$ on $\mathrm{Mat}_{n,k}$ (induced by the polar decomposition of the matrix $x$), which is proved in \cite{Radon} and \cite{Zhang}. As a set up, identify $\Omega_{k,k}$ with the set of symmetric positive-definite $k$ by $k$ matrices viewed as a subspace of $\mathbb{R}^{k+1 \choose 2}$. With this identification we have a measure on $\Omega_{k,k}$ given by $d\widetilde{r}(r)=\mathrm{Det}(r)^{-(k+1)/2}dr$ where $dr$ is the Lebesgue measure. Now we will invoke the real version of lemma 3.1 in G. Zhang's paper:\cite{Zhang}.\\ 
\textit{\textbf{Lemma 3.1 in \cite{Zhang}:} Let $dx$, $dv$ and $d\widetilde{r}(r)$ be the normalized measures on $\mathrm{Mat}_{n,k}$, $V_{n,k}$ and $\Omega_{k,k}$ defined above. Almost every $x\in \mathrm{Mat}_{n,k}$ may be decomposed as:
\begin{equation}\nonumber
x=vr^{1/2}
\end{equation} where $x\in V_{n,k}$ and $r\in \Omega_{k,k}$. Under this decomposition the measure $dx$ is given by:
\begin{equation}\nonumber
dx=C_0\mathrm{Det}(r)^{n/2}dvd\widetilde{r}(r)=C_0\mathrm{Det}(r)^{(n-k-1)/2}dvdr,
\end{equation}
namely, 
\begin{equation}\nonumber
\int\displaylimits_{\mathrm{Mat}_{n,k}}f(x)dx=C_0\int\displaylimits_{\mathrm{V}_{n,k}}\int\displaylimits_{\Omega_{k,k}}f(vr^{1/2})dvdr,
\end{equation}
where
\begin{equation}\nonumber
C_0=\frac{\pi^{nk/2}}{\Gamma_k(n/2)}
\end{equation}
}\\
With this setup, we are now ready to prove Theorem \ref{mainResultTheorem}.  
\begin{Proof}
The expectation that we are computing may be summarized as:
\begin{equation}\nonumber
\frac{1}{(\pi)^{3\cdot n/2}}\int\displaylimits_{U\in\mathrm{Mat}_{n,3}}\left(\frac{1}{(4\pi)^2}\int\displaylimits_{\mathbb{T}^2}\int\displaylimits_{\mathbb{T}^2}\frac{\mathrm{Det}(U^TA(s,t))\mathrm{Det}(U^TA(s',t'))}{\|U^t\va_3(s,t)\|^3\|U^t\va_3'(s',t')\|^3}e^{-\mathrm{Tr}U^TU}d\mu d\mu'\right)e^{-\mathrm{Tr}U^TU}dU
\end{equation}
where in the above we have $d\mu=dsdt$ and $d\mu'=ds'dt'$. Using the decomposition in the previous lemma, we obtain:
\begin{align*}\nonumber
&\frac{C_0}{(\pi)^{3\cdot n/2}}\int\displaylimits_{v\in V_{n,3}}\int\displaylimits_{r\in\Omega_{3,3}}\left(\frac{1}{(4\pi)^2}\int\displaylimits_{\mathbb{T}^2}\int\displaylimits_{\mathbb{T}^2}\frac{\mathrm{Det}((r^{1/2})^Tv^TA(s,t))\mathrm{Det}((r^{1/2})^Tv^TA(s',t'))}{\|(r^{1/2})^Tv^T\va_3(s,t)\|^3\|(r^{1/2})^Tv^T\va_3'(s,t)\|^3}e^{-\mathrm{Tr}U^TU}d\mu d\mu'\right)dvd\widetilde{r}(r)\nonumber\\
&\frac{1}{\Gamma_3(n/2)}\int\displaylimits_{v\in V_{n,3}}\left(\frac{1}{(4\pi)^2}\int\displaylimits_{(s,t)\in\mathbb{T}^2}\int\displaylimits_{(s',t')\in\mathbb{T}^2}\frac{\mathrm{Det}(v^TA(s,t))\mathrm{Det}(v^TA(s',t'))}{\|v^T\va_3(s,t)\|^3\|v^T\va_3'(s,t)\|^3}e^{-\mathrm{Tr}(r)}d\mu d\mu'\right)dvd\widetilde{r}(r)
\end{align*}
In the above, the dependence on $r$ in the integrand was removed using the fact that the linking number is invariant under ambient homeomorphisms of $\mathbb{R}^3$. Lastly, the result follows by integrating over $\Omega_{3,3}$ using equation (2.2) in \cite{Zhang} and then applying Lemma \ref{codim1lemma}. 
\end{Proof}
        \section[Integration Over the Configuration Space]{Integration Over the Configuration Space}
        \label{sec:LabelForChapter3:Section4}            
Using (\ref{codim1lemma}) we may find some conditions on the initial pair of space curves so as to bound the second moment of the linking number. We will see through a few examples later that these are reasonable conditions. As in Theorem (\ref{mainResultTheorem}), set:
\begin{align*}\nonumber
\va_1(t)&=\dot{\vr}_1(t)\textrm{, }\va_2(s)=\dot{\vr}_2(s)\textrm{, }\va_3(s,t)=\vr_2(s)-\vr_1(t)\\
\va_1'(t')&=\dot{\vr}_1(t')\textrm{, }\va_2'(s')=\dot{\vr}_2(s')\textrm{, }\va_3'(s',t')=\vr_2(s')-\vr_1(t')\\
A(s,t)&=\big[\dot{\vr}_1(t),\dot{\vr}_2(s),\vr_2(s)-\vr_1(t)\big]\\
A'(s',t')&=\big[\dot{\vr}_1(t'),\dot{\vr}_2(s'),\vr_2(s')-\vr_1(t')\big]
\end{align*}
\begin{RestateTheorem}{
With the definitions above, let $v_1=\textrm{max}_{t\in\mathbb{S}^1}\|\dot{\vr}_1(t)\|$, $v_2=\textrm{max}_{s\in\mathbb{S}^1}\|\dot{\vr}_2(s)\|$, $k=\mathrm{min}_{(s,t)\in\mathbb{T}^2}\|\vr_2(t)-\vr_1(s)\|$, and 
\begin{equation}\nonumber
C=\int\displaylimits_{(s,t,s',t')\in\mathbb{T}^2\times\mathbb{T}^2}\frac{1}{\sqrt{\|\va_3(s,t)\|^2\|\va_3'(s',t')\|^2-(\va_3(s,t)\cdot\va_3'(s',t'))^2}}dsdtds'dt',
\end{equation}
then
\begin{equation}\nonumber
\langle\mathrm{Lk}^2\rangle\leq\frac{1}{(4\pi)^2}\frac{4Cv_1^2v_2^2}{\pi k^2}
\end{equation}}{Lk2Bound}
\end{RestateTheorem}
 \begin{Proof}
The relevant configuration space to integrate over will be the set of pairs of points, one pair per component of the link. Since
\begin{eqnarray}\nonumber
\langle\mathrm{Lk}^2\rangle=\frac{1}{(4\pi)^2}\int\displaylimits_{(s,t,s',t')\in\mathbb{T}^2\times\mathbb{T}^2}I_{\langle\mathrm{Lk}^2\rangle}(A(s,t),A'(s',t'))dsdtds'dt'
\end{eqnarray}
we may find an upper bound on $\langle\mathrm{Lk}^2\rangle$ by first bounding the function $I_{\langle\mathrm{Lk}^2\rangle}(A,A')$ and then computing the integral over the configuration space. Note that in deriving a bound on $I_{\langle\mathrm{Lk}^2\rangle}(A(s,t),A(s',t'))$, we suppress the dependence on $s$, $s'$, $t$ and $t'$, since the parametrizations of the components of the link have no bearing on the bounds we derive. Now, recall from \eqref{GaussianResult} that:
\begin{align*}
&\frac{\pi}{4}I_{\langle\mathrm{Lk}^2\rangle}(A,A')=\nonumber\\
&\frac{(\va_3\cdot\va_3')-\sqrt{\|\va_3\|^2\|\va_3'\|^2-(\va_3\cdot\va_3')^2}\sin^{-1}\left(\frac{\va_3\cdot\va_3'}{\|\va_3\|\|\va_3'\|}\right)}{(\va_3\cdot\va_3')^3\sqrt{\|\va_3\|^2\|\va_3'\|^2-(\va_3\cdot\va_3')^2}}\mathrm{Det}\left(\begin{bmatrix}
    \va_1\cdot\va_1'       & \va_1\cdot\va_2' & \va_1\cdot\va_3'   \\
    \va_2\cdot\va_1'       & \va_2\cdot\va_2' & \va_2\cdot\va_3'   \\
    \va_3\cdot\va_1'       & \va_3\cdot\va_2' & \va_3\cdot\va_3'   \\
\end{bmatrix}\right)\nonumber
&\nonumber
\end{align*}
which can be rewritten as:
\begin{align}\label{GaussianResult_rewritten}
&I_{\langle\mathrm{Lk}^2\rangle}(A,A')=\nonumber\\
&\frac{4}{\pi}\left(\frac{1}{(\va_3\cdot\va_3')^2\sqrt{\|\va_3\|^2\|\va_3'\|^2-(\va_3\cdot\va_3')^2}}-\frac{\sin^{-1}\left(\frac{\va_3\cdot\va_3'}{\|\va_3\|\|\va_3'\|}\right)}{(\va_3\cdot\va_3')^3}\right)\mathrm{Det}\left(\begin{bmatrix}
    \va_1\cdot\va_1'       & \va_1\cdot\va_2' & \va_1\cdot\va_3'   \\
    \va_2\cdot\va_1'       & \va_2\cdot\va_2' & \va_2\cdot\va_3'   \\
    \va_3\cdot\va_1'       & \va_3\cdot\va_2' & \va_3\cdot\va_3'   \\
\end{bmatrix}\right)
&
\end{align}
In \eqref{GaussianResult_rewritten} we notice that:
\begin{equation}\nonumber
\frac{1}{(\va_3\cdot\va_3')^2\sqrt{\|\va_3\|^2\|\va_3'\|^2-(\va_3\cdot\va_3')^2}}\geq\frac{\sin^{-1}\left(\frac{\va_3\cdot\va_3'}{\|\va_3\|\|\va_3'\|}\right)}{(\va_3\cdot\va_3')^3}\geq 0
\end{equation}
and:
\begin{equation}\nonumber
|\mathrm{Det}\left(\begin{bmatrix}
    \va_1\cdot\va_1'       & \va_1\cdot\va_2' & \va_1\cdot\va_3'   \\
    \va_2\cdot\va_1'       & \va_2\cdot\va_2' & \va_2\cdot\va_3'   \\
    \va_3\cdot\va_1'       & \va_3\cdot\va_2' & \va_3\cdot\va_3'   \\
\end{bmatrix}\right)|\leq\prod\displaylimits_{i=1}^3\|\va_i\|\|\va_i'\|
\end{equation}
In this way, we have an upper bound for 
$I_{\langle\mathrm{Lk}^2\rangle}(A,A')$ which may be expressed as:
\begin{align}\label{starthere}
|I_{\langle\mathrm{Lk}^2\rangle}(A(s,t),A(s',t'))|\leq\frac{4}{\pi}\left(\frac{(\va_3\cdot\va_3')-\sqrt{\|\va_3\|^2\|\va_3'\|^2-(\va_3\cdot\va_3')^2}\sin^{-1}\left(\frac{\va_3\cdot\va_3'}{\|\va_3\|\|\va_3'\|}\right)}{(\va_3\cdot\va_3')^3\sqrt{\|\va_3\|^2\|\va_3'\|^2-(\va_3\cdot\va_3')^2}}\right)\prod\displaylimits_{i=1}^3\|\va_i\|\|\va_i'\|.
\end{align}
In the above, we see two types of possible singularities, one which seems to occur when $\va_3\cdot\va_3'=0$ and another which occurs when $\|\va_3\|^2\|\va_3'\|^2-(\va_3\cdot\va_3')^2=0$. In the former case, we note that in \eqref{starthere}, that both the numerator and denominator of the term appearing in the parentheses vanish to order $3$ in the angle between $\va_3$ and $\va_3'$ whenever $\va_3$ and $\va_3'$ become orthogonal. In addition, we also have a bound:
\begin{equation}\nonumber
\frac{1}{3\cdot\|\va_3\|^2\|\va_3'\|^2} \leq \frac{(\va_3\cdot\va_3')-\sqrt{\|\va_3\|^2\|\va_3'\|^2-(\va_3\cdot\va_3')^2}\sin^{-1}\left(\frac{\va_3\cdot\va_3'}{\|\va_3\|\|\va_3'\|}\right)}{(\va_3\cdot\va_3')^3}\leq\frac{1}{\|\va_3\|^2\|\va_3'\|^2},
\end{equation}
from which it follows that:
\begin{align*}\nonumber
|I_{\langle\mathrm{Lk}^2\rangle}(A,A')|&\leq\frac{4}{\pi}\frac{1}{\|\va_3\|\|\va_3'\|}\frac{\prod\displaylimits_{i=1}^2\|\va_i\|\|\va_i'\|}{\sqrt{\|\va_3\|^2\|\va_3'\|^2-(\va_3\cdot\va_3')^2}}\nonumber\\
&\leq\frac{4}{\pi}\cdot\frac{v_1^2v_2^2}{k^2}\cdot\frac{1}{\sqrt{\|\va_3\|^2\|\va_3'\|^2-(\va_3\cdot\va_3')^2}}
\end{align*}
Finally, integrating over the configuration space and using our parameters $v_1$, $v_2$, $k_1$, $k_2$ and $C$, we obtain the stated upper bound on $\langle\mathrm{Lk}^2\rangle$.
\end{Proof}

The integrand in the definition of $C$ has a singularity whenever $\va_3(s,t)$ and $\va_3'(s',t')$ are parallel or anti-parallel, and we will now demonstrate how to stipulate further conditions on the input data so that $C$ is finite. To this end, split the integral in to two parts, one close to the diagonal where $s=s'$ and $t=t'$ (which we will denote $\mu_1(\epsilon)$) and another away from the diagonal (denoted $\mu_1(\epsilon)^c$), by defining 
\begin{equation}\nonumber
\mu_1(\epsilon)=\{(s',t',s,t)\in\mathbb{T}^4\textrm{ where }|s-s'|<\epsilon\textrm{ and }|t-t'|<\epsilon\}
\end{equation}
 and writing $C=C_1+C_2$ where:
\begin{align*}\nonumber
C_1&=\int\displaylimits_{s=0}^{2\pi}\int\displaylimits_{t=0}^{2\pi}\int\displaylimits_{|s-s'|<\epsilon}\int\displaylimits_{|t-t'|<\epsilon}\frac{1}{\sqrt{\|\va_3(s,t)\|^2\|\va_3'(s',t')\|^2-(\va_3(s,t)\cdot\va_3'(s',t'))^2}}ds'dt'dsdt\\
C_2&=\int\displaylimits_{(s,t,s',t')\in\mu_1^c(\epsilon)}\frac{1}{\sqrt{\|\va_3(s,t)\|^2\|\va_3'(s',t')\|^2-(\va_3(s,t)\cdot\va_3'(s',t'))^2}}ds'dt'dsdt>0
\end{align*}
Next define that 
\begin{align*}\nonumber
\eta_1(\epsilon)&=\mathrm{min}_{(s,t,s',t')\in\mu_1}\sqrt{\frac{\|\va_3(s,t)\|^2\|\va_3'(s',t')\|^2-(\va_3(s,t)\cdot\va_3'(s',t'))^2}{(s-s')^2+(t-t')^2}}\\
\eta_2(\epsilon)&=\mathrm{min}_{(s,t,s',t')\in\mu_1^C}\sqrt{\|\va_3(s,t)\|^2\|\va_3'(s',t')\|^2-(\va_3(s,t)\cdot\va_3'(s',t'))^2}
\end{align*}
When $\eta_1(\epsilon)>0$, then we may bound $C_1$ as:
\begin{align*}\nonumber
C_1\leq\int\displaylimits_{s=0}^{2\pi}\int\displaylimits_{t=0}^{2\pi}\int\displaylimits_{|s-s'|<\epsilon}\int\displaylimits_{|t-t'|<\epsilon}\frac{1}{\eta_1(\epsilon)}\frac{1}{\sqrt{(s-s')^2+(t-t')^2}}ds'dt'dsdt
\end{align*}
Finally, when we stipulate that $v_1$ and $v_2$ are finite, then we have a more precise bound on $\langle\mathrm{Lk}^2\rangle$:
\begin{equation}\label{lk2bound}
\langle\mathrm{Lk}^2\rangle\leq\textrm{min}_{\epsilon}\frac{1}{(4\pi)^2}\frac{4v_1^2v_2^2}{\pi k^2}\big(C_1+\frac{\mathrm{vol}(\mu_1^c(\epsilon))}{\eta_2(\epsilon)}\big)
\end{equation}
\begin{Remark}
If we further restrict the data to be as in \eqref{orthogonalData}, then we can get an even more explicit bound on $\langle\mathrm{Lk}^2\rangle$. That is, starting from \eqref{starthere}, we have:
\begin{equation}\nonumber
\langle\mathrm{Lk}^2\rangle\leq\frac{4}{\pi}\frac{(\sum\displaylimits_{k=1}^N k^2c_k^2)(\sum\displaylimits_{k=1}^N k^2d_k^2)}{\sum\displaylimits_{k=0}^N c_k^2+\sum\displaylimits_{k=0}^N d_k^2}\frac{1}{(4\pi)^2}\int\displaylimits_{0}^{2\pi}\int\displaylimits_{0}^{2\pi}\int\displaylimits_{0}^{2\pi}\int\displaylimits_{0}^{2\pi}\frac{dsdtds'dt'}{\sqrt{(\sum\displaylimits_{k=0}^N c_k^2+\sum\displaylimits_{k=0}^N d_k^2)^2-(F(s,t,s',t'))^2}},
\end{equation}
where $F(s,t,s',t')=\vr_2(s)\cdot\vr_2(s')+\vr_1(t)\cdot\vr_1(t')$. If we consider again the limit as $n\rightarrow\infty$, then we will have:
\begin{equation}\nonumber
\langle\mathrm{Lk}^2\rangle\leq\frac{4}{\pi}\frac{\|\textbf{c}'\|_{l^2}^2\|\textbf{d}'\|_{l^2}^2}{\|\textbf{c}\|_{l^2}^2+\|\textbf{d}\|_{l^2}^2}\frac{1}{(4\pi)^2}\int\displaylimits_{0}^{2\pi}\int\displaylimits_{0}^{2\pi}\int\displaylimits_{0}^{2\pi}\int\displaylimits_{0}^{2\pi}\frac{dsdtds'dt'}{\sqrt{(\|\textbf{c}\|_{l^2}^2+\|\textbf{d}\|_{l^2}^2)^2-(F(s,t,s',t'))^2}},
\end{equation}
Further simplifying, if we take $c_k=d_k=\frac{1}{k^{\alpha}}$, then $F(s,t,s',t')$ can be written as a sum of polylogarithms. Moreover, it's clear that $\langle\mathrm{Lk}^2\rangle$ will therefore be bounded when $C$ is and when the sequence $\textbf{c}$ decays faster than $n^{(-3-\epsilon)/2}$. In the special case of $\alpha=1$ 
\end{Remark}
\begin{Remark}
A final remark to be made here, is that our result \eqref{mainResultTheorem} can be used to find the second moment of the degree of $2-$sphere self maps, however in these cases the configuration space allows for stronger singularities, and a more detailed analysis of convergence must be undertaken.
\end{Remark}

    \newchapter{Second Moments for Higher Dimensional Linking Integrals}{Second Moments for Higher Dimensional Linking Integrals}{Second Moments for Higher Dimensional Linking Integrals}
    \label{sec:LabelForChapter4}

        \section[Integration Over Codimension-$1$ Subspaces]{Integration Over Codimension-$1$ Subspaces}
        \label{sec:LabelForChapter4:Section1}

As in \cite{higherlink}, given two closed, disjoint, oriented manifolds $M_1$ and $M_2$ of respective dimensions $m$ and $n$ which are submanifolds of $\mathbb{R}^{N=m+n+1}$, then one may generalize the classic Gauss linking integral as:
\begin{equation}\
\mathrm{Lk}(M_1,M_2)=\frac{(-1)^{m+1}}{\mathrm{vol}(\mathbb{S}^{N-1})}\int_{M_1\times M_2}\frac{\mathrm{Det}(\textbf{x}-\textbf{y},\frac{d\textbf{x}}{ds_1},...,\frac{d\textbf{x}}{ds_m}\frac{d\textbf{y}}{dt_1},...,\frac{d\textbf{y}}{dt_n})}{\|\textbf{x}-\textbf{y}\|^N}ds_1...ds_m dt_1...dt_n,
\end{equation}
where $\textbf{x}(s_1,..,s_m)$ and $\textbf{y}(t_1,..,t_n)$ are local coordinates on the manifolds $M_1$ and $M_2$. This way, if we take as starting data two manifolds in $\mathbb{R}^{N'}$, where $N'=m+n+2$, then we may generate a random link of manifolds in $\mathbb{R}^{N}$ by picking an $N=m+n+1$ dimensional subspace of $\mathbb{R}^{N'}$ at random, and then orthogonally projecting the manifolds to the subspace. To compute the second moment of the linking number, we may mimic the method in the previous section almost line for line, except for the computation of the $I_{ii}$ integrals, and for these we will need to make an assumption, namely that $m+n+2$ is even. 
Now let $M_1$ and $M_2$ be manifolds in $\mathbb{R}^{N'}$ with local coordinates $\textbf{x}(s_1,s_2,...,s_m)$ and $\textbf{y}(t_1,t_2,...,t_n)$. As above, to make the notation more compact, we will make the following conventions:
\begin{align*}\nonumber
\va_i(s_1,s_2,...,s_m)=\frac{\partial\textbf{x}(s_1,s_2,...,s_m)}{\partial s_i}\textrm{ for $1\leq i\leq m$}\\
\va_j(t_1,t_2,...,t_n)=\frac{\partial\textbf{y}(t_1,t_2,...,t_n)}{\partial t_{j-m}}\textrm{ for $m+1\leq j\leq m+n$}\\
\va_{m+n+1}(s_1,s_2,...,s_m,t_1,t_2,...,t_n)=\textbf{y}(t_1,..,t_n)-\textbf{x}(s_1,..,s_m)
\end{align*}
with similar identifications for the vectors $\va_i'$. Also define that:
\begin{equation}\nonumber
\big[\va_1\colon\va_2\colon...\colon\va_{N}\big]_i=\sum_{j_1,j_2,..,j_N=1}^{N+1}\epsilon_{ij_1j_2...j_N}a_{1j_1}
a_{2j_2}...a_{Nj_n}
\end{equation}

 We will now prove the following:
\begin{Theorem}\nonumber
Given two closed, disjoint, oriented manifolds $M_1$ and $M_2$ of respective dimensions $m$ and $n$ in $\mathbb{R}^{m+n+2}$, where $m+n$ is even, then the value of $\mathrm{Lk}^2$ averaged over all orthogonal projections to $m+n+1$ dimensional subspaces, is given by the following integral:
\begin{equation}\nonumber
\frac{1}{\mathrm{vol}(\mathbb{S}^{N-1})^2}\int_{M_1\times M_2}I_{\langle\mathrm{Lk}^2\rangle}(s_1,...,s_m,t_1,..,t_n)ds_1...ds_m dt_1...dt_n
\end{equation}
where
\begin{align*}\label{eq:result4THM}
&I_{\langle\mathrm{Lk}^2\rangle}(A,A')\|\va_N\|^N\|\va'_N\|^N=\\
&I_{1,1}(a,b)(\big[\vb_2\colon\vb_3\colon...\colon\vb_{m+n+2}\big]\cdot\big[\va_1\colon\va_2\colon...\colon\va_{m+n+1}\big])(\big[\vb_2\colon\vb_3\colon...\colon\vb_{m+n+2}\big]\cdot\big[\va'_1\colon\va'_2\colon...\colon\va'_{m+n+1}\big])+\\
&I_{2,2}(a,b)(\big[\vb_1\colon\vb_3\colon...\colon\vb_{m+n+2}\big]\cdot\big[\va_1\colon\va_2\colon...\colon\va_{m+n+1}\big])(\big[\vb_1\colon\vb_3\colon...\colon\vb_{m+n+2}\big]\cdot\big[\va'_1\colon\va'_2\colon...\colon\va'_{m+n+1}\big])+\\
&I_{3,3}(a,b)\sum_{i=3}^{m+n+2}\vb_i\cdot\big[\va_1\colon\va_2\colon...\colon\va_{m+n+1}\big]\vb_i\cdot\big[\va'_1\colon\va'_2\colon...\colon\va'_{m+n+1}\big]
\end{align*}
and $\{\vb_i\}_{i=1,...,m+n+2}$ is a basis for $\mathbb{R}^{m+n+2}$ such that $\vb_1,\vb_2\in\mathrm{span}\{\va_{m+n+1},\va'_{m+n+1}\}$ and $I_{3,3}(m+n)$ is a function of $a=\vb_1\cdot\frac{\va_{m+n+1}}{\|\va_{m+n+1}\|}$ and $b=\vb_2\cdot\frac{\va_{m+n+1}}{\|\va_{m+n+1}\|}$.
\end{Theorem}
\begin{Proof}
Take $A=(\va_1,\va_2,...,\va_{N})$, $A'=(\va_1',\va_2',...,\va_{N}')$, $dV=\mathrm{exp}(-\|\vV\|^2/2)d\vV$, and define:
\begin{equation}\nonumber
I_{\langle\mathrm{Lk}^2\rangle}(A,A')=\frac{1}{(2\pi)^{N'/2}}\int_{\mathbb{R}^{N'}}\frac{\mathrm{Det}((\va_1,\va_2,...,\va_{N},v/\|v\|)}{\|\mathrm{proj}_{v^{\perp}}(\va_{N})\|^{N}}\frac{\mathrm{Det}((\va'_1,\va'_2,...,\va'_{N},v/\|v\|)}{\|\mathrm{proj}_{v^{\perp}}(\va_{N}')\|^{N}}dV.
\end{equation}
As in the previous section, define that:
\begin{equation}\nonumber
f(A,A')=\frac{\mathrm{det}(A^TA)^{1/2}\mathrm{det}(A'^TA')^{1/2}}{\|\va_{N}\|^{N}\|\va_{N}'\|^{N}}
\end{equation}
so that by applying Gram-Schmidt we obtain:
\begin{equation}\nonumber
I_{\langle\mathrm{Lk}^2\rangle}(A,A')=\frac{f(A,A')}{(2\pi)^{N'/2}}\int_{\mathbb{R}^{N'}}\frac{\mathrm{Det}(\vq_1,\vq_2,...,\vq_{N},v/\|v\|)}{\|\mathrm{proj}_{v^{\perp}}(\vq_{N})\|^{N}}\frac{\mathrm{Det}((\vq_1',\vq_2',...,\vq_{N}',v/\|v\|)}{\|\mathrm{proj}_{v^{\perp}}(\vq_{N}')\|^{N}}dV
\end{equation}
Similar to the argument in the second section, choose a basis such that:
\begin{align*}\nonumber
\vb_1=\frac{\vq_{m+n+1}+\vq_{m+n+1}'}{\|\vq_{m+n+1}+\vq_{m+n+1}'\|}\\
\vb_2=\frac{\vq_{m+n+1}-\vq_{m+n+1}'}{\|\vq_{m+n+1}-\vq_{m+n+1}'\|}\\
\textrm{ and }\vb_3,...,\vb_{m+n+2}\in\mathrm{span}^\perp\{\vb_1,\vb_2\}
\end{align*}
This choice is made to simplify the denominators in $I(A,A')$. The relevant Gaussian integrals that appear in computing $I(A,A')$ upon expanding the product of the determinants and integrating are of the form:
\begin{equation}\nonumber
I_{i,j}=\frac{1}{(2\pi)^{N'/2}}\int_{\mathbb{R}^{N'}}\frac{\|\textbf{v}\|^{2(m+n)}v_iv_jdV}{((bv_1-av_2)^2+v_3^2+...+v_{N'}^2)^{N/2}((bv_1+av_2)^2+v_3^2+...+v_{N'}^2)^{N/2}}
\end{equation}
Expanding the determinants we have:
\begin{align*}\nonumber
&\frac{I_{\langle\mathrm{Lk}^2\rangle}(A,A')}{f(A,A')}=\\
&I_{1,1}(a,b)(\big[\vb_2\colon\vb_3\colon...\colon\vb_{m+n+2}\big]\cdot\big[\vq_1\colon\vq_2\colon...\colon\vq_{m+n+1}\big])(\big[\vb_2\colon\vb_3\colon...\colon\vb_{m+n+2}\big]\cdot\big[\vq'_1\colon\vq'_2\colon...\colon\vq'_{m+n+1}\big])+\\
&I_{2,2}(a,b)(\big[\vb_1\colon\vb_3\colon...\colon\vb_{m+n+2}\big]\cdot\big[\vq_1\colon\vq_2\colon...\colon\vq_{m+n+1}\big])(\big[\vb_1\colon\vb_3\colon...\colon\vb_{m+n+2}\big]\cdot\big[\vq'_1\colon\vq'_2\colon...\colon\vq'_{m+n+1}\big])+\\
&\sum_{i=3}^{m+n+2}I_{i,i}(a,b)\vb_i\cdot\big[\vq_1\colon\vq_2\colon...\colon\vq_{m+n+1}\big]\vb_i\cdot\big[\vq'_1\colon\vq'_2\colon...\colon\vq'_{m+n+1}\big]
\end{align*}
In a similar fashion to the case considered in the previous chapter, the terms including $I_{1,1}$ and $I_{2,2}$ cancel out precisely when $m+n=2$. When $m+n>2$ no clear cancellation occurs. Moreover, it is clear that $I_{i,j}=0$ for $i\neq j$ and $I_{3,3}=I_{4,4}=...=I_{m+n+2,m+n+2}$. With these facts we may simplify the above equation to write:
\begin{align*}\nonumber
&I_{\langle\mathrm{Lk}^2\rangle}(A,A')\|\va_N\|^N\|\va'_N\|^N=\\
&I_{1,1}(a,b)(\big[\vb_2\colon\vb_3\colon...\colon\vb_{m+n+2}\big]\cdot\big[\va_1\colon\va_2\colon...\colon\va_{m+n+1}\big])(\big[\vb_2\colon\vb_3\colon...\colon\vb_{m+n+2}\big]\cdot\big[\va'_1\colon\va'_2\colon...\colon\va'_{m+n+1}\big])+\\
&I_{2,2}(a,b)(\big[\vb_1\colon\vb_3\colon...\colon\vb_{m+n+2}\big]\cdot\big[\va_1\colon\va_2\colon...\colon\va_{m+n+1}\big])(\big[\vb_1\colon\vb_3\colon...\colon\vb_{m+n+2}\big]\cdot\big[\va'_1\colon\va'_2\colon...\colon\va'_{m+n+1}\big])+\\
&I_{3,3}(a,b)\sum_{i=3}^{m+n+2}\vb_i\cdot\big[\va_1\colon\va_2\colon...\colon\va_{m+n+1}\big]\vb_i\cdot\big[\va'_1\colon\va'_2\colon...\colon\va'_{m+n+1}\big]
\end{align*}
\end{Proof}
 Now we will focus on the functional form of the function $I_{3,3}(a,b)$, and to do so, we will reduce it to an integral over a $3$-sphere as in the previous section. We have that:
\begin{equation}\nonumber
I_{3,3}=\frac{1}{(2\pi)^{N'/2}}\int_{\mathbb{R}^{N'}}\frac{\|\textbf{v}\|^{2(m+n)}v_3^2dV}{((bv_1-av_2)^2+v_3^2+...+v_{N'}^2)^{N/2}((bv_1+av_2)^2+v_3^2+...+v_{N'}^2)^{N/2}}
\end{equation}
Now change the coordinates $v_4,v_5,...,v_{m+n+2}$ to $(m+n-1)$-dimensional spherical coordinates so that $v_4^2+v_5^2+...+v_{m+n+2}^2=r^2$ and $dv_4dv_5...dv_{m+n+2}=r^{m+n-2}drdvol(\mathbb{S}^{m+n-2})$ so that the above integral becomes:
\begin{equation}\nonumber
I_{3,3}=A_1\int_{-\infty}^\infty\int_{-\infty}^\infty\int_{-\infty}^\infty\int_{0}^\infty\frac{(v_1^2+v_2^2+v_3^2+r^2)^{(m+n)}v_3^2dV(r)}{((bv_1-av_2)^2+v_3^2+r^2)^{N/2}((bv_1+av_2)^2+v_3^2+r^2)^{N/2}}
\end{equation}
where $dV(r)=\mathrm{exp}(-(v_1^2+v_2^2+v_3^2+r^2)/2)r^{m+n-2}drdv_1dv_2dv_3$ and $A_1=\frac{\mathrm{vol}(\mathbb{S}^{N'-4})}{(2\pi)^{N'/2}}$.
Now make the change of coordinates $r\rightarrow v_4$, and realize that since $m+n+2$ was chosen to be even, then the whole integrand is even and we obtain that:
\begin{equation}\nonumber
I_{3,3}=\frac{A_1}{2ab}\int_{\mathbb{R}^4}\frac{(v_1^2/b^2+v_2^2/a^2+v_3^2+v_4^2)^{(m+n)}v_3^2v_4^{m+n-2}}{((v_1-v_2)^2+v_3^2+v_4^2)^{N/2}((v_1+v_2)^2+v_3^2+v_4^2)^{N/2}}dV'
\end{equation}
where again we have that $a=\vq_N\cdot\vb_1$ and $b=\vq_N\cdot\vb_2$ and used a further change of coordinates: $v_1\rightarrow v_1/b$, and $v_2\rightarrow v_2/a$ and have written that $dV'=\mathrm{exp}(-(v_1^2/a^2+v_2^2/b^2+v_3^2+v_4^2)/2)$. Finally, change to toroidal coordinates to obtain:
\begin{equation}\nonumber
I_{3,3}=\frac{A_1}{2ab}\int_0^{\pi/2}\int_0^{2\pi}\int_0^{2\pi}\int_{0}^\infty\frac{k_1^{m+n}\sin^{N}(\sigma)\cos(\sigma)\sin^{N'-4}(\phi)\cos^2(\phi)dV'}{(1-\cos^4(\sigma)\sin^2(2\theta))^{N/2}}
\end{equation}
where again $k_1(\theta,\sigma)=\sin^2\big(\sigma)+\cos^2(\sigma)((\cos(\theta)/b)^2+(\sin(\theta)/a)^2\big)$ and moreover $dV'=\mathrm{exp}(-k_1r^2/2)r^{m+n+1}drd\theta d\phi d\sigma$. To integrate out the radial dependence, first change $r\rightarrow \frac{r}{\sqrt{k_1}}$:
\begin{eqnarray*}\nonumber
I_{3,3}=\\
\frac{A_1}{2ab}\int_0^{\pi/2}\int_0^{2\pi}\int_0^{2\pi}\int_{0}^\infty\frac{k_1^{(N'-4)/2}\sin^{N}(\sigma)\cos(\sigma)\sin^{N'-4}(\phi)\cos^2(\phi)\mathrm{exp}(-r^2/2)dV'}{(1-\cos^4(\sigma)\sin^2(2\theta))^{N/2}}
\end{eqnarray*}
where now $dV'=r^{N}drd\theta d\phi d\sigma$. Notice that for the case of links when $m=n=1$ that this integral is especially easy. Now integrate over the radial coordinate:
\begin{equation}\nonumber
I_{3,3}=A_2\int_0^{\pi/2}\int_0^{2\pi}\int_0^{2\pi}\frac{k_1^{(N'-4)/2}\sin^{N}(\sigma)\cos(\sigma)\sin^{N'-4}(\phi)\cos^2(\phi)}{(1-\cos^4(\sigma)\sin^2(2\theta))^{N/2}}d\theta d\phi d\sigma
\end{equation}
where the new constant $A_2$ is given by:
\begin{equation}\nonumber
 A_2=\frac{\mathrm{vol}(\mathbb{S}^{m+n-2})\big(\frac{m+n}{2}\big)!2^{(m+n)/2}}{2ab(2\pi)^{(m+n+2)/2}}
\end{equation}
Lastly, define:
\begin{equation}\nonumber
I(\theta,\phi,\sigma,m,n)=A_2\frac{k_1^{(N'-4)/2}\sin^{N}(\sigma)\cos(\sigma)\sin^{N'-4}(\phi)\cos^2(\phi)}{(1-\cos^4(\sigma)\sin^2(2\theta))^{N/2}}.
\end{equation}
Given the sufficiently simple form above, we may then compute some values of $I_{3,3}(m+n)$ for different values of $m$ and $n$ such that $m+n$ is even. Here is a list of a few of them, calculated using a symbolic integration in Wolfram Mathematica:
\begin{equation}\nonumber
I_{3,3}(m+n=2)=\frac{1}{\pi ab}\textrm{ (see section 2) }
\end{equation}
\begin{equation}\nonumber
I_{3,3}(m+n=4)=\frac{1+4a^2b^2}{9\pi a^3b^3}
\end{equation}
\begin{equation}\nonumber
I_{3,3}(m+n=6)=\frac{9b^4+2a^2b^2(5+16b^2)+a^4(9+32b^2+128b^4)}{450\pi a^5b^5}
\end{equation}
\begin{align*}\nonumber
I_{3,3}(m+n&=8)\\
&=\frac{15b^6+3a^2b^4(7+20b^2)+3a^6(1+4b^2)(5+64b^4)+a^4b^2(21+56b^2+192b^4)}{3675\pi a^7b^7}
\end{align*}
Notice that this computation is very similar to computation in \ref{codim1lemma}, and uses the fact that the subspaces being projected on to are codimension $1$ so that the Grassmannian is identified with a sphere. With this identification we were able to compute $\langle\mathrm{Lk}^2(M,N)\rangle$ using the Lebesgue measure on the sphere, however, this assumption on the codimension can be removed with an argument very much similar to that in \ref{theRESULT}. Given this computation, it is then feasible to bound $\langle\mathrm{Lk}(M,N)^2\rangle$ by bounding the configuration space integrals in a way analogous to the method at the end of section 3.3. 
\begin{Remark}
As this calculation was only done for subspaces of codimension $1$, we would need to extend to higher codimension in order to get results analogous to our main theorem \ref{mainResultTheorem} obtained for the linking number of dimension $1$ curves in $\mathbb{R}^3$. We will undertake this computation in a future work. 
\end{Remark}


    \newchapter{Numerical Study: Petal Diagrams}{Numerical Study: Petal Diagrams}{Numerical Study: Petal Diagrams}
    \label{sec:LabelForChapter5}
    
        \section[Model Details]{Model Details}
        \label{sec:LabelForChapter5:Section1} 
In this chapter, we will explore a very special case of the input curve $\vr(t)$, inspired by \cite{Hass} in order to show how the results in the previous chapters may be used. The model considered in their paper is called the Petaluma model and is motivated by the observation in \cite{Adams} that an embedding of a knot may be arranged so that there is a projection to a plane $P=\vV^\perp$ where the knot diagram obtained is a rose with $n$ petals. One may enhance this diagram to account for the crossing data by labeling the strands with the heights through which they pass through the axis determined by $\vV$. With this observation in hand, the Petaluma model is defined by fixing a petal diagram and then choosing permutations of the heights at random.\\
 We will now focus on finding a space curve  $\vr(t)$ in some $\mathbb{R}^N$ such that our random projection model can approximate the Petaluma model. To do so, assume $k$ is odd (the even case follows similarly) and subdivide the interval $\big[0,\pi\big]$ in to $2k$ equal length subintervals, $\{U_i=\big[t_i,t_{i+1}\big]\}$, and define a space curve $\vr(t):\mathbb{S}^1\rightarrow\mathbb{R}^{2+k}$ such that: 
\begin{equation}\label{dataintegral}
\vr(t,\epsilon,k)=\cos(kt)\cos(t)\ve_1+\cos(kt)\sin(t)\ve_2+\textbf{R}(t,\epsilon,k)
\end{equation}
where:
\begin{equation}\label{taper}
\textbf{R}(t,\epsilon,k)=\sum_{i=1}^k\big(\mathds{1}_{U_{2i-1}}(t)\frac{2k\epsilon(t-t_i)}{\pi}+\mathds{1}_{U_{2i}}(t)\frac{2k\epsilon(t_{i+1}-t)}{\pi}\big)\ve_{i+2}
\end{equation}
and $\mathds{1}_{U_i}$ is the indicator function on the interval $U_i$. Intuitively, in the first two coordinates we have the parametrization for a rose with $k$ petals, and in the remaining $k$ coordinates we have a linear function that for the $i^{th}$ strand is supported in the $(2+i)th$ coordinate and runs from $0$ at the outermost part of the  strand to $\epsilon$ at the center and then back to $0$. For the case of polygonal knots we may take a piecewise linear approximation of the rose diagram in the first two coordinates. The piecewise linear petal diagram is related to the grid model also considered in \cite{Hass}. 

\begin{figure}\label{petals}
\centering
\includegraphics[scale=0.5]{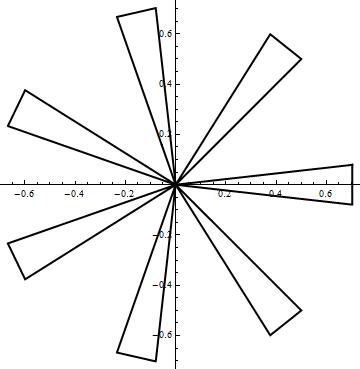}
\includegraphics[scale=0.5]{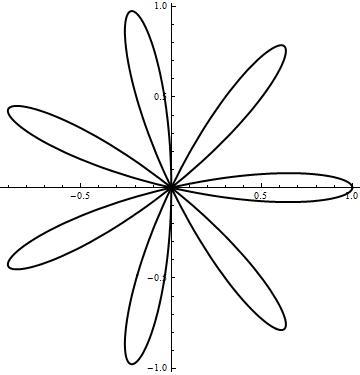}
\caption{Petal diagrams in the first $2$ coordinates. The left diagram is used in the PL case, while the right is used in our approximation of the Petaluma model.}
\end{figure}     
           
        \section[Expected Total Curvature]{Expected Total Curvature}
        We will be interested in computing the expected curvature in both cases, the latter of which is calculated as the sum of the turning angles. 
Given this particularly simple form for $\vr(t)$, we then have that:
\begin{equation}\nonumber
\|\vr'(t,\epsilon,k)\|^2=\frac{4\epsilon^2k^2}{\pi^2}+\frac{1}{2}((1+k^2)-(k^2-1)\cos(2kt))
\end{equation}
\begin{equation}\nonumber
\|\vr''(t,\epsilon,k)\|^2=\frac{1}{2}(1+6k^2+k^4+(k^2-1)^2\cos(2kt))
\end{equation}
\begin{equation}\nonumber
(\vr'(t,\epsilon,k)\cdot\vr''(t,\epsilon,k))^2=\frac{1}{4}k^2(k^2-1)^2\sin(2kt)^2
\end{equation}
\begin{equation}\nonumber
I_{\langle\kappa(C)\rangle}(\vr'(t,\epsilon,k),\vr''(t,\epsilon,k))=\frac{\sqrt{\|\vr'(t,\epsilon,k)\|^2\|\vr''(t,\epsilon,k)\|^2-(\vr'(t,\epsilon,k)\cdot\vr''(t,\epsilon,k))^2}}{\|\vr'(t,\epsilon,k)\|^2}
\end{equation}
For small $\epsilon$, this integral will get closer and closer to simply computing the curvature of the petal diagram. Numerical calculations show that $\int_0^\pi I_{\langle\kappa(C)\rangle}(\vr'(t,0,k),\vr''(t,0,k))=\langle\kappa(\mathrm{C})\rangle\approx\pi(k+1)dt$, and for small enough values of epsilon we have:
\begin{equation}
\int_0^\pi I_{\langle\kappa(C)\rangle}(\vr'(t,\epsilon,k),\vr''(t,\epsilon,k))\leq\pi(k+1)
\end{equation}
Interestingly, in the case when $k=3$ and $\epsilon=.5$ (though this can be further tuned), then we may get an idea about the approximate density of the unknot. That is, exactly like in the proof of corollary 25 in \cite{cantarella} (also see \cite{Shonkwiler1}) which uses the Fary-Milnor theorem (\cite{Milnor}), if we let $x$ denote the fraction of knots with curvature greater than $4\pi$, then
\begin{equation}\nonumber
\langle\kappa(C)\rangle > 4\pi x +2\pi(1-x)
\end{equation}
and solving for $x$ we see that:
\begin{equation}\nonumber
x<\frac{\langle\kappa(C)\rangle}{2\pi}-1
\end{equation}
When $\epsilon=.5$ a numerical integration gives that $\langle\kappa(C)\rangle\approx 9.72$ so that $x<.54$, that is, at least approximately $48$ percent are unknotted.
\Remark{The same results hold when considering a piecewise linear approximation to the petal diagram as discussed in the above remark. It would be interesting to compute the total torsion as well, especially given its relation to the self-linking number.} 
        \label{sec:LabelForChapter5:Section2}
    
        \section[$\langle\mathrm{Lk}^2\rangle$]{$\langle\mathrm{Lk}^2\rangle$}

We may modify the curve $\vr(t)$ defined in the first section to define a model for random links. To do so, we define the following two space curves:
\begin{eqnarray}\label{linkdata}
\vr_1(t,\epsilon,k)=\cos(kt)\cos(t)\ve_1+\cos(kt)\sin(t)\ve_2+\textbf{R}(t,\epsilon,k)\nonumber\\
\vr_2(t',\epsilon,k)=\left(\begin{array}{cc}
\cos(\frac{k-2}{k}\pi) & -\sin(\frac{k-2}{k}\pi)\\
\sin(\frac{k-2}{k}\pi) & \cos(\frac{k-2}{k}\pi)\end{array}\right)\left(\begin{array}{c}
\cos(kt')\cos(t')\\
\cos(kt')\sin(t')\end{array}\right)+\textbf{R}(t',\epsilon,k)\nonumber
\end{eqnarray}
in $\mathbb{R}^{2+2k}$, again where $\textbf{R}(t',\epsilon,k)$ is as in \eqref{taper}, and $k$ zeroes are appended to the end of the vector $\vr_1(t)$, and where another $k$ zeroes are inserted after the second position in $\vr_2(t)$. That is, we take each component to have an equal number of petals and the second component is a rotation of the first in the first two coordinates, as is shown in figure 5.3.1.
\begin{figure}[H]\label{petallink}
\centering
\includegraphics[scale=0.5]{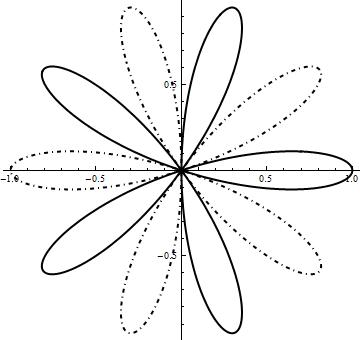}
\caption{Petal diagram with $k=5$ petals per component.}
\end{figure} 

 In tables 5.1 and 5.2 we illustrated two ways of computing $\langle\mathrm{Lk}^2\rangle$, the first of which was found by sampling $10^6$ links using our model (subspaces were chosen by sampling Gaussian random matrices $U\in\mathrm{Mat}_{n,3}$), computing the linking number for each sample by taking a piecewise linear approximation of the link, and then using the algorithm in \cite{Arai} to compute the linking number, and then finding the sum of the squared standard deviation and the squared mean of our samples. The second was found by computing the integral in Theorem \ref{mainResultTheorem} using a Monte Carlo integration in Wolfram Mathematica (also with $10^6$ integrand evaluations).
\begin{center}
\begin{table}[H]\begin{tabular}{ |l|l|l| } 
\hline \multicolumn{3}{ |c| }{$\langle\mathrm{Lk}^2\rangle$ Comparison with $\epsilon=1$} \\ 
\hline 
k (number of petals/ component) & Sampled $\langle\mathrm{Lk}^2\rangle$  & Monte Carlo $\langle\mathrm{Lk}^2\rangle$  \\ 
\hline 3 & 0.789 & $0.834\pm 0.0255$ \\ 
\hline 5 & 2.211 & $2.3946\pm 0.126$ \\ 
\hline 7 & 4.299 & $4.83\pm 0.346$ \\ 
\hline 9 & 7.042 & $8.139\pm 0.803$ \\ 
\hline \end{tabular}\label{tab:mysecondtable}\caption{$\langle\mathrm{Lk}^2\rangle$ Comparison with $\epsilon=1$}\end{table}

\begin{table}[H]\begin{tabular}{ |l|l|l| } 
\hline 
\multicolumn{3}{ |c| }{$\langle\mathrm{Lk}^2\rangle$ Comparison with $\epsilon=.1$} \\ 
\hline k (number of petals/ component) & Sampled $\langle\mathrm{Lk}^2\rangle$  & Monte Carlo $\langle\mathrm{Lk}^2\rangle$  \\ 
\hline 3 & 0.543 & $0.56\pm 0.148$ \\ 
\hline 5 & 1.744 & $1.97\pm 0.984$ \\ 
\hline 7 & 3.611 & $3.675\pm 1.825$ \\ 
\hline 9 & 6.170 & $8.39\pm 2.838$ \\ \hline \end{tabular} \label{tab:myfirsttable}
\caption{$\langle\mathrm{Lk}^2\rangle$ Comparison with $\epsilon=.1$}\end{table}
\end{center} 
\begin{Remark}
Notice that in the second table, that the values of $\langle\mathrm{Lk}^2\rangle$ have larger errors compared to the first table. We expect this is due to the numerical instability that arises when $\epsilon$ is small since the strands become closer and force the denominator in the expression for $I_{\langle\mathrm{Lk}^2\rangle}$ to become close to zero, unlike in the results in the first table. It would be interesting to find a stable numerical method for computing these configuration space integrals. Alternate methods are considered in the appendix.
\end{Remark}
     
        \label{sec:LabelForChapter5:Sectionn}
        
        \newchapter{Future Directions}{Future Directions}{Future Directions}
    \label{sec:LabelForChapter6}
    
In this work we have found a way to compute the second moment of the linking number for a rather general model of random links. Of course, it would be interesting to find the moments of other knot and link invariants arising from configuration space integrals. One such example, the writhe of a knot, although not a knot invariant, provides some insight in to the steps that would be involved. To compute the writhe, one starts with a differentiable closed curve $\gamma(t)$ in $\mathbb{R}^3$, and then computes an analogue of the Gauss linking integral:
\begin{equation}\label{writhe}
\mathrm{Wr}(\gamma)=\frac{1}{4\pi}\int\displaylimits_0^{2\pi}\int\displaylimits_0^{2\pi}\frac{(\dot{\vr}_1(t)\times\dot{\vr}_2(t'))\cdot(\vr_2(t')-\vr_1(t))}{\|\vr_2(t')-\vr_1(t))\|^3}dt'dt
\end{equation}
The form of this integral is especially well suited to the analysis in Chapter $3$, and in fact an almost identical result would be obtained, however some care would need to be taken in proving that the order of integration may be reversed due to the singularity in the integrand in \eqref{writhe}.\\

\section[Finite Type Invariants]{Finite Type Invariants}

After understanding the singularities involved in computing $\langle\mathrm{Wr}^2\rangle$, then other knot invariants, such as finite type Vassiliev invariants \cite{Kontsevich}\cite{VassilievInbook}\cite{Vassiliev2}, which can be computed with configuration space integrals (see \cite{BARNATAN1995423}, \cite{BottTaubes}, and \cite{dthurston} for example) may be explored. For example, in \cite{LW} the authors discuss a knot invariant, $v_2=\rho_1(\gamma)+\rho_2(\gamma)$ for a smooth embedding $\gamma(t):\mathbb{S}^1\rightarrow\mathbb{R}^3$, by integrating over the configuration spaces: 
\begin{eqnarray*}
&\Delta_4=\{(t_1,t_2,t_3,t_4)|0<t_1<t_2<t_3<t_4<1\}\textrm{, and }\nonumber\\
&\Delta_3=\{(t_1,t_2,t_3,\textbf{z})|0<t_1<t_2<t_3<1,\textbf{z}\in\mathbb{R}^3-\{(\gamma(t_1),\gamma(t_2),\gamma(t_3))\}\},\nonumber
\end{eqnarray*}   
where:
\begin{eqnarray}\label{Yintegral}
&\rho_1(\gamma)=-\frac{1}{32\pi^3}\int\displaylimits_{\Delta_3(\gamma)}\mathrm{Det}\big[E(\textbf{z},
t_1),E(\textbf{z},t_2),E(\textbf{z},t_3)\big]d\textbf{z}dt_1dt_2dt_3\\
&\textrm{ with }E(\textbf{z},t)=\frac{(\textbf{z}-\gamma(t))\times\gamma'(t)}{\|\textbf{z}-\gamma(t)\|^3}\nonumber
\end{eqnarray}  
and
\begin{eqnarray*}\nonumber
&\rho_2(\gamma)=\frac{1}{8\pi^2}\int\displaylimits_{\Delta_4}\frac{\mathrm{Det}\big[\gamma(t_3)-\gamma(t_1)
,\gamma'(t_3),\gamma'(t_1)\big]}{\|\gamma(t_3)-\gamma(t_1)
\|^3}\frac{\mathrm{Det}\big[\gamma(t_4)-\gamma(t_2)
,\gamma'(t_4),\gamma'(t_2)\big]}{\|\gamma(t_4)-\gamma(t_2)\|^3}dt_1dt_2dt_3dt_4\nonumber
\end{eqnarray*}
It would be interesting to compute quantities like $\langle v_2\rangle$ and $\langle v_2^2\rangle$ for the model of random knots we have considered in this work. Unlike the linking number, $v_2(\gamma)=v_2(\gamma')$, so the first moment will be different from zero.  Interestingly, the integrand appearing in the definition of $\rho_2(\gamma)$ has similar tensorial properties as the integrand involved in computing $\mathrm{Lk}^2$ and so we may compute $\langle\rho_2(\gamma)\rangle$ by setting
\begin{align*} 
\va_1(t_1)=\gamma'(t_1)\textrm{ , }\va_2(t_3)=\gamma'(t_3)\textrm{ , } \va_3(t_3,t_1)=\gamma(t_3)-\gamma(t_1)\\
\va_1'(t_2)=\gamma'(t_2)\textrm{ , }\va_2'(t_4)=\gamma'(t_4)\textrm{ , } \va_3'(t_2,t_4)=\gamma(t_4)-\gamma(t_2)\\
A(t_1,t_3)=\big[\va_1(t_1),\va_2(t_3),\va_3(t_1,t_3)\big]\textrm{ and } A'(t_2,t_4)=\big[\va_1'(t_2),\va_2'(t_4),\va_3'(t_2,t_4)\big] 
\end{align*}
so that:
\begin{eqnarray}
\langle \rho_2(\gamma)\rangle =\int_{\Delta_4}I_{\langle\mathrm{Lk}^2\rangle}(A(t_1,t_3),A'(t_2,t_4))dt_1dt_2dt_3dt_4
\end{eqnarray}
where $I_{\langle\mathrm{Lk}^2\rangle}(A,A')$ is the same function as in the main result \eqref{GaussianResult}. It would also be interesting to find $\langle\rho_2(\gamma)\rangle$ in terms of an integral over a configuration space of a function of the initial data $A$ and $A'$, since doing so would allow us to find $\langle v_2\rangle =\langle\rho_1(\gamma)\rangle+\langle\rho_2(\gamma)\rangle$. We note here two important facts which might lead to simplifying computations of $\langle\rho_1(\gamma)\rangle$:

The integral over $\textbf{z}\in\mathbb{R}^3-\{(\gamma(t_1),\gamma(t_2),\gamma(t_3))\}$ appearing in the configuration space $\Delta_3$ can be computed in closed form (see \cite{GUADAGNINI1990581} and its applications in \cite{2014arXiv1401.1154F} and \cite{MCVassiliev}), and reduced to an integral over:
\begin{equation}
\Delta_3'=\{(t_1,t_2,t_3)|0<t_1<t_2<t_3<1\}.
\end{equation}
In particular, using the result and a modification of the notation used in \cite{MCVassiliev}, we can define the following functions:
\begin{flalign}
&C_1(\vx,\vy,\vz)=\frac{2\pi}{\|\vx\|\|\vy\|\|\vz\|}&\nonumber\\
&C_2(\vx,\vy,\vz)=\frac{1}{\|\vx\|\|\vy\|+\vx\cdot \vy}&\nonumber\\
&C_3(\vx,\vy,\vz)=\|\vx\|+\|\vy\|-\|\vz\|&\nonumber\\
&&\nonumber\\
\cdot &F(\vx,\vy,\vz,\tilde{\vx},\tilde{\vy},\tilde{\vz})=&\nonumber\\
& C_1C_2C_3\Biggl((\tilde{\vy}\cdot\tilde{\vz})(\tilde{\vx}\cdot \vz)+(\tilde{\vx}\cdot\tilde{\vz})(\tilde{\vy}\cdot \vy)-(\tilde{\vx}\cdot\tilde{\vy})(\tilde{\vz}\cdot \vx))\Biggr) &\label{check0}\\
-&C_1C_2^2C_3\Biggl(\mathrm{Det}\left[\tilde{\vy},\vx,\vy\right]\mathrm{Det}\left[\tilde{\vz},\tilde{\vx},\frac{\vx\|\vy\|+\vy\|\vx\|}{\|\vy\|}\right]&\nonumber\\
&+\mathrm{Det}\left[\tilde{\vz},\vx,\vy\right]\mathrm{Det}\left[\tilde{\vy},\tilde{\vx},\frac{\vx\|\vy\|+\vy\|\vx\|}{\|\vx\|}\right]\Biggr)&\label{check1}\\
&+C_1C_2\Biggl( \mathrm{Det}\left[\tilde{\vy},\vx,\vy\right]\cdot\mathrm{Det}\left[\tilde{\vz},\tilde{\vx},\vy\frac{\|\vz\|-\|\vx\|}{\|\vy\|^2}+\vz\frac{\|\vx\|+\|\vy\|}{\|\vz\|^2}\right]&\nonumber\\
&+\mathrm{Det}\left[\tilde{\vz},\vx,\vy\right]\mathrm{Det}\left[\tilde{\vy},\tilde{\vx},\vx\frac{\|\vz\|-\|\vy\|}{\|\vx\|^2}-\vz\frac{\|\vx\|+\|\vy\|}{\|\vz\|^2}\right]\Biggr)\label{check2},
\end{flalign}
so that the integral \eqref{Yintegral} can be rewritten as:
\begin{align*}
\rho_1(\gamma)&=-\frac{1}{32\pi^3}\int\displaylimits_{\Delta_3(\gamma)}\mathrm{Det}\big[E(\textbf{z},
t_1),E(\textbf{z},t_2),E(\textbf{z},t_3)\big]d\textbf{z}dt_1dt_2dt_3\\
&=-\frac{1}{32\pi^3}\int\displaylimits_{\Delta_3'}F\Biggl(\gamma(t_2)-\gamma(t_3),\gamma(t_1)-\gamma(t_3),\gamma(t_2)-\gamma(t_1),\dot{\gamma(t_1)},\dot{\gamma(t_2)},\dot{\gamma(t_3)}\Biggr)dt_1t_2dt_3
\end{align*}
We note that despite the complicated form of the function\footnote{For convenience, we have included an implementation of $F(\vx,\vy,\vz,\tilde{\vx},\tilde{\vy},\tilde{\vz})$ as a Wolfram Mathematica \cite{Mathematica} function, as it may be useful for experimental applications.} $F(\vx,\vy,\vz,\tilde{\vx},\tilde{\vy},\tilde{\vz})$, it is still a function of the pairwise inner products and norms of the vectors $\vx$, $\vy$, $\vz$, $\tilde{\vx}$, $\tilde{\vy}$, and $\tilde{\vz}$. Unfortunately however, the function $F$ does not have the same tensorial properties as the integrands appearing in the computation of $\langle\mathrm{Lk}^2\rangle$ and $\langle\rho_2\rangle$. 
With this setup, we may then define the function:
\begin{equation}\label{YGaussian}
I_{\langle\rho_1\rangle}(\va_1,\va_2,\va_3,\tilde{\va_1},\tilde{\va_2},\tilde{\va_3})=\int\displaylimits_{U\in\mathrm{Mat}_{n,3}}\frac{F(U^T\va_1,U^T\va_2,U^T\va_3,U^T\tilde{\va_1},U^T\tilde{\va_2},U^T\tilde{\va_3})}{(2\pi)^{3\cdot n/2}}e^{\frac{-Tr(U^TU)}{2}}dU.
\end{equation}
so that for a space curve $\gamma$ with parametrization $\vr(t):\mathbb{S}^1\rightarrow\mathbb{R}^n$, we have:
\begin{equation}\label{Yexp}
\langle\rho_1(\gamma)\rangle=-\frac{1}{32\pi^3}\int\displaylimits_{\Delta_3'}I_{\langle\rho_1\rangle}\Biggl(\vr(t_2)-\vr(t_3),\vr(t_1)-\vr(t_3),\vr(t_2)-\vr(t_1),\dot{\vr(t_1)},\dot{\vr(t_2)},\dot{\vr(t_3)}\Biggr)dt_1t_2dt_3
\end{equation}

At this point, the real difficulty lies in computing a closed form for the function $I_{\langle\rho_2\rangle}$. We do not expect this to be possible to do in general, however with a view towards deriving bounds for the configuration space integral appearing in \eqref{Yexp}, there are ways to simplify the Gaussian integral \eqref{YGaussian} significantly. To begin, we normalize $\va_1$ and $\va_2$ to obtain $\vq_1=\va_1/\|\va_1\|$ and $\vq_2=\va_2/\|\va_2\|$, and define a basis (similar to \eqref{specialbasis}) by: 
\begin{equation}\label{specialbasis2}
\vb_1=\frac{\vq_1+\vq_2}{\|\vq_1+\vq_2\|}\textrm{ , }
\vb_2=\frac{\vq_1-\vq_2}{\|\vq_1-\vq_2\|}\textrm{ and }\vb_3\textrm{,}\vb_4,...,\vb_n\in\mathrm{span}\{\vq_1,\vq_2'\}^\perp,
\end{equation}
With respect to the basis \eqref{specialbasis2}, the vectors $\vq_1$ and $\vq_2$ are represented as:
\begin{equation}\nonumber
\vq_1=a\ve_1+b\ve_2\textrm{ and }\vq_2=a\ve_1-b\ve_2,
\end{equation}
where $a=\cos(\phi/2)$ and $b=\sin(\phi/2)$ with $\phi$ the angle between $\vq_1$ and $\vq_2$. Additionally, due to the nature of the inputs to $I_{\langle\rho_2\rangle}$ in \eqref{Yexp} when computing $\langle\rho_2(\gamma)\rangle$,  we can also assume that $\va_1-\va_2=\va_3$, and derive the form of $\vq_3=\va_3/\|\va_3\|$ in the same basis:
\begin{equation}\nonumber
\vq_3=\frac{a(\|\va_1\|-\|\va_2\|)\ve_1+b(\|\va_1\|+\|\va_2\|)\ve_2}{\|\va_3\|}.
\end{equation}
In this basis, we then have simplified expressions for the quantities $\|U^T\va_1\|$, $\|U^T\va_2\|$, and $\|U^T\va_3\|$ appearing in the definitions of $C_1$, $C_2$, and $C_3$ along with the quantities in the determinants in \eqref{check1} and \eqref{check2}, and will only depend on the first two entries in each column of the Gaussian random matrix $U$.  Additionally, in \eqref{check1} we have:
\begin{flalign*}
&\mathrm{Det}\left[\tilde{\vy},\vx,\vy\right]\mathrm{Det}\left[\tilde{\vz},\tilde{\vx},\frac{\vx\|\vy\|+\vy\|\vx\|}{\|\vy\|}\right]+\mathrm{Det}\left[\tilde{\vz},\vx,\vy\right]\mathrm{Det}\left[\tilde{\vy},\tilde{\vx},\frac{\vx\|\vy\|+\vy\|\vx\|}{\|\vx\|}\right]&\\
&=\|\vx\|^2\|\vy\|\mathrm{Det}\left[\tilde{\vy},\frac{\vx}{\|\vx\|},\frac{\vy}{\|\vy\|}\right]\mathrm{Det}\left[\tilde{\vz},\tilde{\vx},\frac{\vx}{\|\vx\|}+\frac{\vy}{\|\vy\|}\right]+\|\vx\|\|\vy\|^2\mathrm{Det}\left[\tilde{\vz},\frac{\vx}{\|\vx\|},\frac{\vy}{\|\vy\|}\right]\mathrm{Det}\left[\tilde{\vy},\tilde{\vx},\frac{\vx}{\|\vx\|}+\frac{\vy}{\|\vy\|}\right]&,
\end{flalign*}
As such, the term $\frac{\vx}{\|\vx\|}+\frac{\vy}{\|\vy\|}$ will be proportional to the basis element $\vb_1$ when computing $I_{\langle\rho_2(\gamma)\rangle}$ in the basis \eqref{specialbasis2}, as will $\vx/\|\vx\|$ and $\vy/\|\vy\|$. On top of these simplifications, we may also make $\tilde{\vx}$ and $\tilde{\vy}$ (or $\tilde{\vx}$ and $\tilde{\vz}$) orthogonal without changing the values of the determinants, and then use the $O(n)$ invariance to align the vectors along the axes of the standard basis. Similar simplifications can be applied to \eqref{check2} in order to reduce the Gaussian integrals to simplified integrals over a reduced number of variables. All of the above techniques lead to the following identity:
\begin{equation}\nonumber
I_{\langle\rho_1\rangle}(\mathbf{x},\mathbf{y},\mathbf{z},\tilde{\mathbf{x}},\tilde{\mathbf{y}},\tilde{\mathbf{z}}) = -4\Big[(\mathbf{x}\cdot\tilde{\mathbf{x}})(\tilde{\mathbf{y}}\cdot\tilde{\mathbf{z}}) - (\mathbf{y}\cdot\tilde{\mathbf{y}})(\tilde{\mathbf{z}}\cdot\tilde{\mathbf{x}}) - (\mathbf{z}\cdot\tilde{\mathbf{z}})(\tilde{\mathbf{x}}\cdot\tilde{\mathbf{y}})\Big]\cdot K_1(\mathbf{z},-\mathbf{y})
\end{equation}
where

\begin{align*}
K_1(\mathbf{a},\mathbf{b}) &:= \frac{C}{\varepsilon\,\gamma_{ab}\,\delta} - \frac{\arctan(\varepsilon/C)}{\varepsilon^2} - \frac{\arctan(\gamma_{ab}/C)}{\gamma_{ab}^2} - \frac{\arctan(\delta/C)}{\delta^2}, \\
\quad \alpha^2 &:= \|\mathbf{a}\|^2,\quad \beta^2 := \|\mathbf{b}\|^2,\quad \gamma_{ab} := \mathbf{a}\cdot\mathbf{b},\quad C := \sqrt{\alpha^2\beta^2 - \gamma_{ab}^2},\quad \varepsilon := \alpha^2-\gamma_{ab},\quad \delta := \beta^2-\gamma_{ab}.
\end{align*}

\begin{Remark}
Notice that if instead of considering random projections of a fixed space curve to $3$ dimensional subspaces we considered projections to $2$ dimensional subspaces, then we would obtain a model of random plane curves, which coincides with the model we considered when computing the second moment of the winding number. In \cite{LW} it is shown that for immersions $\gamma:\mathbb{S}^1\rightarrow\mathbb{R}^2$ that $I_X(\gamma)=0$ and that $I_Y(\gamma)$ is an invariant of unicursal plane curves, which motivates finding quantities like $\langle I_Y(\gamma)\rangle$ and $\langle I_Y(\gamma)^2\rangle$. Notice that in the model considered in the previous section, that as $\epsilon$ approaches $0$, then the distribution of the permutations of the heights seems to approach the uniform distribution, and also the projections become closer to being planar. In this way, the distribution of $I_Y(\gamma)$ should approach the distribution of $v_2$. Another important point is that when we are projecting the curve in $\mathbb{R}^n$ to $\mathbb{R}^2$, then the terms \eqref{check1} and \eqref{check2} approach $0$ and only the term \eqref{check0} remains in the definition of the function $F(\vx, \vy, \vz, \tilde{\vx}, \tilde{\vy}, \tilde{\vz})$. 
\end{Remark} 
Beyond invariants of knots and links in $\mathbb{R}^3$, there are also configuration space integrals for computing linking numbers of links in hyperbolic space, for which a computation akin to the one in chapter 3 seems both interesting and feasible, provided one found a reasonable integro-geometric space like $Gr(n,3)$ to integrate over. 
\newpage
\section*{Acknowledgement}
This is an expanded version of my Ph.D. dissertation at the University of California, Davis. I would like to thank my adviser, Eric Babson, for his help, encouragement, and guidance over the years in the completion of this work.
\newpage
\appendix

\newchapter{Petal Diagrams (Python Implementation)}{Petal Diagrams (Python Implementation)}{Petal Diagrams (Python Implementation)}
In this appendix we provide both Python \cite{Rossum:1995:PRM:869369} and Wolfram Mathematica \cite{Mathematica} code to generate the initial space curves appearing in Chapter $5$, and the results in Tables \ref{tab:mysecondtable} and \ref{tab:myfirsttable}. The Python implementation is based largely off of NumPy methods \cite{numpy}, and the Wolfram Mathematica implementation closely follows the Python one, but allows for a parallelization of Monte Carlo trials using Mathematica's ParallelTable functionality.

\section{Python}
First we have the definitions of the two initial space curves in $\mathbb{R}^{2+2k}$ where $k$ is the number of petals per component:
\begin{lstlisting}
import numpy as np
def line(t,t1,t2,a,b):
    return (1/(t2-t1))*((t2-t)*a+(t-t1)*b)
def inputcurve1(x,n,eps):
#Petal diagram in first two coordinates:
    listing=[np.cos(n*x)*np.cos(x),np.cos(n*x)*np.sin(x)]
    for i in range(0,n):
            listing.append(np.piecewise(x, [x>=(np.pi/(2*n))*(2*i) and 
            x < (np.pi/(2*n))*(2*i+1), x>= (np.pi/(2*n))*(2*i+1) and 
            x < (np.pi/(2*n))*(2*i+2)], [lambda t: 
	    line(t,(np.pi/(2*n))*(2*i), (np.pi/(2*n))*(2*i+1),0,eps), 
	    lambda t: line(t,(np.pi/(2*n))*(2*i+1), 
	    (np.pi/(2*n))*(2*i+2),eps,0)]))
    for i in range(0,n):
        listing.append(0)
    return np.array(listing)
def inputcurve2(x,n,eps):
#Rotated petal diagram in first two coordinates:
    listing=[(np.cos((n-2)*np.pi/n))*np.cos(n*x)*np.cos(x)-
             (np.sin((n-2)*np.pi/n))*np.cos(n*x)*np.sin(x),
             (np.sin((n-2)*np.pi/n))*np.cos(n*x)*np.cos(x)+
             (np.cos((n-2)*np.pi/n))*np.cos(n*x)*np.sin(x)]
    for i in range(0,n):
        listing.append(0)
    for i in range(0,n):
            listing.append(np.piecewise(x, [x>=(np.pi/(2*n))*(2*i) 
            and x < (np.pi/(2*n))*(2*i+1), x>= (np.pi/(2*n))*(2*i+1) 
            and x < (np.pi/(2*n))*(2*i+2)], [lambda t: 
            line(t,(np.pi/(2*n))*(2*i), (np.pi/(2*n))*(2*i+1),0,eps), 
            lambda t: line(t,(np.pi/(2*n))*(2*i+1),
	  (np.pi/(2*n))*(2*i+2),eps,0)]))
    return np.array(listing)  
\end{lstlisting}
Once the initial space curves are defined, then we can sample numerous random links by sampling matrices $U\in\mathrm{Mat}_{n,3}$ with $N(0,1)$ entries. In what follows, we sample the $n$ Gaussian matrices, $U$, make the columns orthonormal by applying the Gram-Schmidt algorithm to columns to generate the new matrix $Q$, and then we project the initial space curves on to the span of the columns of $Q$ to get a link in $\mathbb{R}^3$. We then use the algorithm in \cite{Arai} to compute the linking number. The second moment of the linking number for these $n$ samples is returned.
\begin{lstlisting}
#Algorithm for computing the linking number, from reference [3]
def linking(v1,v2):
    aa=len(v1);
    bb=len(v2);
    II=0;
    def solidangle(a,b,c):
                return 2*math.atan2(np.linalg.det([a,b,c]),
                        np.linalg.norm(a)*np.linalg.norm(b)*
                        np.linalg.norm(c)+ np.dot(a,b)*np.linalg.norm(c)+
                        np.dot(c,a)*np.linalg.norm(b)+
                        np.dot(b,c)*np.linalg.norm(a))
    for i in range(0,aa):
        for j in range(0,bb):
            a=np.array(v2[j])-np.array(v1[i])
            b=np.array(v2[j])-np.array(v1[(i+1)%aa])
            c=np.array(v2[(j+1)%bb])-np.array(v1[(i+1)%aa])
            d=np.array(v2[(j+1)%bb])-np.array(v1[i])
            II=II+solidangle(a,b,c)+solidangle(c,d,a)
    return  II/(4*np.pi)
#Gram-Schmidt for triples of vectors
def GS(U):
    q1=U[0]/np.linalg.norm(U[0])
    q2=(U[1]-np.dot(U[1],q1)*q1)/np.linalg.norm(U[1]-np.dot(U[1],q1)*q1)
    q3=(U[2]-(np.dot(U[2],q1))*q1-(np.dot(U[2],q2))*q2)/np.linalg.norm(U[2]
        -(np.dot(U[2],q1))*q1-(np.dot(U[2],q2))*q2)
    return [q1,q2,q3]
#In what follows, n is the number of samples, k is the number of petals
def samples(n,k):
    LinkingNumberList=[]
    for i in range(0,n):
        UU=np.transpose(np.random.normal(0,1,(2*k+2,3)))
        Q=GS(UU)
        #here we set epsilon=1
        def R11(t):
            return np.dot(Q,inputcurve1(t,k,1))
        def R22(t):
            return np.dot(Q,inputcurve2(t,k,1))
        vectorlist1=[]
        vectorlist2=[]
        #Samples 20 points from each component of link
        #to feed in to the linking function above
        for i in range(0,20):
            vectorlist1.append(R11((np.pi/20)*i))
            vectorlist2.append(R22((np.pi/20)*i))
        LinkingNumberList.append(linking(vectorlist1,vectorlist2))
    return (np.average(LinkingNumberList)**2+np.std(LinkingNumberList)**2)
\end{lstlisting}
As an example calculation, we ran $\mathrm{samples}(10000,3)$ (that is, we sampled $10000$ orthonormal frames and projected the initial space curves corresponding to $k=3$ petals per component) to get that $\langle\mathrm{Lk}^2\rangle\approx .7181$. Now compare it to numerically integrating the result in \ref{mainResultTheorem}:
\begin{lstlisting}
#Numerical derivatives, n is number of components
#of vector function being differentiated:
def derivativen(f,n):
    def df(x, h=0.1e-8):
        lists=[]
        for i in range(0,n):
            lists.append(( f(x+h/2)[i] - f(x-h/2)[i] )/h)
        return lists
    return df
def r1(t):
    return inputcurve1(t,3,1)
def r2(t):
    return inputcurve2(t,3,1)
dr1 = derivativen(r1,8)
dr2 = derivativen(r2,8)
def r3(s,t):
    return np.subtract(r2(t),r1(s))
def integrandLK2(a1,a2,a3,a11,a22,a33):
    return (2.0/np.pi)*(((np.dot(a3,a3))*(np.dot(a33,a33))-
             np.dot(a3,a33)**2.0)*np.linalg.det(np.dot([a1,a2,a3],
             np.transpose([a11,a22,a33])))+np.dot(a3,a33)*
             np.linalg.det(np.dot([a3,a11,a22,a33],
             np.transpose([a33,a1,a2,a3]))))/(((np.dot(a3,a3))*
             (np.dot(a33,a33)))*(((np.dot(a3,a3))*
             (np.dot(a33,a33))-np.dot(a3,a33)**2.0))**(1.5))
def inta(s,t,ss,tt):
    return (1/(16.0*np.pi**2.0))*integrandLK2(dr1(s),dr2(t),r3(s,t),
                                              dr1(ss),dr2(tt),r3(ss,tt))
#The standard MonteCarlo integration, here n
#is the number of samples from the integration domain
def mcintegrate(n):
    valuelist=[]
    for i in range(0,n):
        s=random.uniform(0,np.pi)
        t=random.uniform(0,np.pi)
        ss=random.uniform(0,np.pi)
        tt=random.uniform(0,np.pi)
        valuelist.append(inta(s,t,ss,tt))
#Returns approximation of integral along with an error estimate
    return [(np.pi**4)*np.average(valuelist),
            (np.pi**4)*(((np.dot(valuelist,valuelist)/n)-
            (np.average(valuelist))**2)/n)**(0.5)]
\end{lstlisting}
Running $\mathrm{mcintegrate(500000)}$ here we obtain $\langle\mathrm{Lk}^2\rangle\approx .7525$ with an error estimate of $.017$

\section{\textit{Wolfram Mathematica} }
\lstset{language=Mathematica}
\lstset{basicstyle={\sffamily\footnotesize},
  breaklines=true,
  captionpos={b},
  frame={lines},
  rulecolor=\color{black},
  framerule=0.5pt,
  columns=flexible,
  tabsize=2
}
 \begin{lstlisting}[language=Mathematica,caption={Definitions of components of Petaluma links}]
k=3; (*Number of petals*)

(*linear interpolation*)
line[t_,t1_,t2_,a_,b_]:=(1.0/(t2-t1))*((t2-t)*a+(t-t1)*b);

(*define indicator*)
indicator[x_,{a_,b_}]:=If[(x>=a)&&(x<b),1,0];

(*define strand intervals*)
intervals=Flatten[Reap[Do[Sow[{{(Pi*(2*i))/(2*k),(Pi*(2*i+1))/(2*k)},{(Pi*(2*i+1))/(2*k),(Pi*(2*i+2))/(2*k)}}],{i,0,k-1}]][[2]][[1]],1];

(*Definition of r1(t) and r2(t) with smooth base petal diagram*)
input1[t_,n_,eps_]:=PadRight[{Cos[n*t]*Cos[t],Cos[n*t]*Sin[t]},2*(n+1)]+Sum[(indicator[t,intervals[[2*i-1]]]*line[t,intervals[[2*i-1]][[1]],intervals[[2*i-1]][[2]],0,eps]+indicator[t,intervals[[2*i]]]*line[t,intervals[[2*i]][[1]],intervals[[2*i]][[2]],eps,0])*UnitVector[2*(n+1),i+2],{i,1,n}];

input2[s_,n_,eps_]:=PadRight[{Cos[(n-2)*Pi/n]*Cos[n*s]*Cos[s]-Sin[(n-2)*Pi/n]*Cos[n*s]*Sin[s],Sin[(n-2)*Pi/n]*Cos[n*s]*Cos[s]+Cos[(n-2)*Pi/n]*Cos[n*s]*Sin[s]},2*(n+1)]+Sum[(indicator[s,intervals[[2*i-1]]]*line[s,intervals[[2*i-1]][[1]],intervals[[2*i-1]][[2]],0,eps]+indicator[s,intervals[[2*i]]]*line[s,intervals[[2*i]][[1]],intervals[[2*i]][[2]],eps,0])*UnitVector[2*(n+1),i+2+n],{i,1,n}];

(*Sample points from r1(t), r2(t) to create PL Petaluma diagrams*)
vertices1=Reap[Do[If[Mod[nn,4]!=0,Sow[input1[(nn*Pi)/(4*k),k,1][[1;;2]]]],{nn,0,4*k-1}]][[2]][[1]];
vertices2=Reap[Do[If[Mod[nn,4]!=0,Sow[input2[(nn*Pi)/(4*k),k,1][[1;;2]]]],{nn,0,4*k-1}]][[2]][[1]];

(*Used for reparametrizing segments of input space curves so that
each segment is traversed in equal time*)
intervalPolygonal=Reap[Do[Sow[{(Pi*(i))/(3*k),(Pi*(i+1))/(3*k)}],{i,0,3*k-1}]][[2]][[1]];

(*Finally, define input space curves*)
curve1[t0_,n0_,eps0_]:=Module[{t=t0,n=n0,eps=eps0},
inputp1[t_]:=Sum[indicator[t,intervalPolygonal[[i+1]]]*line[t,intervalPolygonal[[i+1]][[1]],intervalPolygonal[[i+1]][[2]],vertices1[[Mod[i,3*n]+1]],vertices1[[Mod[i+1,3*n]+1]]],{i,0,3*n-1}];

inputPolygonal1[t_]:=PadRight[{inputp1[Mod[t-Pi/(2*3*n),Pi]][[1]],inputp1[Mod[t-Pi/(2*3*n),Pi]][[2]]},2*(n+1)]+Sum[(indicator[t,intervals[[2*i-1]]]*line[t,intervals[[2*i-1]][[1]],intervals[[2*i-1]][[2]],0,eps]+indicator[t,intervals[[2*i]]]*line[t,intervals[[2*i]][[1]],intervals[[2*i]][[2]],eps,0])*UnitVector[2*(n+1),i+2],{i,1,n}];
inputPolygonal1[t]
];

curve2[t0_,n0_,eps0_]:=Module[{t=t0,n=n0,eps=eps0},
inputp2[t_]:=Sum[indicator[t,intervalPolygonal[[i+1]]]*line[t,intervalPolygonal[[i+1]][[1]],intervalPolygonal[[i+1]][[2]],vertices2[[Mod[i,3*n]+1]],vertices2[[Mod[i+1,3*n]+1]]],{i,0,3*n-1}];

inputPolygonal2[t_]:=PadRight[{inputp2[Mod[t-Pi/(2*3*n),Pi]][[1]],inputp2[Mod[t-Pi/(2*3*n),Pi]][[2]]},2*(n+1)]+Sum[(indicator[t,intervals[[2*i-1]]]*line[t,intervals[[2*i-1]][[1]],intervals[[2*i-1]][[2]],0,eps]+indicator[t,intervals[[2*i]]]*line[t,intervals[[2*i]][[1]],intervals[[2*i]][[2]],eps,0])*UnitVector[2*(n+1),i+2+n],{i,1,n}];
inputPolygonal2[t]
];

(*Compute input space curve derivatives*)
dcurve1[t_,n_,eps_]:=D[curve1[tt,n,eps],tt]/.tt->t;
dcurve2[t_,n_,eps_]:=D[curve2[tt,n,eps],tt]/.tt->t;
  \end{lstlisting}

With the input space curves defined, we then create our various Monte Carlo sampling processes. Note that we allow each function to have a seed set for the random generation so that results may be repeated:
\vspace{.2in}
 \begin{lstlisting}[language=Mathematica,caption={Definitions of Monte Carlo Sampling Methods}]
(*First define the integrand*)
integrandLK2[a1_,a2_,a3_,a11_,a22_,a33_]:=Quiet[1/(4*Pi)^2*4/Pi*(a3.a33-Sqrt[(a3.a3)*(a33.a33)-(a3.a33)^2]*ArcSin[a3.a33/(Norm[a3]*Norm[a33])])/((a3.a33)^3*Sqrt[(a3.a3)*(a33.a33)-(a3.a33)^2])*Det[{a1,a2,a3}.Transpose[{a11,a22,a33}]]]

(*Next define the basic sampling*)
mcSampleGen[dcurve1_,dcurve2_,curve1_,curve2_,eps_,n_]:=Reap[Do[SeedRandom[i];
uu=RandomVariate[UniformDistribution[{0,Pi}],4];
Sow[integrandLK2[dcurve1[uu[[1]],eps],dcurve2[uu[[2]],eps],curve2[uu[[2]],eps]-curve1[uu[[1]],eps],dcurve1[uu[[3]],eps],dcurve2[uu[[4]],eps],curve2[uu[[4]],eps]-curve1[uu[[3]],eps]]],{i,1,n}]][[2]][[1]]

(*Parallelized version*)
mcSampleGenPARALLEL[dcurve10_,dcurve20_,curve10_,curve20_,eps0_,n0_]:=
Module[{dcurve1=dcurve10,dcurve2=dcurve20,curve1=curve10,curve2=curve20,eps=eps0,n=n0},
batches=QuotientRemainder[n,4];
table=ParallelTable[mcSampleGen[dcurve1,dcurve2,curve1,curve2,eps,If[i!=4,batches[[1]],batches[[1]]+batches[[2]]],i],{i,1,4}];
table=Flatten[table];
cleanSamples=DeleteCases[table,Indeterminate];
cleanSamples
]

(*Monte Carlo integration is finally computed with the following lines. First obtain MC samples of the integrand:*)
samples = mcSampleGenPARALLEL[dcurve1, dcurve2, curve1, curve2, 1, 1000000];

(*Remove samples resulting in indeterminate values, these arise at singularities of the denominator which ultimately end up being controlled by the numerator*)
cleanSamples = DeleteCases[samples, Indeterminate];

Pi^4*{Mean[cleanSamples], Sqrt[(Mean[cleanSamples^2] - Mean[cleanSamples]^2)/ Length[cleanSamples]]}
\end{lstlisting}

The result of this numerical integration scheme, with $10^6$ Monte Carlo samples is $0.83\pm0.026$. We compare this to the alternative approach given by sampling Gaussian random matrices $U\in\mathrm{Mat}_{n,3}$ and computing the linking number using the algorithm from \cite{Arai}. 
\vspace{.2in}
 \begin{lstlisting}[language=Mathematica,caption={Alternate Definitions of Monte Carlo Sampling Methods}]
(*First define the basic sample generator which will be parallelized later*)

gaussianSampleGen[curve1_, curve2_, n_, nsamples_, eps0_, seed0_] := 
 Module[{curve1test = curve1, curve2test = curve2, eps = eps0, 
   seed = seed0},

(*Compute solid angle subtended by set of 3 vectors*)
  solidangle[a_, b_, c_] := 
   2*ArcTan[
     Norm[a]*Norm[b]*Norm[c] + (a.b)*Norm[c] + (c.a)*Norm[b] + (b.c)*
       Norm[a], Det[{a, b, c}]];
  
(*Obtain points for PL approximation of knot*)

  knotpoints = 
   Union[Reap[Do[Sow[(Pi*(2*i + 1))/(2*3*n)], {i, 0, 3*n - 1}]][[2]][[
     1]], Reap[Do[Sow[(i*Pi)/(2*n)], {i, 0, 2*n - 1}]][[2]][[1]]];
  vertices1h = 
   Reap[Do[Sow[curve1test[knotpoints[[nn]], eps]], {nn, 1, 
        Length[knotpoints]}]][[2]][[1]];
  vertices2h = 
   Reap[Do[Sow[curve2test[knotpoints[[nn]], eps]], {nn, 1, 
        Length[knotpoints]}]][[2]][[1]];
  
  MM = Length[vertices1h];
  NN = Length[vertices2h];
  
  SeedRandom[seed]; (*Random Seed set here*)
  
(* Return list of linking numbers*)
  samples = Reap[Do[
       U = RandomVariate[NormalDistribution[], {3, 2*(n + 1)}];
       klist = Transpose[U.Transpose[vertices1h]];
       llist = Transpose[U.Transpose[vertices2h]];
       linking = 1/(4*Pi)*Total[Reap[Do[Sow[Total[Reap[Do[
                    
                    alpha = 
                    llist[[Mod[jj, NN] + 1]] - 
                    klist[[Mod[ii, MM] + 1]];
                    
                    beta = 
                    llist[[Mod[jj, NN] + 1]] - 
                    klist[[Mod[ii + 1, MM] + 1]];
                    
                    gamma = 
                    llist[[Mod[jj + 1, NN] + 1]] - 
                    klist[[Mod[ii + 1, MM] + 1]];
                    
                    delta = 
                    llist[[Mod[jj + 1, NN] + 1]] - 
                    klist[[Mod[ii, MM] + 1]];
                    
                    Sow[solidangle[alpha, beta, gamma] + 
                    solidangle[gamma, delta, alpha]], {jj, 0, 
                    NN - 1}]][[2]][[1]]]], {ii, 0, MM - 1}]][[2]][[
           1]]];
       Sow[linking];, {sample, 1, nsamples}]][[2]][[1]]]

(*Create parallelized version of above*)

gaussianSampleGenPARALLEL[curve10_, curve20_, n0_, nsamples0_, eps0_] := 
 Module[{curve1 = curve10, curve2 = curve20, n = n0, nsamples = nsamples0, eps = eps0},
  batches = QuotientRemainder[nsamples, 4];
  table = 
   ParallelTable[
    gaussianSampleGen[curve1, curve2, n, 
     If[i != 4, batches[[1]], batches[[1]] + batches[[2]]], eps, 
     i], {i, 1, 4}];
  Flatten[table]
  ]
\end{lstlisting}

Finally we obtain our samples of the linking number and compute the second moments of the linking number:

 \begin{lstlisting}[language=Mathematica,caption={Alternate Definitions of Monte Carlo Sampling Methods}]
(*Obtain samples*)
gaussiansamples = gaussianSampleGenPARALLEL[curve1, curve2, 3, 1000000, 1];

(*Compute second moment*)
StandardDeviation[gaussiansamples]^2 + Mean[gaussiansamples]^2
\end{lstlisting}

The result is $0.789$, and is close to the value we obtained by numerical integration over the relevant configuration space.

\newpage
\newchapter{\textit{Wolfram Mathematica} Numerical Implementation of the Degree $2$ Vassiliev Invariant}{\textit{Wolfram Mathematica} Numerical Implementation of the Degree $2$ Vassiliev Invariant}{\textit{Wolfram Mathematica} Numerical Implementation of the Degree $2$ Vassiliev Invariant}

Begin by defining the function $F(\vx,\vy,\vz,\vx',\vy',\vz')$ appearing \cite{MCVassiliev}:
\lstset{language=Mathematica}
\lstset{basicstyle={\sffamily\footnotesize},
  breaklines=true,
  captionpos={b},
  frame={lines},
  rulecolor=\color{black},
  framerule=0.5pt,
  columns=flexible,
  tabsize=2
}
  \begin{lstlisting}[language=Mathematica,caption={$\rho_2(\gamma)$ integrand}]
F[x0_,y0_,z0_,xx0_,yy0_,zz0_]:=
Module[{x=x0,y=y0,z=z0,xx=xx0,yy=yy0,zz=zz0},
Av=x;
Bv=y;
Cv=z;
a=Norm[Av];
b=Norm[Bv];
c=Norm[Cv];
C1=(2*Pi)/(a*b*c);
C2=1/(a*b+Av.Bv);
C3=a+b-c;
t1=(yy.zz)*(xx.Cv)+(xx.zz)*(yy.Bv)-(xx.yy)*(zz.Av);
t2=Det[{yy,Av,Bv}]*Det[{zz,xx,(Av*b+Bv*a)/b}]+
Det[{zz,Av,Bv}]*Det[{yy,xx,(Bv*a+Av*b)/a}];

t3=Det[{yy,Av,Bv}]*Det[{zz,xx,Bv*(Norm[Cv]-Norm[Av])/Norm[Bv]^2+(Cv)*(Norm[Av]+Norm[Bv])/Norm[Cv]^2}]+
Det[{zz,Av,Bv}]*Det[{yy,xx,Av*(Norm[Cv]-Norm[Bv])/Norm[Av]^2-(Cv)*(Norm[Av]+Norm[Bv])/Norm[Cv]^2}];
(*Returns:*)
(-1/(32*Pi^3))*(C1*C2*C3*t1-C1*C2^2*C3*t2+C1*C2*t3)
]
  \end{lstlisting}

Now with the integrand appearing in the time-ordered integral $\rho_2(\gamma)$ defined, we can compute the value of $v_2$ for specific examples of parametrized knots. For example, if we take a trefoil knot parametrized by:
\begin{lstlisting}[language=Mathematica,caption={Trefoil parametrization}]
trefoil[t_] := {Sin[t] + 2*Sin[2*t], Cos[t] - 2*Cos[2*t], -Sin[3*t]};
x[s_] := trefoil[s];
y[t_] := trefoil[t];
z[u_] := trefoil[u];
w[v_] := trefoil[v];
xx[s_] := D[x[ss], ss] /. ss -> s;
yy[t_] := D[y[tt], tt] /. tt -> t;
zz[u_] := D[z[uu], uu] /. uu -> u;
ww[v_] := D[w[vv], vv] /. vv -> v;
\end{lstlisting}
then we can compute the integrals $\rho_1(\gamma)$ and $\rho_2(\gamma)$, and compare their sum to the value of $v_2$ on the trefoil which is given by $v_2(\mathrm{trefoil})=\frac{23}{12}$:

\begin{lstlisting}[language=Mathematica,caption={Time-ordered numerical integrals},mathescape]
R = ImplicitRegion[(x < y) && (y < z), {{x, 0, 2*Pi}, {y, 0, 2*Pi}, {z, 0, 2*Pi}}];

I1=NIntegrate[F[y[tt] - x[ss], z[uu] - x[ss], y[tt] - z[uu], xx[ss],  yy[tt], zz[uu]], 
Element[{uu, tt, ss},R], Method -> {"MonteCarlo", MaxPoints -> 10000000}];

Rr = ImplicitRegion[(t1 < t2) && (t2 < t3) && (t3 < t4), {{t1, 0, 
     2*Pi}, {t2, 0, 2*Pi}, {t3, 0, 2*Pi}, {t4, 0, 2*Pi}}];

I2=NIntegrate[1/(8*Pi^2)*Det[{xx[t1], zz[t3], (z[t3] - x[t1])/Norm[z[t3] - x[t1]]^3}]*
  Det[{yy[t2], ww[t4], (w[t4] - y[t2])/Norm[w[t4] - y[t2]]^3}], 
Element[{t1, t2, t3, t4}, Rr], 
 Method -> {"MonteCarlo", MaxPoints -> 10000000}];

I1+I2
\end{lstlisting}

In the above we obtained $I_1=-0.140208$ and $I_2=2.05305$, resulting in $v_2=1.91222$, which is close to the value of $23/12=1.91\overline{6}$. Similar experiments with the unknot (where $v_2=-\frac{1}{12}$) yield equally close results.

We note that the domain of integration is complicated, yet still relatively low dimensional which opens the possibility to use Quasi Monte Carlo sampling and low discrepancy sequences to speed up the rate of convergence of the above integrals - we leave this approach to a future work.

\newpage
\bibliographystyle{plain}
\bibliography{RandomKnotsBib4} 
\end{document}